 \documentclass{article}

\usepackage{amsmath,epsfig}
\usepackage{graphics,amsmath,amsfonts,amssymb,epsfig,latexsym,amsfonts,graphicx,amssymb}
\usepackage{subfigure}   
\usepackage{cite}
\usepackage{bm}
\usepackage{makecell}
\usepackage{color}
\usepackage{diagbox} 
\usepackage{stfloats}  
\usepackage{footnote}
\usepackage{algorithm} 
\usepackage{algorithmic} 
\usepackage{multirow}  
\usepackage{hyperref}  

\usepackage{hyperref}
\usepackage{url}
\usepackage{amsmath}
\usepackage{bm}
\usepackage{amsfonts}
\usepackage{algorithm}
\usepackage{algorithmic}
\usepackage{enumitem,graphicx}
\usepackage{color,cite}

\newtheorem{lemma}{Lemma}
\newtheorem{prop}{Proposition}
\newtheorem{thm}{Theorem}

\newtheorem{coro}{Corollary}[section]

\newtheorem{claim}{Claim}[section]
\newtheorem{conjecture}{Conjecture}[section]

\newtheorem{rmk}{Remark}[section]

\newcommand{\st}{{\rm s.t.}}

\newcommand{\rmax}{\rm max}

\newcommand{\black}{ \color{black}}
\newcommand{\dR}{\mathbb{R}}
\newcommand{\dC}{\mathbb{C}}

\title{ \bf On the Efficiency of Random Permutation for ADMM and Coordinate Descent
}
\author{
	Ruoyu Sun $\quad \quad  $  Zhi-Quan Luo $\quad \quad  $ Yinyu Ye
	\thanks{ This paper is a strengthened
		version of a previous technical report ``On the expected convergence of randomly permuted ADMM'' appeared on arxiv on April 2015, with several new results, mainly the ones on the expected convergence rate of RP-CD and RP-ADMM. }	
	\thanks{Ruoyu Sun is with  University of Illinois at Urbana-Champaign, USA. Part of the work was done when this author was a student at University of Minnesota and a postdoc at Stanford University. Email: ruoyus@illinois.edu.}
	\thanks{Zhi-Quan Luo is with the Chinese University of Hong Kong, Shenzhen, China. He is also affiliated with
	University of Minnesota, Minneapolis, MN 55455, USA.
		Email: luozq@cuhk.edu.cn. }
	\thanks{Yinyu Ye is with the Department of Management Science and Engineering, School of Engineering, Stanford University, USA; and International Center of Management Science and Engineering, School of Management and Engineering, Nanjing
		University, China. Email: yyye@stanford.edu. }
}
\date{Dec 31, 2018}

\begin{document}

\maketitle

%
%

\begin{abstract}
Random permutation is observed to be powerful for optimization algorithms: for multi-block ADMM (alternating direction method of multipliers), while the classical cyclic version divergence,  the randomly permuted version converges in practice; for BCD (block coordinate descent), the randomly permuted version is typically faster than other versions. In this paper, we provide strong theoretical evidence that random permutation has positive effects on ADMM and BCD, by analyzing randomly permuted ADMM (RP-ADMM) for solving linear systems of equations, and randomly permuted BCD (RP-BCD) for solving unconstrained quadratic problems. First, we prove that RP-ADMM converges in expectation for solving systems of linear equations. The key technical result is that the spectrum of the expected update matrix of RP-BCD lies in $(-1/3, 1)$, instead of the typical range $(-1, 1)$. Second, we establish expected convergence rates of RP-ADMM for solving linear sytems and RP-BCD for solving unconstrained quadratic problems. This expected rate of RP-BCD is $O(n)$ times better than the worst-case rate of cyclic BCD, thus establishing a gap of at least $O(n)$ between RP-BCD and cyclic BCD. To analyze RP-BCD, we propose a conjecture of a new matrix AM-GM (algebraic mean-geometric mean) inequality, and prove a weaker version of it. 

\end{abstract}

 \thispagestyle{empty}
\newpage
\tableofcontents
\newpage


\section{Introduction}
{\black
	A simple yet powerful idea for solving large-scale computational problems is to iteratively solve smaller subproblems. The applications of this idea include coordinate descent (CD), POCS (Projection onto Convex Sets), SGD (Stochastic Gradient Descent). 
They are well suited for large-scale \emph{unconstrained} optimization problem (see, e.g. Wright \cite{wright2015coordinate}, for a recent survey of CD) since it decomposes a large problem into small subproblems.
 The decomposition idea is crucial for huge problems due to both the cheap per-iteration cost and small memory requirement. Moreover, this idea is ``orthogonal'' to other large-scale optimization ideas such as first-order methods (using only gradient information) and random projection, and thus can be easily combined with other ideas.  
 
This paper is motivated by a natural question: how should we extend the decomposition idea to solve problems \emph{with constraints}?
We consider a constrained minimization problem with a convex objective function and linear constraints (this is for motivation; our analysis is for a much simpler version):
\begin{equation}\label{optimization, sum fi}
\begin{split}
   \min_{x_1, \dots, x_n} \quad    &  f(x_1, x_2, \dots, x_n), \\
  \st \quad   &   A_1 x_1 + \dots + A_n x_n = b,  \\
      \quad   &   x_i \in \mathcal{X}_i, \ i=1,\dots, n,
\end{split}
\end{equation}
where $A_i \in \dR^{ N \times d_i }, b \in \dR^{ N \times 1 }, \mathcal{X}_i \subseteq \dR^{d_i} $ is a closed convex set,
$i=1,\dots, n$, and $f: \dR^{d_1 + d_2 + \dots + d_n} \rightarrow \dR $ is a closed convex function.
Many  machine learning and engineering problems can be cast into linearly-constrained optimization problems
with two blocks (see Boyd et al. \cite{boyd2011distributed} for many examples) or more than two blocks (e.g. linear programming, robust principal component analysis, composite regularizers for structured sparsity; see Chen et al. \cite{chen2014direct} and Wang et al. \cite{wang2014parallel} for more examples).

To apply the decomposition idea to a constrained problem, one possible way is to form the augmented Lagrangian function and perform coordinate descent for the primal problem and a gradient step for the dual problem, i.e. combining BCD with augmented Lagrangian method, to obtain the so-called alternating direction method of multipliers (ADMM). ADMM  was originally proposed in
Glowinski and Marroco \cite{glowinski1975approximation} (see also Chan and Glowinski \cite{chan1978finite}, Gabay and Mercier \cite{gabay1976dual}) to solve problem \eqref{optimization, sum fi} when there are only two blocks (i.e. $n=2$) and the objective function is separable. It is natural and computationally beneficial to extend the original ADMM directly to solve the general $n$-block problem \eqref{optimization, sum fi} via the following procedure:
\begin{equation}\label{cyclic ADMM}
\begin{cases}
x_{1}^{k+1} = \arg \min_{ x_1 \in \mathcal{X}_{1} }  \mathcal{L}( x_{1},  x_{2}^{k}, \dots, x_n^k ; \mu^k ),   \\
\quad \quad \quad  \vdots   \\
x_{n}^{k+1} = \arg \min_{ x_n \in \mathcal{X}_n } \mathcal{L}( x_1^{k+1},  \dots , x_{n-1}^{k+1}, x_n ; \mu^k ),  \\
\mu^{k+1} = \mu^k - \beta(  A_1 x_1^{k+1} + \dots +  A_n x_n^{k+1} - b ),
\end{cases}
\end{equation}
where the augmented Lagrangian function
\begin{equation}\label{augmented Lag func}
\mathcal{L}(x_1, \dots, x_n ; \mu ) = f(x_1, \dots, x_n) - \mu^T (\sum_i A_i x_i  - b) + \frac{\beta}{2} \| \sum_i A_i x_i - b \|^2.
\end{equation}



The convergence of the direct extension of ADMM to multi-block case had been an open question, until a counter-example was recently given in Chen et al. \cite{chen2014direct}. More specifically, Chen et al. \cite{chen2014direct} showed that even for the simplest scenario where the objective function is $0$ and the number of blocks is $3$, ADMM can be divergent for a certain choice of $A = [A_1, A_2, A_3]$.
There are several proposals to overcome the drawback (see, e.g., \cite{he2012alternating,he2012convergence,
hong2012linear, han2012note,chen2013convergence,he2013proximal,he2013full,deng2013parallel, lin2014convergence,hong2014block,cai2014direct,sun2014convergent,lin2014global,han2014augmented,li2014schur,li2015convergent,lin2015iteration,deng2017parallel}), but they either need to restrict the range of original problems being solved, add additional cost in each step of computation, or limit the stepsize in updating the Lagrange multipliers. These solutions typically slow down the performance of ADMM for solving most practical problems.
Moreover, it is not clear how to compare the convergence speed of these algorithms as they typically contain different parameters.
 One may ask whether a ``minimal'' modification of cyclic multi-block ADMM \eqref{cyclic ADMM} can lead to convergence,
 and whether we can provide some convergence speed analysis that is easy to interpret. 

One of the simplest modifications of \eqref{cyclic ADMM} is to add randomness to the update order.
Randomness  has been very useful in the analysis of block coordinate descent (BCD) methods
and stochastic gradient descent (SGD) methods.
In particular, a recent work Sun and Ye \cite{sun2016worst} showed that  randomized CD (R-CD) can be up to $O(n^2)$ times faster than
cyclic CD (C-CD) for quadratic minimization in the worst case,
where $n$ is the number of variables \footnote{Rigorously speaking, these two bounds are not directly comparable since the result for the randomized version only holds with high probability, while the result for the cyclic version always holds; anyhow, this $O(n^2)$ gap is still meaningful if ignoring this difference between deterministic and randomized algorithm.}.
 Another example is the comparison of IAG (Incremental Aggregated Gradient) in Blatt et al. \cite{blatt2007convergent} and its randomized version SAG (Stochastic Average Gradient) \cite{schmidt2013minimizing}: it turns out that the introduction of randomness leads to better iteration complexity bounds. There is also some study on randomly permuted version of pure SGD
 \cite{gurbuzbalaban2015random}. 
These examples show that randomization may improve the algorithm in theory and in practice. 

 It is important to note that the iteration complexity bounds for randomized algorithms are usually established for independent randomization (sampling with replacement),  while in practice, random permutation (sampling without replacement) has been reported to exhibit faster convergence (e.g. Shalev et al. \cite{shalev2013stochastic}, Recht and Re \cite{recht2013parallel}, Sun \cite{sun15thesis}). Interestingly, our simulation shows that for solving linear system of equations, randomly permuted ADMM (RP-ADMM) always converges, but independently randomized versions of ADMM  can be divergent even for Gaussian data. 
  Therefore, we focus on the analysis of RP-ADMM in this paper.
 
Random permutation is known to be notoriously difficult to analyze. Even for unconstrained quadratic minimization, the convergence rate of RP-BCD is poorly understood. 
Many existing works treated cyclic BCD and RP-BCD together \cite{tseng2001convergence,beck2013convergence,sun2015improved}, and thus
the best known convergence rate of RP-BCD for general convex problems are in fact the same as that of C-BCD \cite{sun2015improved}. However, in light of a recent study which established an up to $O(n^2)$ gap between cyclic CD and R-CD \cite{sun2016worst}, it is unlikely that RP-CD has the same convergence rate as C-CD since that would imply RP-CD could be $O(n^2)$-times slower than R-CD. 
For the special example that demonstrates the gap between C-CD and R-CD,  it was shown recently that RP-CD is faster than R-CD \footnote{This paper appeared after the first version of the current paper.} \cite{wright2017analyzing}. 
However, the general quadratic case seems to be quite difficult, probably due to its close connection to a matrix AM-GM (algebraic mean-geometric mean) inequality \cite{recht2012beneath}, the difficulty of which is essentially to prove an  inequality in non-commutative algebra.  


\subsection{Summary of Contributions}


We consider two extremes of a general RP-ADMM: i) the objective is zero, i.e., RP-ADMM for solving a linear system; ii) the constraint is zero and the objective is a quadratic function, i.e., RP-BCD for solving quadratic minimization.
Due to the lack of understanding of random permutation for quadratic minimization as discussed previously, we restrict to the two cases in this paper. 

The first result of this paper is the expected convergence of RP-ADMM for solving linear systems. More specifically, when the objective function is zero and the constraint is a non-singular square linear system of equations, the expected output of randomly permuted ADMM converges to the unique primal-dual optimal solution.
A major technical result in this proof is that the eigenvalues of the expected iteration matrix of RP-BCD for quadratic problems lie in $(-1/3, 1)$, instead of the typical range $(-1,1)$.

The second result is about the expected convergence rate of RP-ADMM for solving linear systems and RP-BCD for solving quadratic problems. 
We show that RP-BCD for a convex quadratic minimization problem with equal diagonal entries has expected iteration complexity $O( n \frac{\lambda_{\text{avg}} }{ \lambda_{\min}  } \log (1/\epsilon) )$,
where $ \lambda_{\text{avg} } $ and $ \lambda_{\min}   $ are the average eigenvalue and the minimum eigenvalue of the coefficient matrix, and one ``iteration'' here means a cycle of updating all blocks.
 This improves an existing bound of $O( n^2  \frac{\lambda_{\text{avg}} }{ \lambda_{\min}  } \log (1/\epsilon) ) $ for RP-BCD  by a factor of $n$.
  Built on this result, we further show that RP-ADMM for solving linear systems achieves the same expected iteration complexity bound $O( n \frac{\lambda_{\text{avg}} }{ \lambda_{\min}  } \log (1/\epsilon) )$. 
  
Technically, we provide a simple and clean proof of the expected convergence, by
 applying a classical result on the eigenvalues of Jordan product. 
  For proving the expected convergence rate, we propose a new variant of the matrix AM-GM inequality conjecture, and prove a weaker version of this conjecture.

Our result shows that random permutation may be a good answer to the question ``how to apply the decomposition idea to solve constrained problems''. 
As multi-block BCD is widely used for large-scale unconstrained problems,
we expect multi-block RP-ADMM to be a good candidate for large-scale linearly constrained problems.
Our result provides one of the few direct analyzes of random permutation  in optimization algorithms, and offers an explanation of the mysterious gap between RP-ADMM and cyclic ADMM. As reflected by the proof,
the intuition is that random permutation provides ``3-level symmetrization'' that adjusts the spectrum of the update matrix.
Based on the analysis for RP-ADMM, we are able to improve the best known complexity of RP-BCD for equally-diagonal quadratic problems by a factor of $n$, when expressing the complexity only in terms of the quantity $\frac{ \lambda_{\text{avg}} }{ \lambda_{\min}  }$.


\subsection{Related Works}


This paper is a stronger version of a previous technical report Sun et al. \cite{sun2015expected} which was not published. Another related work is the paper Chen et al. \cite{chen2015convergence}, which modifies the proof of 
\cite{sun2015expected} to make it work with a quadratic objective function.


We highlight a few novel contributions of the current paper (neither in the original technical report \cite{sun2015expected} nor in the paper \cite{chen2015convergence}). 

  	  (i) The current paper provides a much simpler proof for the result of expected convergence. 
  	  	
 (ii)  The current paper provides the first convergence rate analysis of RP-ADMM.
 See Theorem \ref{Thm 4} and the proof in  Section \ref{subsec: thm 4 proof}, Section \ref{subsec: lemma 5 proof} and Section	\ref{subsec: proof of Claim of spec radius}.
 
 (iii) The current paper provides an improved convergence rate analysis of RP-BCD,
See Theorem \ref{Thm 3} and the proof in Section
\ref{subsec: thm 3 proof} and Section \ref{subsec: proof of AMGM weak}. 

  (iv) The current paper introduces a theory-motivated algorithm Bernoulli-ADMM, which reduces the sampling time yet still achieves the expected convergence. This update order has not appeared before even in other algorithm setups to our knowledge. 
See Section \ref{sec: BR ADMM} and Proposition \ref{prop: BR-ADMM convergence}. 
	

Besides the technical contributions, we want to emphasize that the current paper is not just adding new result to our previous technical report \cite{sun2015expected}, but actually completes a missing step of the story. 
  From a mathematical point of view, the most striking consequence  of our original proof is
  that the spectral radius of RP-BCD lies in a smaller region $(-1/3, 1)$. 
  It is natural to think that this fundamental fact should have an impact on the analysis of original RP-BCD. 
  Our current paper fills this gap by showing that this result can help build an $O(n)$ gap between the (expected) onvergence rate of RP-BCD and cyclic BCD. 
  A general message is that on one hand, to understand constrained optimization we have to understand unconstrained optimization (analyzing ADMM reduces to analyzing BCD); on the other hand, analyzing constrained optimization helps improve the understanding of unconstrained optimization (the analysis of ADMM leads to progress in BCD).
  We find this  interaction between unconstrained optimization (BCD) and  constrained optimization (ADMM)
  fascinating. 
 The whole story is only revealed in the current paper, but not in the previous technical report \cite{sun2015expected}  or Chen et al. \cite{chen2015convergence}.
 
 Besides the above unique aspects, the current paper inherits 
  some interesting numerical findings from the technical report Sun et al. \cite{sun2015expected}
which do not appear in Chen et al. \cite{chen2015convergence}.
 We find that cyclic ADMM diverges with probability 1 for many random distributions of data, thus showing that the seemingly surprising divergence behavior reported in \cite{chen2014direct} is quite common. However, it is easy to miss this finding if one uses the Gaussian distribution to generate data. 
  Another interesting finding is that the independently randomized version of ADMM diverges with probability 1  for Gaussian data but not for the counter-example in \cite{chen2014direct}, preventing us from analyzing the independently randomized version. 
 Without these findings, the motivation of studying RP-ADMM would be less clear.  See Section \ref{subsec: randomized ADMM} and Section \ref{sec: simulation}.

\subsection{Notation and Organization }
\emph{Notation.} For a matrix $X$, we denote $X(i,j)$ as the $(i,j)$-th entry of $X$,
 $\text{eig}(X)$ as the set of eigenvalues of $X$,
$\rho(X)$ as the spectral radius of $X$ (i.e. the maximum modulus of the eigenvalues of $X$),
$\| X\|$ as the spectral norm of $X$, and $X^T$ as the transpose of $X$.
When $X$ is block partitioned, we use $X[i,j]$ to denote the $(i,j)$-th block of $X$.
When $X$ is a real symmetric matrix, let $\lambda_{\max}(X)$ and $\lambda_{\min}(X)$ denote the maximum and minimum  eigenvalue of $X$ respectively. %
For two real symmetric matrices $X_1$ and $X_2$, $X_1 \succ X_2$ (resp.\ $X_1 \succeq X_2$) means $X_1 - X_2$ is positive definite (resp.\ positive
semi-definite). We use $I_{m }$ to denote the identity matrix with dimension $m$, and
we will simply use $I$ when it is clear from the context what the dimension is.
For square matrices $U_i \in \dR^{u_i \times u_i}, i=1,\dots, k$,
we denote $\text{Diag}( U_1, U_2, \dots, U_k )$ as the block-diagonal matrix with $U_i$
being the $i$-th diagonal block.


\emph{Organization.} In Section \ref{sec: model and algorithm}, we present three versions of randomized ADMM, with an emphasis on
RP-ADMM. In Section \ref{sec: main result}, we present our main results Theorem \ref{Thm 1}, Theorem \ref{Thm 2}
 and their proofs. 
 The subsequent sections are devoted to the proofs of the two technical results Lemma \ref{lemma 1} and Lemma \ref{lemma 2},
 which are used in the proof of Theorem \ref{Thm 2}.
 In particular, the proof of Lemma \ref{lemma 1} is given in Section \ref{section: proof of Lemma 1},
 and the proof of Lemma \ref{lemma 2} is given in Section \ref{sec: proof of Lemma 2 for n-block}.


\section{ Algorithms }\label{sec: model and algorithm}
In this section, we will present both randomly permuted and independently randomized versions of  ADMM for solving \eqref{optimization, sum fi},
and specialize RP-ADMM for solving a square system of equations.
We also present a rather novel algorithm Bernoulli-randomized ADMM (motivated by our proof). 

	\subsection{Randomly Permuted ADMM}
 In this subsection, we first propose RP-ADMM for solving the general optimization problem \eqref{optimization, sum fi},
then we present the update equation of RP-ADMM for solving a linear system of equations. 

Define $\Gamma$ as 
\begin{equation}\label{Gamma def}
\Gamma \triangleq  \{ \sigma \mid \sigma \text{ is a permutation of } \{1,\dots, n\} \}.
\end{equation}
 At each round, we draw a permutation $\sigma$ of $\{1,\dots, n \}$ uniformly at random from $\Gamma$, and update the primal variables in the order of
 the permutation, followed by updating the dual variables in a usual way.
 Obviously, all primal and dual variables are updated exactly once at each round.
 See Algorithm \ref{Algorithm: n-block RP-ADMM} for the details of RP-ADMM.
 Note that with a little abuse of notation, the function $\mathcal{L}( x_{\sigma(1)}, x_{\sigma(2)}, \dots, x_{\sigma(n)}; \mu  )$
 in this algorithm should be understood as  $\mathcal{L}( x_1, x_2, \dots, x_n; \mu  )$.
 For example, when $n=3$ and $\sigma = (231)$, $\mathcal{L}( x_{\sigma(1)}, x_{\sigma(2)}, x_{\sigma(3)}; \mu  ) = \mathcal{L}( x_2, x_3, x_1; \mu )$ should be understood as $\mathcal{L}( x_1, x_2, x_3; \mu  )  $.
\begin{algorithm}[htb]
\caption{\normalsize $n$-block Randomly Permuted ADMM (RP-ADMM) } 
\label{Algorithm: n-block RP-ADMM}
\begin{algorithmic}
\STATE Initialization: $ x_i^0 \in \dR^{d_i \times 1}, i=1,\dots, n; \ \mu^0 \in \dR^{ N \times 1}$.
\STATE Round $k$ ($k = 0, 1,2,\dots$):
\STATE 1) Primal update.
\STATE  $\text{} \quad  $  Pick a permutation $\sigma$ of $\{1, \dots, n \}$  uniformly at random.
\STATE  $\text{} \quad  $  For $ i =1,\dots, n$,  compute $x_{\sigma(i) }^{k+1} $ by 
\begin{equation}\label{RP-ADMM, primal update, general objective}
 x_{\sigma(i) }^{k+1} = \arg \min_{ x_{\sigma(i)} \in \mathcal{X}_{\sigma(i)} } \mathcal{L}( x_{\sigma(1) }^{k+1}, \dots, x_{\sigma(i-1) }^{k+1}, x_{\sigma(i)},
 x_{\sigma(i+1) }^{k}, \dots,  x_{\sigma(n) }^{k}; \mu^k)
 \end{equation}
\vspace{-0.4cm}
\STATE 2) Dual update.  Update the dual variable by 
  \begin{equation}\label{RP-ADMM, dual update}
   \mu^{k+1} = \mu^k - \beta ( \sum_{i=1}^n A_i     x_i^{k+1} - b ) .
    \end{equation}
\end{algorithmic}
\end{algorithm}
\vspace{-0.1cm}

\subsubsection{Optimization Formulation of Solving a Linear System of Equations}\label{subsec: RP-ADMM linear system}
Consider a special case of \eqref{optimization, sum fi} where
$f_i = 0$, $\mathcal{X}_i = \dR^{d_i}, \forall i$ and $ N = \sum_i d_i $ (i.e. the constraint is a square system of equations).
Then problem \eqref{optimization, sum fi} becomes
\begin{equation}\label{optimization formulation}
\begin{split}
  \min_{x \in \dR^N }  \quad   &    0 ,   \\
    \st  \quad  &   A_1 x_1 + \cdots +  A_n x_n = b,
\end{split}
\end{equation}
where $A_i \in \dR^{N \times d_i}, x_i \in \dR^{d_i \times 1}, b\in \dR^{N \times 1}$.
Solving this feasibility problem (with $0$ being the objective function) is equivalent to solving
a linear system of equations
\begin{equation}\label{linear system; n-block}
A x = b,
\end{equation} where $A = [A_1,\dots, A_n] \in \dR^{N \times N}, x = [x_1^T, \dots, x_n^T]^T \in \dR^{N \times 1}, b \in \dR^{N \times 1}$.

Throughout this paper, we assume $A$ is non-singular.
Then the unique solution to \eqref{linear system; n-block} is $ x = A^{-1} b$,
and problem \eqref{optimization formulation} has a unique primal-dual optimal solution  $(x, \mu ) = (A^{-1}b, 0)$.
The augmented Lagrangian function \eqref{augmented Lag func} for the optimization problem \eqref{optimization formulation} becomes
\begin{equation}\label{Aug Lag for linear system}
  \mathcal{L}(x, \mu ) = -\mu^T( Ax - b ) + \frac{ \beta }{2}\| Ax - b \|^2.
\end{equation}
Throughout this paper, we assume $\beta  = 1$; note that our algorithms and results can be extended to any $\beta > 0$ by simply scaling $\mu$.

\subsubsection{Example of $3$-block ADMM}
Before presenting the update equation of general RP-ADMM for solving \eqref{optimization formulation},
 we consider a simple case $N=n = 3, d_i = 1, \forall i$ and $\sigma = (123)$, and let $a_i = A_i \in \dR^{3 \times 1}$.
 The update equations \eqref{RP-ADMM, primal update, general objective} and \eqref{RP-ADMM, dual update} can be rewritten as
\begin{equation}\nonumber
\begin{split}
   -a_1^{T} \mu^k & + a_1^T(a_1 x_1^{k+1} + a_2 x_2^{k} + a_3 x_3^k - b) = 0, \\
   -a_2^{T} \mu^k & + a_2^T(a_1 x_1^{k+1} + a_2 x_2^{k + 1} + a_3 x_3^k - b) = 0, \\
   -a_3^{T} \mu^k & + a_3^T(a_1 x_1^{k+1} + a_2 x_2^{k + 1} + a_3 x_3^{k+1 } - b ) = 0,  \\
   (a_1     x_1^{k+1}  & +   a_2 x_2^{k + 1} + a_3 x_3^{k+1 } - b) + \mu^{k+1} - \mu^k = 0.
\end{split}
\end{equation}
Denote $ y^k = [x_1^{k}; x_2^{k}; x_3^{k}; (\mu^k)^T ]  \in \dR^{6 \times 1}$, then the above update equation becomes
\begin{equation}\label{update matrix of 123}
   \begin{bmatrix}
   a_1^T a_1 &     0     &      0     & 0   \\
   a_2^T a_1 & a_2^T a_2 &      0     & 0   \\
   a_3^T a_1 & a_3^T a_2 & a_3^T a_3  & 0   \\
      a_1    &   a_2     &      a_3   & I_{3\times3} \\
   \end{bmatrix}
   y^{k+1} =
   \begin{bmatrix}
   0 & -a_1^T a_2 &  -a_1^Ta_3     & a_1^T   \\
   0 &     0      &  -a_2^Ta_3     & a_2^T   \\
   0 &     0      &      0         & a_3^T   \\
   0 &     0      &       0        &  I_{3\times3} \\
   \end{bmatrix}
    y^k +  \begin{bmatrix}
  A^T b  \\ b
   \end{bmatrix}.
\end{equation}
Define
\begin{equation}\label{L123 3-block ADMM}
 L \triangleq  \begin{bmatrix}
   a_1^T a_1 &     0     &      0        \\
   a_2^T a_1 & a_2^T a_2 &      0      \\
   a_3^T a_1 & a_3^T a_2 & a_3^T a_3     \\
   \end{bmatrix},  \quad
    R \triangleq  \begin{bmatrix}
   0 & -a_1^T a_2 &  -a_1^Ta_3       \\
   0 &     0      &  -a_2^Ta_3       \\
   0 &     0      &      0           \\
   \end{bmatrix}.
   \end{equation}
The relation between $L$ and $R$ is $$ L - R = A^T A. $$  Define
   \begin{equation}
     \bar{L} \triangleq \begin{bmatrix}
       L  &  0   \\
       A  & I_{3\times3} \\
   \end{bmatrix}, \quad
    \bar{R} \triangleq
   \begin{bmatrix}
  R       & A^T   \\
   0      &  I_{3\times3} \\
   \end{bmatrix},  \quad \bar{b} = \begin{bmatrix}
  A^T b  \\ b
   \end{bmatrix}
   \end{equation}
  then the update equation  \eqref{update matrix of 123}  becomes  $ \bar{L} y^{k+1} = \bar{R} y^k + \bar{ b } $, i.e.
\begin{equation}\label{update matrix of 123, concise}
      y^{k+1} = (\bar{L})^{-1} \bar{R} y^k + \bar{L}^{-1} \bar{b}.
 \end{equation}
As a side remark, reference Chen et al. \cite{chen2014direct} provides a specific example of $A \in \dR^{3 \times 3}$ so that
 $\rho( (\bar{L})^{-1} \bar{R} ) > 1$, which implies the divergence of the above iteration if the update order $\sigma=(123)$ is used all the time.
  This counterexample disproves the convergence of cyclic 3-block ADMM.

\subsubsection{General Update Equation of RP-ADMM}

In general, for the optimization problem \eqref{optimization formulation}, the primal update \eqref{RP-ADMM, primal update, general objective} becomes
 \begin{equation}\label{RP-ADMM, primal update, temp}
 -A_{ \sigma(i) }^{T} \mu^k  + A_{\sigma(i) }^T( \sum_{j=1}^i   A_{\sigma(j) } x_{\sigma(j) }^{k+1} +  \sum_{l= i+1}^n   A_{\sigma(l) } x_{\sigma(l) }^{k} - b ) = 0  , \ i=1,\dots, n.
 \end{equation}
Replacing $\sigma(i), \sigma(j), \sigma(l)$ by $i,j, l$, we can rewrite the above equation as
 \begin{equation}\label{RP-ADMM, primal update}
  -A_{ i }^{T} \mu^k  + A_{ i }^T(
  \sum_{ \sigma^{-1}(j) \leq \sigma^{-1}(i) }   A_{ j } x_{j }^{k+1} +  \sum_{\sigma^{-1}(l) > \sigma^{-1}(i)}   A_{ l } x_{ l }^{k} - b ) = 0,
  \ i=1, \dots, n,
 \end{equation}
 where $\sigma^{-1}$ denotes the inverse mapping of a permutation $\sigma$, i.e. $\sigma(i) = t  \Leftrightarrow i = \sigma^{-1}(t)$.
Denote the output of Algorithm \ref{Algorithm: n-block RP-ADMM}  after round $(k-1)$ as
\begin{equation}\label{definition of yk}
  y^k \triangleq [x^k; \mu^k] = \begin{bmatrix} x_1^k ; \dots ; x_n^k  ; \mu^k \end{bmatrix} \in \dR^{2N \times 1}.
\end{equation}
The update equations of Algorithm \ref{Algorithm: n-block RP-ADMM} for solving \eqref{optimization formulation},
i.e. \eqref{RP-ADMM, primal update} and \eqref{RP-ADMM, dual update},
can be written in the matrix form as (when the permutation is $\sigma$ and $\beta = 1$)
\begin{equation}\label{update eqn of n-block-RADMM}
   y^{k+1} = \bar{L}_{\sigma}^{-1} \bar{R}_{\sigma} y^{k} + \bar{L}_{\sigma}^{-1} \bar{b},
\end{equation}
where   $\bar{L}_{\sigma}, \bar{R}_{\sigma}, L_{\sigma}, R_{\sigma}, \bar{b}$
are defined by
\begin{equation}\label{barL def, n-block}
     \bar{L}_{\sigma} \triangleq \begin{bmatrix}
       L_{\sigma}  &  0   \\
       A  & I_{N \times N} \\
   \end{bmatrix}, \quad
    \bar{R}_{\sigma} \triangleq
   \begin{bmatrix}
  R_{\sigma}       & A^T   \\
   0      &  I_{N \times N} \\
   \end{bmatrix}, \quad   \bar{b} = \begin{bmatrix}
  A^T b  \\ b
   \end{bmatrix},
   \end{equation}
   in which $L_{\sigma} \in \dR^{N \times N}$ has $n \times n$ blocks and the $(i,j)$-th block is defined as
\begin{equation}\label{Lsigma def, temp}
   L_{\sigma}[ i , j ] \triangleq \begin{cases}
 A_i^T A_j &  \sigma^{-1}(j) \leq \sigma^{-1}(i)  ,    \\
0  &  \text{otherwise}.
\end{cases}
 \end{equation}
 and $R_{\sigma}$ is defined as
  \begin{equation}\label{Rsigma, def}
  R_{\sigma} \triangleq L_{\sigma} - A^T A.
  \end{equation}

Another expression of $L_{\sigma}$, equivalent to \eqref{Lsigma def, temp}, is the following:
 \begin{equation}\label{Lsigma def}
   L_{\sigma}[ \sigma(i) , \sigma(j) ] \triangleq \begin{cases}
  A_{\sigma(i)}^T A_{\sigma(j)} &  j \leq i,    \\
0  &  j > i ,
\end{cases}
 \end{equation}
To illustrate the above expression of $L_{\sigma}$, we consider the $n$-coordinate case that $d_i = 1, \forall i$.
In this case, each block $x_i$ is a single coordinate, and each $A_i$ is a vector. Denote $a_i \triangleq A_i \in \dR^{N \times 1}$.
Let $ L_{\sigma}( k,l ) $ denote the $(k, l)$-th entry of the matrix $ L_{\sigma}$, then
the definition \eqref{Lsigma def} becomes
 \begin{equation}\label{Lsigma def, n-coordinate}
   L_{\sigma}( \sigma(i) , \sigma(j) ) \triangleq \begin{cases}
 a_{\sigma(i)}^T a_{\sigma(j)} &  j \leq i,  \\   
0  &  j > i ,
\end{cases}
 \end{equation}
A user-friendly rule for writing $ L_{\sigma} $ is described as follows (use $\sigma = (231)$ as an example).
Start from a zero matrix.
First, find all reverse pairs of $\sigma$; here, we say $(i,j)$ is a reverse pair if $i$ appears after $j$ in $\sigma$.
For the permutation $(231)$, all the reverse pairs are
 $(1,3), (3,2) $ and $(1,2)$.
Second, in the positions corresponding to the reverse pairs,
write down the corresponding entries of $A^T A$, i.e. $a_1^T a_3, a_3^T a_2 $ and $a_1^T a_2$, respectively.
At last, write $a_i^T a_i$ in the diagonal positions. Using this rule, we can write down the expression of $L_{(231)}$ as 
$$
  L_{(231)} =
   \begin{bmatrix}
   a_1^T a_1 &    a_1^T a_2   &     a_1^T a_3        \\
   0 & a_2^T a_2  &      0      \\
   0 & a_3^T a_2 &  a_3^T a_3     \\
   \end{bmatrix}.
$$
A user-friendly rule to quickly check the correctness of an expression of $L_{\sigma}$ is the following (still take $\sigma = (231)$ as an example).
According to the order of the permutation $(231)$,
the $2$nd row, the $3$rd row and the $1$st row should have a strictly decreasing number of zeros ($2$ zeros, $1$ zero and no zero).
In contrast, the $2$nd column, the $3$rd column and the $1$st column should have a strictly  increasing number of zeros.

For the general case that $d_i \geq 1, \forall i$, we can write down the block partitioned  $L_{\sigma}$ in a similar way.
For example, when $n = 3$ and $\sigma = (231)$, we have
$$
  L_{(231)} =
   \begin{bmatrix}
   A_1^T A_1 &    A_1^T A_2   &     A_1^T A_3        \\
   0 & A_2^T A_2 &      0      \\
   0 & A_3^T A_2 &  A_3^T A_3     \\
   \end{bmatrix}.
$$

\subsection{Randomly Permuted BCD}
RP-ADMM is a generalization of RP-BCD.
In fact, when the constraint does not exist, RP-ADMM reduces to RP-BCD.
In this subsection, we present RP-BCD for solving convex quadratic problems. 
Note that RP-ADMM for solving linear systems and RP-BCD fo solving quadratic problems are two extremes of general RP-ADMM: in the former case the objective function is zero, and in the latter case the constraint is zero. Interestingly, the two extreme cases are related as the expected iteration matrix of RP-BCD appears as a component of the expected iteration matrix of RP-ADMM. We will show later that their eigenvalues are closely related. 

Consider a special case of \eqref{optimization, sum fi} where
$f(x) = \frac{1}{2}\| Ax - b \|^2$,  $\mathcal{X}_i = \dR^{d_i}, \forall i$ and there is no constraint.
 With abuse of notation, we use $A$ to denote the coefficient matrix, while in the original formulation $A$ denotes the constraint matrix. We ``recycle'' the notation $A$ so that we can build a connection with RP-ADMM for solving linear systems later. 
Assume  $N = \sum_i d_i$. 
Then problem \eqref{optimization, sum fi} becomes a least-squares problem
\begin{equation}\label{quadratic optimization formulation}
\begin{split}
\min_{x \in \dR^N }  \quad   & \frac{1}{2}  \| Ax - b \|^2 = \frac{1}{2}  \| A_1x_1 + \dots + A_n x_n - b \|^2
\end{split}
\end{equation}
where $A_i \in \dR^{N \times d_i}, x_i \in \dR^{d_i \times 1}, b\in \dR^{N \times 1}$.
Similar to Section \ref{subsec: RP-ADMM linear system}, we assume $A$ is non-singular.
Then the unique solution to \eqref{linear system; n-block} is $ x = A^{-1} b$.

In the augmented Lagrangian function given in \eqref{Aug Lag for linear system}, if we delete the first term which depends on the dual variable $\mu$, we obtain the quadratic function $\frac{1}{2}\| Ax - b\|^2$. Thus if we eliminate the dual variable $\mu$ in the update equations of RP-ADMM, we will obtain the update equations for RP-BCD. 
Suppose $x^k$ is the iterate after the k-th epoch (i.e. go through all coordinates once), and $\sigma$ is the order used in the $k$-th iteration, then, as a simpler version of  \eqref{update eqn of n-block-RADMM},
we have
\begin{equation}\label{update eqn of n-block-RP-CD}
x^{k+1} = L_{\sigma}^{-1} R_{\sigma} x^{k} + L_{\sigma}^{-1} b,
\end{equation}
where $L_{\sigma}$ and $R_{\sigma}$ are defined as in \eqref{Lsigma def} and \eqref{Rsigma, def}, and $\sigma $ is a random permutation.


 \subsection{Residual Trick for Efficient Implementation of ADMM and BCD}\label{subsec: residual trick}
 We note here that when $d_i = 1, \forall i$ (in this case BCD becomes CD),  per-epoch computation time of 
 ADMM and CD (no matter what order) is $O(n^2)$; or in other words, per-coordinate-update time is $O(n)$.
 For instance, updating  $x^{k+1}$ by \eqref{update eqn of n-block-RP-CD} in RP-BCD or updating $y^{k+1}$ by  \eqref{update eqn of n-block-RADMM} in RP-ADMM  only takes time $O(n^2)$. 
 As mentioned in Section 3.1 of \cite{nesterov2012efficiency}, the trick is to keep track of the residual.
 For both efficient practical implementation and calculation of computation complexity, one should use this residual trick, but for the ease of theoretical analysis we use the matrix update forms  \eqref{update eqn of n-block-RADMM} and  \eqref{update eqn of n-block-RP-CD} in this paper; there is no contradiction as our theory only depends on the value of $x^k$ but not the specific procedure to compute $x^k$.
 
For completeness, we briefly explain how this trick works in our settings. 
Suppose $d_i = 1, \forall i$, and we use CD methods to solve \eqref{quadratic optimization formulation} with a certain update order (could be any order, such as cyclic, randomized or randomly permuted). Suppose the coordinate $i$ is picked, then  $x_i$ is updated by
by 
\begin{equation}\label{LS original update}
  x_i^+ =  \frac{ 1 }{ A_i^T A_i } [  A_i^T (  b -  A_{-i} x_{-i} )  ] ,
\end{equation}
where $A_{-i}$ contains all columns of $A$ except $A_i$, $x_{-i}$ contains all elements of $x$ except $x_i$ and represents the current
values, and $x_i^+$ represents the new value. 
A straightforward implementation of \eqref{LS original update} requires multiplying $x_{-i}$ by $ A_{-i} $ which takes $O(n^2)$ operations. With the residual trick (e.g. \cite{nesterov2012efficiency}), we introduce the residual $r =  Ax - b $, and replace \eqref{LS original update} by 
\begin{align*}
	 x_i^+  &  =  x_i   -  \frac{ 1 }{ A_i^T A_i }   A_i^T r  ,   \\
	  r^+    & = r + A_i ( x_i^+ - x_i ).
\end{align*}
Now the calculation of $x_i^T $ and $r^+$ takes time $O(n)$, and thus one epoch of BCD takes time $O(n^2)$.
 The same trick can be applied to the primal update of ADMM;
with this trick, the dual update \eqref{RP-ADMM, dual update} can be rewritten as $\mu^+ = \mu - \beta r $ which takes time $O(n)$,
and thus one epoch of ADMM takes time $O(n^2)$. 

Finally, when $d_i > 1$, similar update equations can still be used except a minor difference that  $  \frac{ 1 }{ A_i^T A_i } $ should be replaced by  $ ( A_i^T A_i )^{-1}  $. In a special case that $d_i = d, \forall i$ and $N = d n$, 
each iteration of BCD takes time $O(N d + d^3)$ and each epoch takes time $O(N^2 + N d^2)$. This cost can be reduced if we use BCGD (i.e. not solving the subproblem exactly but updating each block of variables by a gradient step). In order not to make the paper more complicated, we will not discuss the inexact versions of BCD and ADMM in this paper.  

\subsection{Two Versions of Independently Randomized ADMM}\label{subsec: randomized ADMM}
In this subsection, we present two other versions of randomized ADMM which can be divergent according to simulations.
The failure of these versions makes us focus on analyzing RP-ADMM in this paper.
These versions can be viewed as natural extensions of R-BCD (randomized BCD) \cite{leventhal2010randomized} and
	\cite{nesterov2012efficiency}. 

In the first algorithm, called primal-dual randomized ADMM (PD-RADMM), the whole dual variable is viewed as the $(n+1)$-th block.
In particular, at each iteration, the algorithm draws one index $i$ from $\{1,\dots, n, n+1 \}$, then performs the following update:
if $i \leq n$, update the $i$-th block of the primal variable; if $i = n+1$, update the whole dual variable. The details are given in Algorithm \ref{Algorithm: PD-RADMM}. We have tested PD-RADMM for the counter-example given in Chen et al. \cite{chen2014direct}, and found that
PD-RADMM always diverges (for random initial points).

A variant of PD-RADMM has been proposed in Hong et al. \cite{hong2014block} with two differences:
first, instead of minimizing the augmented Lagrangian $\mathcal{L}$, that algorithm minimizes a strongly convex upper bound of $\mathcal{L}$;
second, that algorithm uses a diminishing dual stepsize.
With these two modifications,  \cite{hong2014block} shows that each limit point of the sequence generated by their algorithm is a primal-dual optimum with probability 1. Note that \cite{hong2014block} also proves the same convergence result for the cyclic version of multi-block ADMM with these two modifications, thus it does not show the benefit of randomization.

\begin{algorithm}[htb]
\caption{\normalsize Primal-Dual Randomized ADMM (PD-RADMM) }
\label{Algorithm: PD-RADMM}
\begin{algorithmic}
\STATE Iteration $t$ ($t = 0,1 , 2, \dots $): 
\STATE $\quad $ Pick $i \in \{ 1, \dots, n, n+1\}$ uniformly at random;
\STATE   $\quad \quad \quad $ If $1 \leq i \leq n$:
\STATE  $\quad  \quad \quad \quad \quad $ $     x_i^{t+1} = \arg \min_{x_i \in \mathcal{X}_i } \mathcal{L}( x_1^{t}, \dots, x_{i-1}^{t},
x_i,  x_{i+1}^{t}, \dots,  x_n^t; \mu^t) , $
\STATE  $\quad  \quad \quad \quad \quad $   $ x_j^{t+1} = x_j^t, \ \forall \  j\in \{1,\dots, n \}\backslash \{ i\}, $
\STATE  $\quad  \quad \quad \quad \quad $   $   \mu^{t+1} = \mu^t.   $
\STATE   $\quad  \quad \quad $ Else If $i = n + 1$:
\STATE  $\quad  \quad \quad \quad \quad $  $   \mu^{t+1} = \mu^t - \beta( \sum_{i=1}^n A_i     x_i^{t+1} - b ) , $
\STATE  $\quad  \quad \quad \quad \quad $  $ x_j^{t+1} = x_j^t, \ \forall \ j\in \{1,\dots, n \}. $
\STATE  $\quad  \quad \quad $ End
\end{algorithmic}
\end{algorithm}

In the second algorithm, called primal randomized ADMM (P-RADMM), we only perform randomization for the primal variables.
In particular, at each round, we first draw $n$ independent random variables $j_1,\dots, j_n$ from the uniform distribution of $\{1,\dots, n \}$ and
update $x_{j_1}, \dots, x_{j_n}$ sequentially, then update the dual variable in the usual way.
The details are given in Algorithm \ref{Algorithm: P-RADMM}.
This algorithm looks quite similar to RP-ADMM as they both update $n$ primal blocks at each round; the difference is that
RP-ADMM samples \emph{without replacement} while this algorithm P-RADMM samples \emph{with replacement}.
In other words, RP-ADMM updates each block exactly once at each round, while
 P-RADMM may update one block more than one times or does not update one block at each round.

 We have tested P-RADMM in various settings.
 For the counter-example given in Chen et al. \cite{chen2014direct}, we found that P-RADMM does converge.
However, if $n \geq 30$ and $A$ is a Gaussian random matrix (each entry is drawn i.i.d. from $\mathcal{N}(0,1)$),
then P-RADMM diverges in almost all cases we have tested. This phenomenon is rather strange since for random Gaussian matrices $A$
the cyclic ADMM actually converges (according to simulations). An implication is that randomized versions do not always outperform their deterministic counterparts in terms of convergence.

Since both Algorithm \ref{Algorithm: PD-RADMM} and Algorithm \ref{Algorithm: P-RADMM}
can diverge in certain cases, we will not further study them in this paper.
In the rest of the paper, we will focus on RP-ADMM (i.e. Algorithm \ref{Algorithm: n-block RP-ADMM}).

\begin{algorithm}[htb]
\caption{\normalsize Primal Randomized ADMM (P-RADMM) } 
\label{Algorithm: P-RADMM}
\begin{algorithmic}
\STATE Round $k$ ($k = 0,1, 2, \dots$):
\STATE 1) Primal update.
\STATE  $\text{} \quad  $  Pick $l_1, \dots, l_n$ independently from the uniform distribution of $\{1, \dots, n \}$.
\STATE  $\text{} \quad  $  For $ i =1,\dots, n$:
\STATE  $ \quad \quad \quad \  $    $t = k n + i - 1$,
   \STATE  $\quad  \quad \quad \ $ $     x_{l_i}^{t+1} = \arg \min_{ x_{l_i} \in \mathcal{X}_{l_i} } \mathcal{L}( x_1^{t}, \dots, x_{ l_i - 1}^{t},
x_{l_i},  x_{ l_i + 1}^{t}, \dots,  x_n^t; \mu^t) , $
\STATE  $\quad  \quad \quad \ $   $ x_{j }^{t+1} = x_j^t, \ \forall \  j\in \{1,\dots, n \}\backslash \{ l_i\}, $
\STATE  $\quad  \quad \quad \ $   $ \mu^{t+1} = \mu^t. $
\STATE  $\text{} \quad  $  End.
\STATE 2) Dual update.
\STATE  $\text{} \quad \quad \quad \quad  $  $  \mu^{(k+1) n} = \mu^{ kn} - \beta( \sum_{i=1}^n A_i     x_i^{(k+1)n} - b ) . $
\end{algorithmic}
\end{algorithm}

\subsection{Bernoulli-Randomized ADMM}\label{sec: BR ADMM}

To implement randomly permuted ADMM, one needs to sample from all blocks without replacement. 
To save the sampling time, we propose another algorithm which we call Bernoulli-randomized ADMM.
This algorithm is motivated by the proof of Theorem \ref{Thm 1}.  
This updating scheme can be applied to other algorithms such as SGD and coordinate descent methods. 

The new update order combines the  well-known double-sweep order and Bernoulli-randomization. The original double-sweep order is $(1,2,...,n-1, n,  n-1, n-2, ...,1)$, meaning that
$x_1, x_2, \dots, x_{n-1}, x_n, x_{n-1}, x_{n-2}, \dots, x_1 $ are updated sequentially in each ``cycle''. It combines the normal cyclic order $(1,2,\dots, n)$ and a reverse order $(n,n-1,\dots, 1)$. We propose the following updating scheme:  add a check box to each block, and in each cycle we perform the following operations.  
\begin{enumerate}

 \item Phase I: go through the blocks $x_1, x_2, \dots, x_n $ one by one sequentially as follows: for each block $x_i$, flip a fair coin and:
   \begin{enumerate}
   	\item if the outcome is ``head'', update the block $x_i$ and check the check box;
   	
   	\item if the outcome is ``tail'', do nothing about $x_i$ and uncheck the check box. 
\end{enumerate}
 \item Phase II: go through the blocks $x_n, x_{n_1}, \dots, x_1$ in the reverse order, and update $x_i$ if the box is unchecked. 

\end{enumerate}

Note that in each cycle we go through each block twice but update each block exactly once so that the number of totally updated blocks remains $n$. For example, when $n  = 5$, $(35421)$ is a possible update order, as shown in the following diagram.
 \begin{table}[!htbp]
	\centering	
	\begin{tabular}{c c c c  c  c  c }
	                         &	&       1                               &           2                                              &              3                           & 4    &     5    \\
	\text{Phase I}   &    \text{   begin}   &    \text{    skip}  $ \rightarrow $    &           \text{    skip}  $ \rightarrow $     &       \quad       3 $ \rightarrow $   &     \text{    skip}  $ \rightarrow $      &     5    \\  
	 &	&                                   &                                                   &                                &  &   $ \downarrow  $  \\
	\text{Phase II}   &    \text{end  }  &  1    &      $ \leftarrow $     2        &  $  \leftarrow $    \text{skip}   & $ \leftarrow $ 4    &    $ \leftarrow $    \text{skip}    \\
	\end{tabular}
\end{table}
Similarly, $(13542)$ is also a possible update order. 
 But $(13524)$ and $(35412)$ are not possible.  The set of all possible update orders is given by
$$
\Gamma_{\text{BR}} \triangleq \{ \sigma \in \Gamma \mid  \exists \, i  \in \{ 1,\dots, n - 1 \}  \text{ such that } \sigma(1) < \sigma(2) < \dots < \sigma(i) 
\text{ and } \sigma(i+1) > \dots > \sigma(n)  \}  ,
$$
where $\Gamma$ is the set of permutations of $\{ 1, 2, \dots, n \}$ as defined in \eqref{Gamma def}. 
In other words, a sequence from $\Gamma_{\text{BR} }$  is a concatenation of an increasing sequence  and a decreasing sequence. 
Note that the permutation $(1,2,...,n)$ is in $\Gamma_{\text{BR} } $ since it can be viewed as the concatenation of an increasing sequence $(1, 2, ..., n-1)$ and a ``decreasing sequence'' $(n)$, and we can let $i = n-1$ in the above definition to cover this case. 
Similarly, the permutation $(n, n-1, \dots, 1)$ is also in $ \Gamma_{\text{BR}} $ as $i = 1$ will cover this case.

The algorithm Bernoulli-randomized ADMM (BR-ADMM) is formally described below. We skip the epoch index $k$ since otherwise the notation would be cumbersome. 
\vspace{-.6cm}


\textit{\begin{algorithm}[htb]
	\caption{\normalsize $n$-block Bernoulli-Randomized ADMM (BR-ADMM) } 
	\label{Algorithm: n-block BR-ADMM}
	\begin{algorithmic}
		\STATE Initialization: $ x_i^0 \in \dR^{d_i \times 1}, i=1,\dots, n; \ \mu^0 \in \dR^{ N \times 1}$.
		\STATE Round $k$ ($k = 0, 1,2,\dots$):
		\STATE 1) Primal update.
		\STATE $\text{} \quad  $ Set $ c_i = 0, i=1,\dots, n $.
		\STATE  $\text{} \quad  $  \textbf{Phase I. }  
		\STATE  $\text{} \quad  $    For $ i =1, 2, \dots, n$:
		\STATE  $\text{} \quad  \quad $  Draw a random variable $\xi \sim \text{Bernnolli}(1/2)$, i.e. $ Pr(\xi = 1) =  Pr(\xi = 0) = 1/2 $.
		\STATE  $\text{} \quad  \quad $  If $\xi = 1$:  set $ c_i = 1$ and update $x_{ i } $ by 
		\begin{equation}\label{BR-ADMM, primal update, general objective}
		   x_{ i } \leftarrow \arg \min_{ x_i  \in \mathcal{X}_{ i } } \mathcal{L} (   x_1, \dots, x_{i-1}, x_i, x_{i+1}, \dots, x_n; \mu ).
		\end{equation}
				\STATE  $\text{} \quad  $  \textbf{Phase II. } 
			\STATE  $\text{} \quad  $ For $ i =n , n -1, \dots, 1$:  if $c_i = 0 $, update $x_{ i } $ by \eqref{BR-ADMM, primal update, general objective}.		 
		\STATE 2) Dual update.  Update the dual variable by 
		\begin{equation}\label{BR-ADMM, dual update}
		\mu \leftarrow \mu - \beta ( \sum_{i=1}^n A_i     x_i - b ) .
		\end{equation}
	\end{algorithmic}
\end{algorithm}}

For solving linear systems of equations, the update formula is the same as \eqref{update eqn of n-block-RADMM}, the update formula of RP-ADMM. The difference is that for RP-ADMM $\sigma$ can be an arbitrary permuation, while for BR-ADMM  there is some restriction on $\sigma$: it has to be a permuation in $\Gamma_{\text{BR} }$.

\section{ Main Results }\label{sec: main result}

\subsection{Expected Convergence of RP-ADMM}
Let $\sigma_i$ denote the permutation used in round $i$ of Algorithm \ref{Algorithm: n-block RP-ADMM}, which is a uniform random variable drawn from the set of permutations $\Gamma$.
After round $k$, Algorithm \ref{Algorithm: n-block RP-ADMM} generates a random output $y^{k+1}$, which depends on the observed draw of the random variable
\begin{equation}\label{def of xi}
  \xi_k = (\sigma_0, \sigma_1, \dots, \sigma_k).
\end{equation}
We will show that the expected iterate (the iterate $y^k$ is defined in \eqref{definition of yk})
\begin{equation}
 \phi^k =  E_{\xi_{k-1} }(y^{k})
\end{equation}
converges to the primal-dual solution of the problem \eqref{optimization formulation}. 
Although the expected convergence does not necessarily imply the convergence in a particular realization, it serves as an
evidence of convergence.
Our proof seems much different from and more difficult than previous proofs for other randomized methods, since random permutation, as well as spectral radius of non-symmetric matrices, are difficult objects to deal with -- not many existing mathematical tools are available to help \footnote{There has been some effort in using random matrix theory to tackle this problem but no progress has been reported to our knowledge. This is partially due to the fact that the desired result seems to be rather tight such that even a small relaxation can lead to failure.}.
Note that the extension of this result to the non-square full column-rank case is simple \footnote{Suppose $A$ 
is an $m \times n$ full column-rank matrix, where $m \geq n$, and the system $Ax = b$ is feasible. The update formula is $ y^{k+1} = (I - L_{\sigma}^{-1} A^T A ) y^k $, which is same as the update formula for solving a square system of equations $ \bar{A} x = b    $, where $\bar{A} \in \dR^{n \times n}$ is the square root matrix of the matrix $A^T A \in \dR^{n \times n}$.
Now the matrix $\bar{A}$ is a square invertible matrix, thus by applying the result for square system of equations, we can obtain the convergence of the sequence $\phi^k = E(y^k)$. }. 

\begin{thm}\label{Thm 1}
$ \text{}$	
Assume the coefficient matrix $A = [A_1,\dots, A_n]$ of the constraint in \eqref{optimization formulation} is a non-singular square matrix.
Suppose Algorithm \ref{Algorithm: n-block RP-ADMM} is used to solve problem \eqref{optimization formulation},
then the expected output converges to the unique primal-dual optimal solution to \eqref{optimization formulation},
i.e.
\begin{equation}\label{convergence to zero}
  \{ \phi^k \}_{k \rightarrow \infty} \longrightarrow  \begin{bmatrix}
       A^{-1}b  \\
       0 \\
   \end{bmatrix}.  
\end{equation}
\end{thm}

Since the update matrix does not depend on previous iterates, we claim (and prove in Section
\ref{sec: proof of Theorem 1}) that Theorem \ref{Thm 1}
holds if the expected update matrix has a spectral radius less than 1, i.e. if the following
Theorem \ref{Thm 2} holds.


\begin{thm}\label{Thm 2} $ \text{}$	
Suppose $A = [A_1, \dots, A_n] \in \dR^{N\times N}$ is non-singular, and $\bar{L}_{\sigma}^{-1} , \bar{R}_{\sigma} $
are defined by \eqref{barL def, n-block} for any permutation $\sigma$.
Define
\begin{equation}\label{M definition, theorem}
   M \triangleq  E_{\sigma}(\bar{L}_{\sigma}^{-1}  \bar{R}_{\sigma}  ) = \frac{1}{n!} \sum_{\sigma \in \Gamma} (\bar{L}_{\sigma}^{-1}  \bar{R}_{\sigma}  ),
 \end{equation}
where the expectation is taken over the uniform random distribution over $\Gamma$, the set of permutations of $\{1,2,\dots, n \}$. 
Then the spectral radius of $M$ is smaller than $1$, i.e.
\begin{equation}
  \rho(M) < 1.
\end{equation}
\end{thm}

\begin{rmk} $ \text{}$	
 For the counterexample in Chen et al. \cite{chen2014direct} where $A = [1,1,1; 1,1,2; 1,2,2]$, it is easy to verify that
$\rho(M_{\sigma}) > 1.02$ for any permutation $\sigma$ of $(1,2,3)$. Interestingly, Theorem \ref{Thm 2} shows that even if each $M_{\sigma}$
is ``bad'' (with spectral radius larger than $1$), the average of them is always ``good'' (with spectral radius smaller than $1$).
\end{rmk}

Theorem \ref{Thm 2} is just a linear algebra result, and can be understood even without knowing the details of the algorithm.
However, the proof of Theorem \ref{Thm 2} is rather non-trivial.
This proof will be provided in Section \ref{sec: proof of Thm 2}, and the technical results used
in this proof will be proved in Section \ref{section: proof of Lemma 1} and Section \ref{sec: proof of Lemma 2 for n-block}.

The convergence rate of RP-ADMM for solving linear systems of equations is closely related to the convergence rate of RP-BCD (randomly permuted BCD) for solving quadratic problems. We will discuss their relation
and how our results in this paper improve our understanding for RP-BCD. 


A similar convergence result holds for BR-ADMM proposed in Section \ref{sec: BR ADMM}, as presented below. The proof is a simple modification of the proof of Theorem \ref{Thm 1},  and can be found
 in Section \ref{appen: proof of BR ADMM}.

\begin{prop}\label{prop: BR-ADMM convergence} $ \text{}$	
	Assume the coefficient matrix $A = [A_1,\dots, A_n]$ of the constraint in \eqref{optimization formulation} is a non-singular square matrix.
	Suppose Algorithm \ref{Algorithm: n-block BR-ADMM} is used to solve problem \eqref{optimization formulation},
	then the expected output converges to the unique primal-dual optimal solution to \eqref{optimization formulation}.
\end{prop}

\subsection{Expected Convergence Rate of RP-ADMM and RP-BCD}\label{sec: discussion of convergence rate}
There is  a close relation between RP-ADMM for solving linear systems and RP-CD for solving quadratic problems (see Lemma \ref{lemma 2}). Thus it is not surprising that we need to understand RP-BCD before understanding RP-ADMM. 
We will first present an expected convergence rate of RP-BCD (in terms of the expected iterates) for solving quadratic problems, which improves the best existing
convergence rate (one type of rates, to be precise) by a factor of $n$ \footnote{Rigorously speaking, this is not a fair comparison as the complexity of C-CD is deterministic complexity.}. 
The result is proved via establishing a weak version of matrix AM-GM inequality. 
This result also establishes a large gap of $O(n)$ between RP-BCD and C-BCD (cyclic BCD).
Second, built upon the result for RP-BCD, we establish a  convergence rate of RP-ADMM which is similar to
RP-BCD and also $n$ times better than that of C-BCD.


The first result is about the expected convergence rate of RP-BCD for the case $A_i^T A_i = I$.
This assumption is made so that the expression is simple, and the case for general $A_i$ is given in the next result.
\begin{thm}\label{Thm 3}{(rate of RP-BCD for quadratic functions with identity diagonal blocks)} $ \text{}$	
	Assume the coefficient matrix $A = [A_1,\dots, A_n]$ is a non-singular square matrix, and
	$
	  A_i^T A_i = I, \forall \; i. 
	  $
	Suppose RP-BCD is used to solve problem \eqref{quadratic optimization formulation}, where $x^k$ denotes the variable
	after $k$ epochs (each epoch represents one cycle of updating all coordinates).
	Denote the unique optimal solution as $x^* = A^{-1}b$. 
	Then
	\begin{equation}\label{CD rate}
	\|  E(x^k)- x^*  \| \leq  \max \left\{ 1 - \frac{ 1}{n} \lambda_{\min}(AA^T)  ,   \frac{1}{3} \right\}^k \| x^0 - x^*\| .
	\end{equation}
\end{thm}

To put this convergence rate result in the context, we consider the simple case that each $d_i = 1$, i.e., each block consists of a single coordinate. In this case, every diagonal entry of $A^T A$ is $1$, thus the average eigenvalue of $A^T A$ is $1$.
Throughout the paper, we consider the total computation complexity 
\footnote{The computation complexity equals the iteration complexity times the per-iteration cost. We do not present iteration complexity since there may be confusion about whether ``one iteration'' means $n$ coordinate updates or $1$ coordinate update. Presenting iteration complexity is better if one considers a general convex problem, but then one needs to discuss the per-iteration cost. We are considering quadratic problems throughout the paper, so we feel it is more clear to stick to computation complexity.}; note that we assume the residual trick as described in   \ref{subsec: residual trick} is always used for all methods. 
 
 Our Theorem \ref{Thm 3} provides an expected computational complexity upper bound $O(n^3 \kappa_{\text{CD}}\log \frac{1}{\epsilon} )$ for RP-CD, since each epoch takes $ O(n^2) $ time and it requires $O(n^2 \log \frac{1}{\epsilon})$ epochs to achieve
 error $\epsilon$ according to \eqref{CD rate}}. 
 It is known that the computational complexity
of R-CD (randomized coordinate descent) to achieve relative accuracy $\epsilon$ 
\footnote{Here, the relative accuracy $\epsilon$ means $\|  E(x^k)- x^*  \|/\| x^0 - x^*\| $
	or $\|  E( f( x^k) )- f^*  \|/\| f(x^0) - f^*\| $. } is
$ O(n^2 \kappa_{\text{CD}} \log \frac{1}{\epsilon} ) $, where $\kappa_{\text{CD} } = \lambda_{\text{avg}} (A^T A)  / \lambda_{\min} (A^T A) =  1 / \lambda_{\min} (A^T A) $ is the ratio of the average eigenvalue over the minimum eigenvalue.  
It was recently shown that in terms of $\kappa_{\text{CD}} $ and $n$ only, the worst-case complexity of C-CD (cyclic CD)
is $O(n^4  \kappa_{\text{CD}} \log \frac{1}{\epsilon} )$, which is $n^2$ times worse than R-CD and $n$ times worse than GD. 
This shows a large gap between C-CD and R-CD in the worst case.

It was widely conjectured that RP-CD is at least as fast as R-CD, but this conjecture is considered to be rather difficult to prove.  For a special class of matrices, recent works \cite{lee2016random,wright2017analyzing} validated the conjecture.
However, to our knowledge, even for a general quadratic function with equal diagonal entries $1$, the previously best known convergence rate of RP-CD is almost the same as C-CD (see \cite{sun2015improved}\cite{sun2016worst}), which can be $n^2$ times worse than that of R-CD. Our Theorem \ref{Thm 3} provides an expected
computational complexity upper bound $O(n^3 \kappa_{\text{CD}}\log \frac{1}{\epsilon} )$ for RP-CD, which is $n$ times faster than C-CD and $n$ times slower than R-CD.
	This improves the best existing rate by a
factor of $n$ \footnote{Note that this ``improvement'' is valid when the convergence rate is characterized by  only $\kappa_{\text{CD} }$ and $n$. It is common to use other parameters such as the maximum eigenvalue to characterize the convergence rate (see \cite{sun2016worst} for a detailed discussion), and our result here does not provide improvement for other kinds of convergence rate.}.
We summarize the comparison of the complexity for C-CD, R-CD and RP-CD in Table \ref{table simple compare}. 
 \begin{table}
 	\centering	  %
 	\caption{ Worst-case computation complexity comparison, using only $\kappa_{\text{CD}} $ as parameter, for equal-diagonal quadratic case (ignore $O(\log \frac{1}{\epsilon})$ factor), and consider the error in the expected iterates for RP-CD}\label{table simple compare}
 	\begin{tabular}{|c|c|c|l|c|c|}
 		\hline
 	     &  GD	&       C-CD          &              R-CD         &   RP-CD (Theorem \ref{Thm 3})  & RP-CD (conjectured)  \\
 		\hline
 	   Computation Complexity   &    $ n^3 \kappa_{\text{CD}} $     &   $ n^4 \kappa_{\text{CD}} $             &           $ n^2 \kappa_{\text{CD}}  $   &   
 		$ n^3 \kappa_{\text{CD}}  $ &      $ n^2 \kappa_{\text{CD}}  $      \\
 		\hline
 	\end{tabular}
 \end{table}
%


The following proposition generalizes Theorem \ref{Thm 3} to the non-identity-diagonal case, i.e., $A_i^T A_i$
does not need to be an identity matrix. 
\begin{prop}\label{prop of rate of general RP-BCD}$\;${(rate of RP-BCD for quadratic functions, with non-identity blocks)} $ \text{}$	
	Assume the coefficient matrix $A = [A_1,\dots, A_n]$ is a non-singular square matrix. 
	Suppose RP-BCD is used to solve problem \eqref{quadratic optimization formulation}.
	Denote $D = \text{diag}( A_1^T A_1 , \dots, A_n^T A_n )$ as a block-diagonal matrix, and the norm
	$ \| z\|_D = \sqrt{z^T D z}$. Then
	\begin{equation}\label{CD rate, general}
	\|  E(x^k)- x^*  \|_D \leq  \max \left\{ 1 - \frac{ 1}{n} \lambda_{\min}( D^{1/2} A^T  A D^{-1/2} )  ,   \frac{1}{3} \right\}^k \| x^0 - x^*\|_D .
	\end{equation}
\end{prop}
The proof of Proposition \ref{prop of rate of general RP-BCD} is given in Section \ref{proof of Prop of rate of RP BCD}.
One can  easily transform the  quantity $ \lambda_{\min}( D^{1/2} A^T A D^{-1/2} )$ to certain quantity that only depends
on the eigenvalues of $A_i^T A_i$ and $A^T A$. However, as noted in \cite{sun2016worst}, it is far from clear how
tight the transformation is, thus we skip the transformation here. In fact, it is related to some open question on the so-called Jacobi-preconditioning. We refer the interested readers to \cite{sun2016worst} for a detailed discussion of the subtle issues in the non-identity-diagonal case. 

At last, we present a result on the expected convergence rate of RP-ADMM for solving linear systems, under the assumption that $A_i^T A_i = I, \; \forall i$. Very similar to Proposition \ref{prop of rate of general RP-BCD}, we can also generalize this result to non-identity-diagonal case, i.e., $A_i^T A_i \neq I$, but to save space we skip the generalization here.  
The proof of Theorem \ref{Thm 4} is given in Section \ref{subsec: thm 4 proof}.
\begin{thm}\label{Thm 4}{(Expected convergence rate of RP-ADMM for linear systems)}
	$ \text{}$	
	Assume the coefficient matrix $A = [A_1,\dots, A_n]$ of the constraint in \eqref{optimization formulation} is a non-singular square matrix and $A_i^T A_i = I_{d_i}$.
	Suppose Algorithm \ref{Algorithm: n-block RP-ADMM} is used to solve problem \eqref{optimization formulation}.
	Denote $y^* =  \begin{bmatrix}
	A^{-1}b  \\
	0 \\
	\end{bmatrix}$ as the unique primal-dual optimal solution to the problem \eqref{optimization formulation},
	then
	\begin{equation}\label{convergence to zero, ADMM rate}
	 \| E(y^k) - y^* \| \leq  \left( 1 - \frac{ 1}{2n} \lambda_{\min}(AA^T)  \right) ^k  \| y^0 - y^* \|. 
	\end{equation}
\end{thm}

This result implies that similar to RP-CD for solving quadratic problems, the complexity of RP-ADMM in terms of the expected iterates for solving linear systems is also at most $$ T_{\text{RP-ADMM}} = O( n^3 \kappa_{\text{CD} } \log(1/\epsilon) ) .$$ 
In light of the fact that C-CD has been shown to only achieve a rate $O(n^4 \kappa_{\text{CD} } \log(1/\epsilon))$ \cite{sun2016worst}, the rate of RP-ADMM we obtain is already quite good. 
 Nevertheless, we conjecture that this complexity upper bound can be improved to $O(n^2 \kappa_{\text{CD}} \log(1/\epsilon) )$, the same as the conjectured complexity  for RP-CD.
But an improved rate of RP-ADMM leads to an improved rate of RP-BCD (this should be clear via the comparison of \eqref{I - QAA eigenvalues upp bound} and \eqref{rho M bound}), thus proving this conjecture is an even more difficult problem
than the long-standing open question on RP-CD.

\subsection{Matrix AM-GM Inequality}\label{sec: discussion of AM GM inequality}

To analyze the convergence rate of randomly permuted algorithms, one major technical challenge is matrix AM-GM (algebraic mean-geometric mean) inequality.  
The following conjecture of matrix AM-GM inequality was proposed in \cite{recht2012beneath}: for any positive semi-definite matrix $A_1, \dots, A_n \in \mathbb{R}^{n  \times n} $, 
\begin{equation}\label{Recht AMGM}
\|	\frac{1}{n!}  \sum_{\sigma = (\sigma_1 ,\dots, \sigma_n) \in \Gamma}    A_{\sigma_n} A_{\sigma_{n-1}} \dots  A_{\sigma_1} \|
\leq  \left\|  \left( \frac{1}{n} \sum_i A_i \right)^n \right\|.  
\end{equation}
The original version is more general: the number of matrices does not need to be the same as the dimension of the matrix. For simplicity, we just present a simpler version here.

The matrix AM-GM inequality is a generalization of the well-known AM-GM inequality: for non-negative numbers $a_1, \dots, a_n$, the geometric mean $(a_1 a_2 \dots a_n)^{1/n}$ is no more than the algebraic mean $ \frac{1}{n} \sum_{i=1}^n a_i $.
When extending this inequality to matrix domain, the non-commutative nature of matrix multiplication makes the problem rather difficult to prove.  

We observe that we only need to prove a matrix AM-GM inequality for projection matrices.
We conjecture that the following matrix AM-GM inequality holds. 
\begin{conjecture}\label{new AMGM ineq}$\;$(matrix AM-GM inequality for projection matrices)$\;$
	Suppose $P_i \in \dR^{N \times N}, i=1,\dots, n $ are projection matrices, then
	\begin{equation}\label{conjecture matrix AM-GM ineq}
	\frac{1}{n!}  \sum_{\sigma = (\sigma_1 ,\dots, \sigma_n) \in \Gamma}    P_{\sigma_n} P_{\sigma_{n-1}} \dots  P_{\sigma_1}  \preceq  \left( \frac{1}{n} \sum_i P_i \right)^n. 
	\end{equation}
\end{conjecture}
Compared with \eqref{Recht AMGM}, our conjecture makes a stronger claim on the relation, but it only applies to projection matrices. 
We have found examples to show that \eqref{conjecture matrix AM-GM ineq} does not hold for general positive semi-definite matrices, but it holds for projection matrices in all of our experiments. 

We are not able to prove the new conjecture -- that would solve the open question 
of the best convergence rate of RP-CD for quadratic problem. Nevertheless,
inspired by the new conjecture, we prove a weaker version (see Lemma \ref{lemma of weak AMGM ineq}), which can lead to an improved convergence rate estimate for RP-CD.



\section{Proof of Main Results}
\subsection{Proof of Theorem \ref{Thm 1}}\label{sec: proof of Theorem 1}
Denote $\sigma_k$ as the permutation used in round $k$, and define $\xi_k$ as in \eqref{def of xi}.
Rewrite the update equation \eqref{update eqn of n-block-RADMM} below (replacing $\sigma$ by $\sigma_k$):
\begin{equation}\label{update eqn, with random index}
   y^{k+1} = \bar{L}_{\sigma_k}^{-1} \bar{R}_{\sigma_k} y^{k} + \bar{L}_{\sigma_k}^{-1} \bar{ b }.
\end{equation}

We first prove \eqref{convergence to zero} for the case $b = 0$. By \eqref{barL def, n-block} we have $\bar{b} = 0$, then \eqref{update eqn, with random index} is simplified to $ y^{k+1} = \bar{L}_{\sigma_k}^{-1} \bar{R}_{\sigma_k} y^{k}$.
Taking the expectation of both sides of this equation in $\xi_k$ (see its definition in \eqref{def of xi}),  and note that
$y^{k}$ is independent of $\sigma_k$, we get
$$
  \phi^{k+1} = E_{\xi_k }( \bar{L}_{\sigma_k}^{-1} \bar{R}_{\sigma_k} y^{k}) = E_{\sigma_k } \left( E_{\xi_{k-1}} ( \bar{L}_{\sigma_k}^{-1} \bar{R}_{\sigma_k} y^k )
  \right)
  = E_{\sigma_k }( \bar{L}_{\sigma_k}^{-1} \bar{R}_{\sigma_k} \phi^k ) = M \phi^{k}.
$$
Since the spectral radius of $M$ is less than 1 by Theorem \ref{Thm 2}, we have that $\{ \phi^k\} \rightarrow 0$, i.e. \eqref{convergence to zero}.

We then prove \eqref{convergence to zero} for general $b$.
Let $y^*=[A^{-1}b;0]$ denote the optimal solution. Then it is easy to verify that
$$
y^{*} = \bar{L}_{\sigma_k}^{-1} \bar{R}_{\sigma_k} y^{*} + \bar{L}_{\sigma_k}^{-1} \bar{ b }
$$ for all $\sigma_k \in \Gamma$ (i.e. the optimal solution is the fixed point of the update equation for any order). Compute the difference between this equation and \eqref{update eqn, with random index}
 and letting $\hat{y}^{k}=y^k-y^*$ , we get
$ \hat{y}^{k+1}=\bar{L}_{\sigma_k}^{-1} \bar{R}_{\sigma_k} \hat{y}^{k}$.
According to the proof for the case $b = 0$, we have
$ E(\hat{y}^k ) \longrightarrow 0 $, which implies $E(y^k)  \longrightarrow y^*$.

\subsection{ Proof of Theorem \ref{Thm 2} }\label{sec: proof of Thm 2}  
The difficulty of proving Theorem \ref{Thm 2} (bounding the spectral radius of $M$ defined in \eqref{M definition, theorem}) is two-fold. First, $M$ is a non-symmetric matrix, and there are very few tools to bound the spectral radius of a non-symmetric matrix.
In fact, spectral radius is neither subadditive nor submultiplicative (see, e.g. Kittaneh \cite{kittaneh2006spectral}).
Note that the spectral norm of $M$ can be much larger than $1$ (there are examples that $\| M\| > 2$),
thus we cannot bound the spectral radius simply by the spectral norm.
Second, although it is possible to explicitly write each entry of $M$ as a function of the entries of $A^T A$,
these functions are very complicated ($n$-th order polynomials) and it is not clear how to utilize this explicit expression.

The proof outline of Theorem \ref{Thm 2} and the main techniques are described below.
In Step 0, we provide an expression of the expected update matrix $M$.
 In Step 1, we establish the relationship
between the eigenvalues of $M$ and the eigenvalues of a simple symmetric matrix $A Q A^T $, where $Q$ is defined  in \eqref{Q first def}.
As a consequence, the spectral radius of $M$ is smaller than one iff the eigenvalues of $A Q A^T $ lie in the region $(0, 4/3)$.
   This step partially resolves the first difficulty, i.e. how to deal with the spectral radius of a non-symmetric matrix.
In Step 2, we show that the eigenvalues of $A QA^T $ do lie in $(0, 4/3)$ using mathematical induction.
The induction analysis circumvents the second difficulty, i.e. how to utilize the relation between $M$ and $A$.

Step 0: compute the expression of the expected update matrix $M$.
Define
\begin{equation}\label{Q first def}
Q \triangleq E_{\sigma}(L_{\sigma}^{-1} ) = \frac{1}{n!}\sum_{\sigma \in \Gamma} L_{\sigma}^{-1}.
\end{equation}
It is easy to prove that $Q$ defined by \eqref{Q first def} is symmetric. In fact, note that $L_{\sigma}^T = L_{\bar{\sigma} }, \forall \sigma \in \Gamma$, where $\bar{\sigma}$ is a reverse permutation of $\sigma$ satisfying $\bar{\sigma}(i) = \sigma(n+1 - i), \forall \ i$,
thus $ Q = \frac{1}{n!}\sum_{\sigma} Q_{\sigma} = ( \frac{1}{n!}  \sum_{\sigma} Q_{ \bar{\sigma}} )^T = Q^T , $
where the last step is because the sum of all $Q_{ \bar{\sigma} }$ is the same as the sum of all $Q_{\sigma}$.

Denote
\begin{equation}
M_{\sigma}  \triangleq \bar{L}_{\sigma}^{-1} \bar{R}_{\sigma} =
 \bar{L}_{\sigma}^{-1}
   \begin{bmatrix}
  R_{\sigma}       & A^T   \\
   0      &  I \\
   \end{bmatrix}.
\end{equation}
Substituting the expression of $\bar{L}_{\sigma}^{-1}$ into the above relation, and replacing $R_{\sigma}$ by $L_{\sigma} - A^T A $, we obtain
\begin{equation}\label{Msigma expression}
M_{\sigma} \overset{ }{=}
\begin{bmatrix}
       L_{\sigma}^{-1}  &  0   \\
       -A L_{\sigma}^{-1}   & I \\
   \end{bmatrix}
   \begin{bmatrix}
 L_{\sigma} - A^T A       & A^T   \\
   0      &  I \\
   \end{bmatrix}
   =
\begin{bmatrix}
  I - L_{\sigma}^{-1} A^T A       & L_{\sigma}^{-1} A^T   \\
   -A + A L_{\sigma}^{-1} A^T A      &  I - A L_{\sigma}^{-1} A^T  \\
   \end{bmatrix}.
\end{equation}
Since $M_{\sigma}$ is linear in $ L_{\sigma}^{-1}$, we have
\begin{equation}\label{M expression, first}
\begin{split}
M = E_{\sigma}( M_{\sigma} )
&  =  \begin{bmatrix}
  I - E_{\sigma}( L_{\sigma}^{-1} ) A^T A       & E_{\sigma}( L_{\sigma}^{-1} ) A^T   \\
   -A + A E_{\sigma}(L_{\sigma}^{-1}) A^T A      &  I - A E_{\sigma}( L_{\sigma}^{-1}) A^T  \\
   \end{bmatrix}   \\
&   = \begin{bmatrix}
  I - Q A^T A       &  Q A^T   \\
   -A + A Q A^T A      &  I - A Q A^T  \\
   \end{bmatrix}.
\end{split}
\end{equation}


Step 1: relate $M$ to a simple symmetric matrix. 
The main result of Step 1 is given below, and the proof of this result is relegated to Section \ref{section: proof of Lemma 1}.
\begin{lemma}\label{lemma 1} $\;$
Suppose $A \in \dR^{ N \times N }$ is non-singular and $ Q \in \dR^{ N \times N }$ is an arbitrary matrix. Define $M \in \dR^{2N \times 2N}$  as
\begin{equation}\label{M expression}
M = \begin{bmatrix}
  I - Q A^T A       &  Q A^T   \\
   -A + A Q A^T A      &  I - A Q A^T  \\
   \end{bmatrix}.
\end{equation}
  Then
  \begin{equation}\label{one-to-one relation of eigen}
     \lambda \in \text{eig}(M)  \Longleftrightarrow   \frac{ (1- \lambda)^2 }{ 1 - 2\lambda } \in \text{eig}(Q A^T A).
  \end{equation}
  Furthermore, when $Q$ is symmetric, we have
  \begin{equation}\label{eig transform}
        \rho( M ) < 1   \Longleftrightarrow   \text{eig}(Q A^T A ) \subseteq (0, \frac{4}{3}).
     \end{equation}
\end{lemma}

Remark: For our problem, the matrix $Q$ as defined by \eqref{Q first def} is symmetric (see the argument after
equation \eqref{Q first def}), thus the relation \eqref{eig transform} indeed holds according to Lemma 1.
For a general non-symmetric $Q$, \eqref{eig transform}  does not need to hold, but the first conclusion \eqref{one-to-one relation of eigen} still holds. 

Step 2: Bound the eigenvalues of $Q A^T A$. The main result of Step 2 is summarized in the following Lemma \ref{lemma 2}.
The proof of Lemma \ref{lemma 2} is given in Section \ref{sec: proof of Lemma 2 for n-block}.
\begin{lemma}\label{lemma 2} $\;$
Suppose $A = [A_1, \dots, A_n] \in \dR^{N \times N}$ is non-singular.
 Define $Q$ as
  \begin{equation}
  Q \triangleq E_{\sigma}(L_{\sigma}^{-1} )= \frac{1}{n ! } \sum_{\sigma \in \Gamma } L_{\sigma}^{-1},
  \end{equation}
  in which $L_{\sigma}$ is defined by \eqref{Lsigma def} and $\Gamma$ is defined by \eqref{Gamma def}.
  Then all eigenvalues of $Q A^T A$ lie in $(0, 4/3)$, i.e.
    \begin{equation}\label{QAA eigen lies in (0, 4/3)}
         \text{eig}(Q A^T A ) \subseteq (0, \frac{4}{3}).
    \end{equation}
\end{lemma}

{\black Remark: The upper bound $\frac{4}{3}$ in \eqref{QAA eigen lies in (0, 4/3)} is probably tight, since we have found numerical examples with
 $ \text{eig}(Q A^T A ) > 1.3333 $.
Now the expected convergence of RP-ADMM seems to be a pleasant coincidence:
Lemma \ref{lemma 1} shows that to prove the expected convergence we need to prove $ \sup_{A} \text{eig}(Q A^T A ) $,
a quantity that can be defined without knowing ADMM, is bounded by $4/3$; Lemma \ref{lemma 2} and numerical experiments
show that this quantity happens to be exactly $4/3$ so that RP-ADMM can converge (in expectation).}

Theorem \ref{Thm 2} follows immediately from Lemma \ref{lemma 1} and Lemma \ref{lemma 2}.

\subsection{Proof of Theorem \ref{Thm 3}}\label{subsec: thm 3 proof}

We first describe the outline of the proof. 
 The expected update matrix of RP-BCD is $I - QA^T A$, and the eigenvalues of this matrix lie in $(-1,1)$. 
 The expected convergence speed of RP-BCD depends on the distance between the eigenvalues and the two extremes $-1$ and $1$.
 Lemma \ref{lemma 2} shows that the distance to $-1$ is at least $1/3$, which is a constant.
 We will show that the distance to $1$ is at least $\lambda_{\min}(A^T A)/n$, by proving a weaker version of matrix AM-GM inequality. 
 Combining the two results, we obtain the expected convergence speed of RP-BCD.

The formal proof is presented below.

According to \eqref{update eqn of n-block-RP-CD}, we have $ x^{k+1} - x^* = ( I - L_{\sigma}^{-1}A^T A ) (x^k - x^* ) $, where $\sigma$ is the randomly picked permutation at the $k$-th epoch. 
Therefore, the expected update formula of RP-BCD for solving the least squares problem is
\begin{equation}\label{RP-BCD expected update}
    E(x^{k+1}) - x^* = ( I - Q A^T A ) ( E(x^k) - x^*) .
\end{equation}
It implies
\begin{equation}\label{linear rate of BCD on radius}
  \| E(x^{k+1}) - x^*  \| \leq \rho( I - Q A^T A ) \| E(x^k) - x^* \|. 
\end{equation}
%
%

Suppose the eigenvalues of $ Q A^T A$ are $\eta_1 \geq \eta_2 \geq  \dots \geq  \eta_n$, then according to 
Lemma \ref{lemma 2}, $$ 4/3 > \eta_1  > \dots > \eta_n > 0. $$ 
The eigenvalues of $I -Q A^T A  $ are 
$$
-\frac{1}{3} < 1- \eta_1 \leq \dots \leq 1 - \eta_n < 1,
$$
thus the  spectral radius of $I - QA^T A$ is
\begin{equation}\label{I - QAA eigenvalues upp bound}
\rho(I - QA^T A) = \max\{   1 - \eta_n,  |1 - \eta_1|    \} \leq  \max\{   1 - \eta_n,  \frac{1}{3}   \} 
= \max\{ \lambda_{\max}(I - QA^T A),   \frac{1}{3} \}.
\end{equation}
An interesting phenomenon occurs here. The spectral radius is either $ 1 - \eta_n$ or $|1- \eta_1|$. In the latter case, $\rho(I - QA^T A) = |1 - \eta_1| \leq 1/3$, implying that $\| E(x^k) - x^*\| \leq \frac{1}{3^k} \| E(x^0)- x^* \|$, or equivalently, the relative error $ | E(x^k) - x^*\|/| E(x^0) - x^*\| $ achieves $\epsilon $ in $ \log 3 \log(1/\epsilon) $ epochs. 
We do not even need to  compute $\eta_1$ since it will only affect the convergence speed when the speed is already very fast. From a theoretical perspective, the improvement from $ \log 3 $  to $\log( 1/(1 - |1 - \eta_1|) )  $ is just an improvment in the constant. Therefore, it is reasonable to ignore $\eta_1$ and focus on the estimate of $1 - \eta_n$. 

To estimate the maximum eigenvalue of $I - QA^T A$ (or equivalently, that of $I - A Q A^T$), we first provide a useful identity that connects  $I - A Q A^T$ and projection matrices $P_i = I - A_i A_i^T $.
\begin{claim}\label{lemma of projection expression}  $  \;$
Suppose $A = [A_1,\dots, A_n]$ is a non-singular square matrix, and
$
A_i^T A_i = I, \forall \; i. 
$
	For a permutation $\sigma = (\sigma_1, \dots, \sigma_n) \in \Gamma$,  $L_{\sigma}$ is defined as in \eqref{Lsigma def},   and $Q_{\sigma} = L_{\sigma}^{-1}$.
	Denote $P_i = I - A_i A_i^T $, $i=1, \dots, n$.
	Then we have
	\begin{subequations}
		\begin{align}
			I - AQ_{\sigma} A^T = P_{\sigma_n} P_{\sigma_{n-1}} \dots  P_{\sigma_1}  ,  \label{AQA' and projection} \\
			I - AQA^T = \frac{1}{n!} \sum_{\sigma = (\sigma_1 ,\dots, \sigma_n) \in \Gamma}    P_{\sigma_n} P_{\sigma_{n-1}} \dots  P_{\sigma_1}.  \label{expected AQA' and projection}
		\end{align}
	\end{subequations}	
\end{claim}

The proof of  Claim \ref{lemma of projection expression}  is given at the end of this subsection. 
 Claim \ref{lemma of projection expression} states that $I - AQA^T$ is exactly equal to $ \frac{1}{n!}  \sum_{\sigma = (\sigma_1 ,\dots, \sigma_n) \in \Gamma}    P_{\sigma_n} P_{\sigma_{n-1}} \dots  P_{\sigma_1} $, thus we only need to estimate the maximal eigenvalue of the latter expression. This is achieved by the following lemma (the proof 
 is given in Section \ref{subsec: proof of AMGM weak}).

\begin{lemma}\label{lemma of weak AMGM ineq}$\;$(weak matrix AM-GM inequality)$\;$
	Suppose $P_i \in \dR^{N \times N}, i=1,\dots, n $ are projection matrices, then
	\begin{equation}\label{weak matrix AM-GM ineq}
	\frac{1}{n!}  \sum_{\sigma = (\sigma_1 ,\dots, \sigma_n) \in \Gamma}    P_{\sigma_n} P_{\sigma_{n-1}} \dots  P_{\sigma_1}  \preceq   \frac{1}{n} \sum_i P_i. 
	\end{equation}
\end{lemma}

The above Lemma \ref{lemma of weak AMGM ineq} and Claim \ref{lemma of projection expression} immediately lead to the following corollary.

\begin{coro}$\;$
Suppose $A = [A_1,\dots, A_n]$ is a non-singular square matrix, and
$
A_i^T A_i = I, \forall \; i.
$
Suppose $P_i = I - A_i A_i^T, \; \forall \; i.$
$L_{\sigma}$ is defined as in \eqref{Lsigma def}, and	$Q = E_{\sigma} (L_{\sigma}^{-1}) $.
Then 
	\begin{equation}\label{I - AQA upper bound}
  I - AQA^T  \preceq   \frac{1}{n} \sum_i P_i. 
\end{equation}
\end{coro}

Note that $ \frac{1}{n} \sum_i P_i = \frac{1}{n} ( nI - \sum_i A_i A_i^T ) = I - \frac{1}{n} AA^T $,
thus \eqref{I - AQA upper bound} implies
$$
  I - AQA^T \preceq I - \frac{1}{n} AA^T,
$$
which implies
 \begin{equation}\label{QAA lam min bound}
   \lambda_{\max}(I - AQA^T ) \leq 1 - \frac{ 1}{n} \lambda_{\min}(AA^T) . 
 \end{equation} 
Substituting into \eqref{I - QAA eigenvalues upp bound}, we get
$$
\rho(I - QA^T A) \leq   \max\{ \lambda_{\max}(I - QA^T A),   \frac{1}{3} \}
\leq \max\{ 1 - \frac{ 1}{n} \lambda_{\min}(AA^T), 1/3  \}. 
$$
Substituting this relation into  \eqref{linear rate of BCD on radius}, we obatain
$$
\| E(x^{k+1}) - x^*  \| \leq   \max\{ 1 - \frac{ 1}{n} \lambda_{\min}(AA^T)  ,   \frac{1}{3} \}   \| E(x^k) - x^* \|. 
$$  \textbf{Q.E.D.}

Remark: There is a coefficient $1/n$ in front of $ \lambda_{\min}(AA^T) $ in \eqref{QAA lam min bound}, and this is why
the complexity of RP-CD we establish is $n$ times worse than the conjectured one in Table \ref{table simple compare}. If Conjecture \ref{new AMGM ineq} holds, then this factor of $1/n$ would be removed
and the conjectured (expected) complexity of RP-CD in Table \ref{table simple compare} would hold.

\subsubsection{Proof of Claim \ref{lemma of projection expression}}\label{subsubsec: proof of dual expression}
We prove \eqref{AQA' and projection} by induction on $n$.
Without loss of generality, we can assume $\sigma = (1,2,\dots , n)$, then
$L_{\sigma} =
\begin{bmatrix}
A_1^T A_1 &    0   &    \dots &    0       \\
A_2^T A_1 & A_2^T A_2 &     \dots &   0      \\
\vdots    &  \vdots       &      \ddots &  \vdots   \\
A_n^T A_1 & A_n^T A_2 &   \dots & A_n^T A_n     \\
\end{bmatrix}.  $
In this case,  \eqref{AQA' and projection} becomes
$$
I - AL_{\sigma}^{-1} A^T  =  P_n P_{n-1} \dots P_1 .
$$

The expression obviously holds for $n = 1$.
Suppose the expression holds for $n-1$, i.e., for $\hat{A} = [A_1 , \dots, A_{n-1}]$,
we have
\begin{equation}\label{P1P2...Pn-1 induction hypothesis}
\hat{Z} \triangleq I - \hat{A} \hat{L}_{\hat{\sigma} }^{-1} \hat{A}^T =  P_{n-1} \dots P_2 P_1,
\end{equation}
where $\hat{\sigma}= (1,2,\dots, n-1)$ is a permutation of $n-1$ elements and
$ \hat{L}_{\hat{\sigma}} $ is the counterpart of $L_{\sigma}$ for $n-1$ blocks defined as
$$
\hat{L}_{ \hat{\sigma} } =
\begin{bmatrix}
A_1^T A_1 &    0   &    \dots &    0       \\
A_2^T A_1 & A_2^T A_2 &     \dots &   0      \\
\vdots    &  \vdots       &      \ddots &  \vdots   \\
A_{n-1}^T A_1 & A_{n-1}^T A_2 &   \dots & A_{n-1}^T A_{n-1}    \\
\end{bmatrix}. 
$$

The two matrices $L_{\sigma}$ and $\hat{L}_{\sigma'}$ are related by
$$
L_{\sigma} = \begin{bmatrix}
\hat{L}_{\hat{\sigma}}  &  0  \\
A_n^T \hat{A}    &     I   \\
\end{bmatrix}  ,
$$
which implies
$$
L_{\sigma}^{-1} = 
\begin{bmatrix}
\hat{L}_{\hat{\sigma}}^{-1}  &  0  \\
-   A_n^T \hat{A}  \hat{L}_{\hat{\sigma}}^{-1}  &     I   \\
\end{bmatrix}  .
$$
Therefore we have
\begin{align*}
A L_{\sigma}^{-1} A^T  =
[\hat{A}, A_n ] \begin{bmatrix}
\hat{L}_{\hat{\sigma}}^{-1}  &  0  \\
-   A_n^T \hat{A}  \hat{L}_{\hat{\sigma}}^{-1}  &     I   \\
\end{bmatrix}   [\hat{A}, A_n ]^T 
&  =  \hat{A} \hat{L}_{\hat{\sigma} }^{-1} \hat{A}^T   
- A_n  A_n^T \hat{A}  \hat{L}_{\hat{\sigma}}^{-1}  \hat{A}^T + A_n A_n^T \\
& =     \hat{Z} - A_n A_n^T \hat{Z} + A_n A_n^T  \\
& = I -   (I - A_n A_n^T ) (I - \hat{Z})   \\
& = I - P_n P_{n-1}  \dots P_1,
\end{align*}
where in the last step we use the induction hypothesis \eqref{P1P2...Pn-1 induction hypothesis}. 
Thus we have proved \eqref{AQA' and projection}.
Summing up \eqref{AQA' and projection} for all possible permutations $\sigma$ and divide by $n!$, we obtain \eqref{expected AQA' and projection}.   $\Box$

\subsection{Proof of Proposition \ref{prop of rate of general RP-BCD} }\label{proof of Prop of rate of RP BCD}
According to \eqref{RP-BCD expected update}, the (expected) update equation of RP-BCD is given by
$ E(x^{k+1}) - x^* = ( I - Q A^T A ) ( E(x^k) - x^*)  = Z( E(x^k) - x^*) $,
where $Z = I - QA^T A = I - E(L_{\sigma}^{-1} A^T A) $.

Consider a new coefficient matrix $\tilde{A} = [ \tilde{A}_1, \dots, \tilde{A}_n ] $ where $\tilde{A}_i = A_i (A_i^T A_i)^{-\frac{1}{2}}$. Clearly $\tilde{A}_i^T \tilde{A}_i = I_{d_i}$. 
Denote the corresponding matrices as $\tilde{L}_{\sigma}, \tilde{Z} .$
Define $
\Lambda \triangleq  \text{Diag}( (A_1^T A_1)^{\frac{1}{2}}, \dots, (A_n^T A_n)^{\frac{1}{2}} ) = D^{1/2}. $
When $\sigma = (1, 2,\dots, n)$, we have
$$ L_{\sigma} =
\begin{bmatrix}
A_1^T A_1 &    0   &    \dots &    0       \\
A_2^T A_1 & A_2^T A_2 &     \dots &   0      \\
\vdots    &  \vdots       &      \ddots &  \vdots   \\
A_n^T A_1 & A_n^T A_2 &   \dots & A_n^T A_n     \\
\end{bmatrix}  , \quad 
\tilde{L}_{\sigma}  =   
\begin{bmatrix}
\tilde{A}_1^T \tilde{A}_1 &    0   &    \dots &    0       \\
\tilde{A}_2^T \tilde{A}_1 &  \tilde{A}_2^T \tilde{A}_2 &     \dots &   0      \\
\vdots    &  \vdots       &      \ddots &  \vdots   \\
\tilde{A}_n^T \tilde{A}_1 & \tilde{A}_n^T \tilde{A}_2 &   \dots & \tilde{A}_n^T \tilde{A}_n     \\
\end{bmatrix}
=  \Lambda^{-1} L_{\sigma} \Lambda^{-1}. 
$$
It is not hard to verify that the above relation $ \tilde{L}_{\sigma}  = \Lambda^{-1} L_{\sigma} \Lambda^{-1} $ is true for any $\sigma$.
Similarly, we have $\tilde{A}^T \tilde{A} = \Lambda^{-1} A^T A \Lambda^{-1} $, thus
$$
\tilde{L}_{\sigma}^{-1} \tilde{A}^T \tilde{A} = \Lambda L_{\sigma}^{-1} \Lambda \Lambda^{-1} A^T A \Lambda^{-1}=  \Lambda L_{\sigma}^{-1} A^T A \Lambda^{-1}.
$$
This implies
$$
\tilde{Z} = E( I - \tilde{L}_{\sigma}^{-1} \tilde{A}^T \tilde{A}  ) = \Lambda  (I - E(L_{\sigma}^{-1} A^T A) ) \Lambda^{-1} = \Lambda Z \Lambda^{-1}.
$$
Consider a sequence $\tilde{x}^k = \Lambda x^k  $ and define $\tilde{x}^* = \Lambda x^*$.
Then from the original update equation we have
$  \Lambda^{-1} ( E( \tilde{x}^{k+1}) - \tilde{x}^* ) =   Z  \Lambda^{-1}   ( E( \tilde{x} ^{k}) - \tilde{x}^* ) $,
i.e., 
$$
E( \tilde{x}^{k+1}) - \tilde{x}^*  = \Lambda Z \Lambda^{-1}( E( \tilde{x} ^{k}) - \tilde{x}^* ) = \tilde{Z} ( E( \tilde{x} ^{k}) - \tilde{x}^* ) .
$$
According to Theorem \ref{Thm 3}, we have
\begin{equation}\label{tilde result rate}
\| E( \tilde{x}^{k}) - \tilde{x}^*\| \leq  \left\{ 1 - \frac{ 1}{n} \lambda_{\min}( \tilde{A}^T  \tilde{A} )  ,   \frac{1}{3} \right\}^k \| \tilde{x}^0 - \tilde{x}^* \|. 
\end{equation}
Note that  $ \| E( \tilde{x}^{k}) - \tilde{x}^*\| = \|  \Lambda ( E(x^k) - x^* )  \| = \sqrt{ (E(x^k) - x^*)^T \Lambda^2  E(x^k) - x^*}
= \| E(x^k) - x^*  \|_D $, and 
$ \tilde{A}^T  \tilde{A} =  \Lambda^{-1} A^T A \Lambda^{-1} =   D^{-1/2} A^T \tilde{A}  D^{-1/2}  $.
Substituting into \eqref{tilde result rate}, we obtain the desired inequality. 

\subsection{Proof of Theorem \ref{Thm 4}}\label{subsec: thm 4 proof}



Now we consider the expected convergence rate of RP-ADMM. 
The difference with the analysis for RP-BCD is that here we need to consider the distance between the eigenvalues
of $I - AQA^T$ with $-1/3$ while for RP-BCD what matters is the distance between the eigenvalues of $I - AQA^T$
and $-1$ which is at least $2/3$ and thus can be ignored.  

\begin{claim}$\;$\label{claim: RP-ADMM spectral radius}
Suppose the minimum and maximum eigenvalues of $QA^T A$ are $0< \tau_{\min} \leq \tau_{\max} < 4/3$.
Then
$$
 \rho(M) = \max\left\{ \sqrt{ (1 - \tau_{\min})_+ }  , ( \tau_{\max} - 1 )_+ + \sqrt{ \tau_{\max} (\tau_{\max} -1 )_+ }  \right\} ,
$$
where $z_+ = max\{ z, 0\}$. Furthermore, we have
 \begin{equation}\label{rho M bound}
 \rho(M) \leq \max\left\{ 1 - \frac{3}{4} (4 - 3 \tau_{\max}) , \; \; 1 - \frac{1}{2}\tau_{\min} \right\}. 
 \end{equation}
\end{claim}
The proof of Claim \ref{claim: RP-ADMM spectral radius} is given in Section  \ref{subsec: proof of Claim of spec radius}.
The next lemma provides a universal estimate of the maximum eigenvalules of $QA^T A$. 
 
 \begin{lemma}\label{lemma of stronger bound} $\;$
 	The maximum eigenvalues of $QA^T A $ is at most $ \frac{4}{3} -  \frac{4}{9} \frac{1}{ n+1 } $, i.e.,
 \begin{equation}\label{tau max bound}
 	  \tau_{\max} = \lambda_{\max}( QA^T A ) \leq \frac{4}{3} -  \frac{4}{9} \frac{1}{ n+1 } .
 \end{equation}
 \end{lemma}
The proof of Lemma \ref{lemma of stronger bound} is given in Section \ref{subsec: lemma 5 proof}
 	
According to \eqref{QAA lam min bound}, which is established in the proof of the expected convergence rate of RP-BCD, we have
\begin{equation}\label{tau min bound}
  \tau_{\min} =  \lambda_{\min}(QA^T A) \geq \frac{1}{n} \lambda_{\min}(A^T A) .
\end{equation}

Substituting the bounds \eqref{tau max bound} and \eqref{tau min bound} into \eqref{rho M bound}, we obtain
\begin{equation}\label{final rho(M) bound}
 \rho(M) \leq \max\left\{ 1 - \frac{3}{4} (4 - 3 \tau_{\max}) , \; \; 1 - \frac{1}{2}\tau_{\min} \right\}
  = \max \left\{ 1 - \frac{1}{n+1},  1 -  \frac{1}{2n}  \lambda_{\min}(A^TA)   \right\}  .
\end{equation}
Since $  \lambda_{\min}(A^TA)  \leq 1 $, $\frac{1}{2n} \leq \frac{1}{n+1}$, this bound can be simplified to 
$$
   \rho(M) \leq 1 -  \frac{1}{2n}  \lambda_{\min}(A^TA) .  \quad \quad \quad   \textbf{Q.E.D.}
$$ 	

 Remark: The eigenvalues of $QA^T A$ lie in the region $(0,4/3)$, which guarantees the expected convergence of RP-ADMM.
To obtain the expected convergence rate, we need to know the distance of the spectrum to the two extremes $0 $ and $4/3 $. We conjecture that the bound can be improved to $\rho(M) \leq 1 - \frac{1}{2} \lambda_{\min}(A^T A) $.
This requires more effort than the conjecture of RP-CD: besides showing $ \tau_{\min} \geq O( \lambda_{\min}(A^T A)  ) ,$
we also need to show $ \tau_{\max} \leq \frac{4}{3} - O( \lambda_{\min}(A^T A)  )  $. This is left as future work.



\section{Proof of Lemma \ref{lemma 1} }\label{section: proof of Lemma 1}

The proof of Lemma \ref{lemma 1} relies on two simple techniques.
The first technique, as elaborated in the Step 1 below, is to factorize $M$ and rearrange the factors.
The second technique, as elaborated in the Step 2 below, is to reduce the dimension by eliminating a variable from the eigenvalue equation.

\textbf{Step 1}: Factorizing $M$ and rearranging the order of multiplication. 
The following observation is crucial: the matrix $M$ defined by \eqref{M expression} can be factorized as
\begin{equation}\nonumber
M = \begin{bmatrix}
  I        &  0   \\
   -A      &  I  \\
   \end{bmatrix}
   \begin{bmatrix}
  Q A^T        &  I   \\
   I     &  A  \\
   \end{bmatrix}
   \begin{bmatrix}
  -A      &  I   \\
   I      &  0  \\
   \end{bmatrix}.
\end{equation}
Switching the order of the products by moving the first component to the last, we get a new matrix
\begin{equation}\label{M prime expression}
M^{\prime} \triangleq
   \begin{bmatrix}
  Q A^T        &  I   \\
   I     &  A  \\
   \end{bmatrix}
   \begin{bmatrix}
  -A      &  I   \\
   I      &  0  \\
   \end{bmatrix}
   \begin{bmatrix}
  I        &  0   \\
   -A      &  I  \\
   \end{bmatrix}
  =  \begin{bmatrix}
  Q A^T        &  I   \\
   I     &  A  \\
   \end{bmatrix}
   \begin{bmatrix}
  -2A      &  I   \\
   I      &  0  \\
   \end{bmatrix}
   =
    \begin{bmatrix}
  I - 2 Q A^T A        &  Q A^T   \\
   -A     &  I  \\
   \end{bmatrix}.
\end{equation}
Note that $\text{eig}(XY) = \text{eig}(YX)$ for any two square matrices, thus $$ \text{eig}(M) = \text{eig}(M^{\prime} ) . $$
To prove \eqref{one-to-one relation of eigen}, we only need to prove
 \begin{equation}\label{one-to-one relation of eigen, M'}
     \lambda \in \text{eig}(M')  \Longleftrightarrow   \frac{ (1- \lambda)^2 }{ 1 - 2\lambda } \in \text{eig}(Q A^T A).
  \end{equation}

\textbf{Step 2}: Relate the eigenvalues of $M'$ to the eigenvalues of $QA^T A$, i.e. prove \eqref{one-to-one relation of eigen, M'}.
This step is simple as we only use the definition of eigenvalues. However, note that, without Step 1, just applying the definition of eigenvalues
of the original matrix $M$ may not lead to a simple relationship as \eqref{one-to-one relation of eigen, M'}.

We first prove one direction of \eqref{one-to-one relation of eigen, M'}:
\begin{equation}\label{one direction of M to QAA}
     \lambda \in \text{eig}(M')  \Longrightarrow   \frac{ (1- \lambda)^2 }{ 1 - 2\lambda } \in \text{eig}(Q A^T A).
  \end{equation}
Suppose $v \in \dC^{2N \times 1}\backslash \{ 0\} $ is an eigenvector of $M'$ corresponding to the eigenvalue $\lambda$, i.e.
$$
   M' v = \lambda v.
$$
Partition $v$ as $ v = \begin{bmatrix}  v_1 \\ v_0 \end{bmatrix}   $, where $ v_1, v_0 \in \dC^{ N \times 1} $. Using the expression of $M'$
 in \eqref{M prime expression}, we can write the above equation as
$$
  \begin{bmatrix}
  I - 2 Q A^T A        &  Q A^T   \\
   -A     &  I  \\
   \end{bmatrix}  \begin{bmatrix}  v_1 \\ v_0 \end{bmatrix}
   = \lambda \begin{bmatrix}  v_1 \\ v_0 \end{bmatrix},
$$
which implies
\begin{subequations}\label{eigen equations}
  \begin{align}
     (I - 2 Q A^T A ) v_1 +  Q A^T v_0 = \lambda v_1, \label{eigen equation a} \\
      -A v_1 + v_0  = \lambda v_0.   \label{eigen equation b}
  \end{align}
\end{subequations}

We claim that \eqref{one direction of M to QAA} holds when $v_1 = 0$.
In fact, in this case we must have $v_0 \neq 0$ (otherwise $v = 0$ cannot be an eigenvector). By \eqref{eigen equation b} we have $ \lambda v_0 =  v_0$, thus $\lambda = 1$.
By \eqref{eigen equation a} we have $0 = Q A^T v_0 = Q A^T A (A^{-1} v_0 )$, which implies $\frac{ (1- \lambda)^2 }{ 1 - 2\lambda } = 0 \in \text{eig}(QA^T A)$,
therefore \eqref{one direction of M to QAA} holds in this case.

We then prove  \eqref{one direction of M to QAA} for the case
\begin{equation}\label{v1 not equal 0}
v_1 \neq 0.
\end{equation}
The equation \eqref{eigen equation b} implies  $ ( 1 - \lambda ) v_0 =  A v_1  $. Multiplying both sides of \eqref{eigen equation a} by $(1 - \lambda)$
and invoking this equation,  we get
$$
  (1 - \lambda) (I - 2 Q A^T A ) v_1 +  Q A^T A v_1  = (1 - \lambda ) \lambda v_1.
$$
This relation can be simplified to
\begin{equation}\label{temp equation}
 (1 - 2 \lambda) Q A^T A  v_1 =  (1 - \lambda)^2 v_1.
\end{equation}
We must have $ \lambda \neq \frac{1}{2}$; otherwise, the above relation implies $ v_1 = 0$, which contradicts \eqref{v1 not equal 0}.
Then \eqref{temp equation} becomes
\begin{equation}\label{QAA eigen equation}
   Q A^T A  v_1 = \frac{ (1 - \lambda)^2 } { 1 - 2 \lambda } v_1.
\end{equation}
Therefore, $ \frac{ (1 - \lambda)^2 } { 1 - 2 \lambda } $ is an eigenvalue of $ Q A^T A $, with the corresponding eigenvector
$ v_1 \neq 0$, which finishes the proof of \eqref{one direction of M to QAA}.

The other direction \footnote{For the purpose of proving Theorem \ref{Thm 2}, we do not need to prove this direction. Here we present the proof since it is quite straightforward
and makes the result more comprehensive. }
\begin{equation}\label{reverse direction of M to QAA}
     \lambda \in \text{eig}(M)  \Longleftarrow   \frac{ (1- \lambda)^2 }{ 1 - 2\lambda } \in \text{eig}(Q A^T A)
  \end{equation}
is easy to prove. Suppose $ \frac{ (1- \lambda)^2 }{ 1 - 2\lambda } \in \text{eig}(Q A^T A)$. We consider two cases.

Case 1: $ \frac{ (1- \lambda)^2 }{ 1 - 2\lambda } = 0  $. In this case $ \lambda = 1$.
Since $0 = \frac{ (1- \lambda)^2 }{ 1 - 2\lambda } \in \text{eig}(Q A^T A) $, there exists $v_0 \in \dC^N \backslash \{0 \}$ such that $ QA^T A v_0 = 0 $ and
Let $v_1 = (0,\dots, 0)^T \in \dC^{N \times 1} $, then
$v_0, v_1 $ and $\lambda = 1$ satisfy \eqref{eigen equations}. Thus $ v = \begin{bmatrix}  v_1 \\ v_0 \end{bmatrix} \in \dC^{2N} \backslash \{0 \}   $ satisfies
$M v = \lambda v$, which implies $ \lambda = 1 \in \text{eig}(M) $.

Case 2: $ \frac{ (1- \lambda)^2 }{ 1 - 2\lambda } \neq 0  $, then $\lambda \neq 1$.
Let $v_1$ be the eigenvector corresponding to $ \frac{ (1- \lambda)^2 }{ 1 - 2\lambda } $ (i.e. pick $v_1$ that satisfies \eqref{QAA eigen equation}),
and define $v_0 = v_1/(1-\lambda)$. It is easy to verify that $ v = \begin{bmatrix}  v_1 \\ v_0 \end{bmatrix}  $ satisfies
$M v = \lambda v$, which implies $ \lambda \in \text{eig}(M) $.

\textbf{Step 3}: When $Q$ is symmetric, prove \eqref{eig transform} by simple algebraic computation.

Since $Q$ is symmetric, we know that $\text{eig}(QA^T A) = \text{eig}(AQA^T) \subseteq \dR$. 
Suppose $\tau \in \dR $
  is an eigenvalue of $QA^T A$, then
  any $\lambda$ satisfying $  \frac{ (1- \lambda)^2 }{ 1 - 2\lambda } = \tau $
  is an eigenvalue of $M$.
This relation can be rewritten as
$ \lambda^2 + 2(\tau - 1) \lambda + (1 - \tau) = 0$, which, as a real-coefficient quadratic equation in $\lambda$, has two roots 
\begin{equation}\label{two roots}
  \lambda_1 = 1 - \tau + \sqrt{\tau(\tau - 1) }, \quad \lambda_2 = 1 - \tau - \sqrt{\tau(\tau - 1) }.
\end{equation}
Note that when $ \tau(\tau - 1)< 0 $, the expression $\sqrt{\tau(\tau - 1) } $ denotes a complex number $ i \sqrt{ \tau(1 -\tau) }  $,
where $i $ is the imaginary unit.
To prove  \eqref{eig transform}, we only need to prove
\begin{equation}\label{|lam| < 1 <==> f(lam) < 4/3}
  \max\{ |\lambda_1|, |\lambda_2| \} <1 \Longleftrightarrow 0 <  \tau < \frac{4}{3}.
\end{equation}
Consider three cases.

Case 1: $\tau < 0$. Then $\tau(\tau - 1) = |\tau| ( |\tau| + 1) > 0$. In this case,
$ \lambda_1 = 1 + |\tau| + \sqrt{ |\tau| ( |\tau| + 1) } > 1 .$ 

Case 2: $ 0 < \tau < 1$. Then $\tau(\tau - 1) < 0$, and \eqref{two roots} can be rewritten as
$$
  \lambda_{1,2} = 1 - \tau \pm i \sqrt{\tau( 1 -\tau) },
$$
which implies $|\lambda_1| = |\lambda_2| = \sqrt{ (1 - \tau)^2 + \tau(1-\tau) } = \sqrt{1 - \tau} < 1$.

Case 3: $ \tau > 1 $. Then $\tau( \tau - 1) > 0$. According to \eqref{two roots}, it is easy to verify $\lambda_1 > 0 > \lambda_2$ and
$$
  |\lambda_{2}| =  \tau - 1  +  \sqrt{\tau( \tau - 1 ) } >  1 - \tau + \sqrt{\tau(\tau - 1) }  =  | \lambda_1 |.
$$
Then we have
\begin{equation}\nonumber
  \max\{ |\lambda_1|, |\lambda_2| \} <1 \Longleftrightarrow  |\lambda_{2}| =  \tau - 1  +  \sqrt{\tau( \tau - 1 ) } < 1
  \Longleftrightarrow  1 < \tau < \frac{4}{3}.
\end{equation}

Combining the conclusions of the three cases immediately leads to \eqref{|lam| < 1 <==> f(lam) < 4/3}.

\section{ Proof of Lemma \ref{lemma 2} }\label{sec: proof of Lemma 2 for n-block}
This section is devoted to the proof of Lemma \ref{lemma 2}. We first give a proof overview in Section \ref{proof overview}.
The formal proof of Lemma \ref{lemma 2} is given in Section \ref{sec: outline proof of Lemma 2}.
The proofs of the technical results involved in the proof are given in the subsequent subsections.

Without loss of generality, we can assume
$$
  A_i^T A_i = I_{d_i \times d_i}, \ i = 1,\dots, n.
$$
To show this, let us write $M_{\sigma}, M$ as $M_{\sigma}(A_1,\dots, A_n)$ and $M(A_1,\dots, A_n)$ respectively,
i.e. functions of the coefficient matrix $(A_1,\dots, A_n)$.
Define $\tilde{A}_i = A_i (A_i^T A_i)^{-\frac{1}{2}}$ and $$
D \triangleq  \text{Diag}( (A_1^T A_1)^{-\frac{1}{2}}, \dots, (A_n^T A_n)^{-\frac{1}{2}} , I_{N \times N} ) . $$ It is easy to verify that
$
  M_{\sigma}(A_1,\dots, A_n) = D^{-1} M_{\sigma}(\tilde{A}_1,\dots, \tilde{A}_n) D,
$
which implies
$$
  M(A_1,\dots, A_n) = D^{-1} M(\tilde{A}_1,\dots, \tilde{A}_n) D.
$$
Thus $\rho( M(A_1,\dots, A_n) ) = \rho( M(\tilde{A}_1,\dots, \tilde{A}_n) )$.
In other words, normalizing $A_i$ to $\tilde{A}_i$, which satisfies
$\tilde{A}_i^T \tilde{A}_i = I_{d_i \times d_i}$, does not change the spectral radius of $M$.

\subsection{Proof Overview}\label{proof overview}

In the proof overview, we discuss a few issues one may encounter when proving the result, and how we resolve these issues. 

The simulations show that $\|Q A^T A \| < \frac{4}{3} \ll \| Q\| \|A^T A \|$, thus we cannot relax $ \|Q A^T A \|$ to the product of $\| Q\|$ and $\|A^T A \|$, and have to 
 treat $QA^T A$ as a single subject.
However, each entry of $QA^T A$  is a complicated function (in fact, a high order polynomial) of the entries of $A^T A$. In other words, $Q$ is like a black box. 
To open the ``black box'', we use a simple expression of  $Z = I - A QA^T $ proved
in Claim \ref{lemma of projection expression}, i.e., $Z = E_{ \sigma } (P_{\sigma_1} \dots, P_{\sigma_n}), $ where  $P_i = I - A_i A_i^T$ is directly related to $A_i$. The problem becomes how to connect
the eigenvalues of $E_{ \sigma } (P_{\sigma_1} \dots, P_{\sigma_n})$ with those of $AA^T = \sum_i A_i A_i^T = n - \sum_i P_i $. 

Although this is a clear linear algebra problem, it is not easy to obtain a lower bound of $E_{ \sigma } (P_{\sigma_1} \dots, P_{\sigma_n})$. 
In fact, even though we know  the eigenvalues of $Z = E_{ \sigma } (P_{\sigma_1} \dots, P_{\sigma_n}) $ are lower bounded by $-1$ because RP-CD converges, it is not clear how to prove this lower bound directly from a linear algebra perspective.

In our solution, we apply two tricks.
The first trick is to view $E_{ \sigma } (P_{\sigma_1} \dots, P_{\sigma_n})$ as an induction formula that connects it and its lower dimensional analogs.  
This is based on a simple observation that any permutation $ (\sigma_1 \sigma_2 \dots \sigma_n) $ can be written as the concatenation of $(\sigma_1 \sigma_2 \dots \sigma_{n-1})$ and 
$\sigma_n$, thus the expression of $Z =E_{ \sigma } (P_{\sigma_1} \dots, P_{\sigma_n}) $
can be decomposed accordingly. 
We then reduce the problem to bounding the eigenvalues of a Jordan product
$ P_n \hat{Z} + \hat{Z} P_n  $, where $P_n$ is a projection matrix and $\hat{ Z }$ is the lower dimensional analog of $Z$. 
The second trick is to apply a formula on the eigenvalues of Jordan product developed by Strang in 1962 \cite{strang1962eigenvalues}. Somewhat surprisingly, his formula exactly leads to the desired lower bound of $-1/3$.

\subsection{Proof of Lemma \ref{lemma 2}}\label{sec: outline proof of Lemma 2}

The proof can be divided into three steps: first
provide an alternative expression of $AQA^T$, then prove an induction formula, and finally apply Strang's formula to perform mathematical induction. This subsection contains the major part of the proof, and the intermediate technical results  will be proved in later subsections. 

\textbf{Step 0}: Expression of $I - AQ A^T$. 
As proved in Claim \ref{lemma of projection expression}, we have a simple expression of the update matrix $I - AQ A^T$
$$ 	I - AQA^T = \frac{1}{n!} \sum_{\sigma = (\sigma_1 ,\dots, \sigma_n) \in \Gamma}    P_{\sigma_n} P_{\sigma_{n-1}} \dots  P_{\sigma_1}.  
$$


\textbf{Step 1}: Induction formula. 

For any $k \in [n]$, define
\begin{equation}\label{Gamma k def}
\Gamma_k \triangleq \{ \sigma' \mid \sigma' \text{ is a permutation of } [n]\backslash \{ k \} \}.  
\end{equation}
For any $\sigma^{\prime} \in \Gamma_k $,  we define $ L_{\sigma'} \in \dR^{(N-d_k) \times (N-d_k)} $
as a $(n-1)\times(n-1)$ block-partitioned matrix, with the $(\sigma'(i) , \sigma'(j))$-th block being 
\begin{equation}\label{Lsigma' def, n-block}
L_{\sigma'}[ \sigma'(i) , \sigma'(j) ] \triangleq \begin{cases}
A_{\sigma'(i)}^T A_{\sigma'(j) }  &  i \geq j,    \\
0  &  i < j,
\end{cases}
\end{equation}
We then define $\hat{Q}_k \in  \dR^{(N-d_k) \times (N-d_k)}  $  by
\begin{equation}\label{hat(Q) def, n-block}
\hat{Q}_k \triangleq \frac{1}{|\Gamma_k|}\sum_{\sigma' \in \Gamma_k } L_{\sigma'}^{-1}, \ k=1, \dots, n.
\end{equation}
Define $W_k$ as the $k$-th block-column of $A^T A$ excluding the block $A_k^T A_k$, i.e.
\begin{equation}\label{w def, n-block}
\begin{split}
W_k = [ A_k^T A_1 , \dots,  A_k^T A_{k-1} , A_k^T A_{k+1} , \dots, A_k^T A_n ]^T , \ \forall k \in [n].
\end{split}
\end{equation}

Based on the expression of $I- AQA^T$ presented before, we build a connection between the update matrix $I - AQA^T$ and its lower dimensional analogs. 
 The proof of Proposition \ref{prop new: induction formula} is given in Section \ref{appen: proof of induction formula, n-block}.
\begin{prop}\label{prop new: induction formula} $\;$
 Define
	$$   Z = I - AQA^T, \quad \hat{Z}_k = I -  \hat{A}_k  \hat{ Q}_k  \hat{A}_k^T ,  $$
	where $Q$ is defined as in \eqref{Q first def}, $\hat{A}_k =  [A_1, \dots, A_{k-1}, A_{k+1}, \dots, A_n]$,
	and $\hat{Q}_k$  is defined in \eqref{hat(Q) def, n-block}, and $P_k = I - A_k A_k^T$.  Then we have
\begin{equation}\label{Z induction}
Z  = \frac{1}{ 2 n} \sum_{k=1}^n  ( P_k \hat{Z}_k  + \hat{Z}_k  P_k   ). 
\end{equation}
\end{prop}

\textbf{Step 2}: Applying Strang's result on Jordan product to perform mathematical induction. 

 It is obvious that the product of two symmetric matrices is not necessarily symmetric,
 so it is common to encounter the symmetrized product $XY + YX$, which is called Jordan product of two matrices $X$ and $Y$.
Our induction formula basically states that $Z $ is the average of the Jordan product of the lower dimensional analog and $P_k$. 

The eigenvalues of the Jordan product of two matrices have been studied before.
The following result is proved in Strang \cite{strang1962eigenvalues}.

\begin{lemma}\label{Strang lemma}(\cite[Theorem 1]{strang1962eigenvalues}; eigenvalues of Jordan product)
	Suppose two symmetric positive-semidefinite matrices $X$ and $Y$ satisfy 
	$$
	  \alpha_{1} I  \preceq  X \preceq  \alpha_{n} I,  \quad  \beta_{1} I  \preceq  Y \preceq  \beta_{n} I,  
	$$
	then the maximal (resp. minimal) eigenvalue of the Jordan product $XY + YX $ are 
	the largest (resp. smallest) of the set
	\begin{equation}\label{set of values}
		\left\{   2 \alpha_i \beta_j,  i,j \in \{ 1, n\},  \;  \frac{ 16 \alpha_1 \alpha_n \beta_1 \beta_n -
		(\beta_1 - \beta_n)^2 (\alpha_1 - \alpha_n )^2 }{ 4(\alpha_1 + \alpha_n) (\beta_1 + \beta_n ) }    \right\}.
	\end{equation}
	
\end{lemma}

	Let us come back to the proof of Lemma \ref{lemma 2}. We use mathematical induction to prove Lemma \ref{lemma 2}.
	For the basis of the induction ($n=1$), Lemma \ref{lemma 2} holds since $Q A^T A = I_{d_1 \times d_1}$.
	Assume Lemma \ref{lemma 2} holds for $n-1$, we will prove Lemma \ref{lemma 2} for $n$.
	
Consider one term of \eqref{Z induction} $ P_k \hat{Z}_k  + \hat{Z}_k  P_k $.
Note	that $P_k = I- A_k A_k^T$ is a projection matrix, since we have assumed $A_k^T A_k = I$.
Combining with the induction hypothesis, we have
$$
0  \preceq  P_k  \preceq I,   \quad   -\frac{1}{3}I     \prec  \hat{Z}_k  \prec I .
$$
Let $\alpha_1 = 0, \alpha_n = 1, \beta_1 = -1/3, \beta_n = 1$, then the set \eqref{set of values}
becomes (keep the repeated values)
$$
 \{ 0, 0, -2/3,  2,    -2/3 \}.
$$
Then by Lemma \ref{Strang lemma} we have
$$
  -\frac{1}{3} I  \preceq    \frac{1}{2}(  P_k \hat{Z}_k  + \hat{Z}_k  P_k )   \preceq I . 
$$
Note that since by the induction hypothesis the eigenvalues of $\hat{Z}_k $  cannot achieve the extreme values of region $(-1/3, 1)$, the eigenvalues  of $ \frac{1}{2}(  P_k \hat{Z}_k  + \hat{Z}_k  P_k )   $ also cannot  \footnote{A more detailed argument is as follows. Since $-I/3  \preceq  \hat{Z}_k$, we can let $\beta_1 = -1/3 + \epsilon $ for a sufficiently small positive number $\epsilon$, while keeping $\alpha_1 = 0, \alpha_n = 1, \beta_n = 1 $. The set \eqref{set of values}
now	becomes	$	\{ 0, 0, -2/3 + 2 \epsilon,  2,    - \frac{ (4/3 - \epsilon)^2 }{ 4( 2/3 + \epsilon ) }  \}. $ Both $ -2/3 + 2 \epsilon$ and $- \frac{ (4/3 - \epsilon)^2 }{ 4( 2/3 + \epsilon ) } $ are strictly larger than $2/3$, thus the extreme value $-2/3$ cannot be achieved.  By a similar argument  the other extreme value $2$ also cannot be achieved.  }.
So we have
$$
-\frac{1}{3} I  \prec    \frac{1}{2}(  P_k \hat{Z}_k  + \hat{Z}_k  P_k )   \prec  I . 
$$
Thus according to \eqref{Z induction} we have
$$
 -\frac{1}{3} I  \prec Z \prec I. 
$$
This finishes the induction step. \textbf{Q.E.D.}

Remark: Where does the magical number $-1/3$ come from? It is actually the strange and complicated term 
$\frac{ 16 \alpha_1 \alpha_n \beta_1 \beta_n -
	(\beta_1 - \beta_n)^2 (\alpha_1 - \alpha_n^2 ) }{ 4(\alpha_1 + \alpha_n) (\beta_1 + \beta_n ) }$ in Strang's result \eqref{set of values}, which occurs due to the special structure of the Jordan product. 


\subsection{Proof of Proposition \ref{prop new: induction formula} (the induction formula) }\label{appen: proof of induction formula, n-block}
It is easy to build an induction formula from the expression \eqref{expected AQA' and projection}.
For example, when $n=3$, the matrix $  \sum_{ \sigma  } P_{\sigma_1} P_{\sigma_2} P_{\sigma_3 }  $
can be decomposed as the sum of 
$    P_1 (P_2 P_3 + P_3 P_2) + (P_2 P_3 + P_3 P_2) P_1  $ and two other similar terms (changing the outside part $P_1$ to $P_2, P_3$ and the inside part $P_2 P_3 + P_3 P_2$ correspondingly).
The inside part  $P_2 P_3 + P_3 P_2$ only involves two matrices, thus is a lower-dimensional analog. 
To make this even easier to see, denote $ X = P_1, Y = P_2, Z = P_3,  $ then
\begin{align*}
  2 \sum_{ \text{permutate} X,Y,Z } XYZ  
= [ X( YZ + ZY ) + ( YZ + ZY) X ] +  [ Y( XZ + ZX ) + ( XZ + ZX) Y ]  \\
+ [ Z(XY + YX) + (XY + YX)Z  ].
\end{align*}

A rigorous argument based on the above intuition is given as follows. 
Applying the formula \eqref{expected AQA' and projection} to the matrix $P_1, \dots, P_{k-1}, P_{k+1},
\dots, P_n$, and by the definition $\hat{A}_k = [A_1, \dots, A_{k-1}, A_{k+1}, \dots, A_n] $
and the definition of $\hat{ Q}_k $ in \eqref{hat(Q) def, n-block}, 
 we have 
$$
  I -  \hat{A}_k  \hat{ Q}_k  \hat{A}_k = \frac{1}{ (n-1)! } \sum_{\sigma = (\sigma_1 ,\dots, \sigma_{n-1}) \in \Gamma_k }   P_{\sigma_{n-1} } P_{\sigma_{n-1}} \dots  P_{\sigma_1}.
$$
We then have
\begin{align*}
 2(	I - AQA^T )  &  = \frac{2}{n!} \sum_{\sigma = (\sigma_1 ,\dots, \sigma_n) \in \Gamma}    P_{\sigma_n} P_{\sigma_{n-1}} \dots  P_{\sigma_1} \\
      &  = \frac{1}{n} \frac{1}{ (n-1)!}  \sum_{k=1}^n  \sum_{\sigma = (\sigma_1 ,\dots, \sigma_{n-1}) \in \Gamma_k }   ( P_k  P_{\sigma_{n-1} } P_{\sigma_{n-1}} \dots  P_{\sigma_1}
	+ P_{\sigma_{n-1} } P_{\sigma_{n-1}} \dots  P_{\sigma_1} P_k )    \\
	 & = \frac{1}{n} \sum_{k=1}^n  ( P_k (  I -  \hat{A}_k  \hat{ Q}_k  \hat{A}_k^T ) + (  I -  \hat{A}_k  \hat{ Q}_k  \hat{A}_k ) P_k   ), 
\end{align*}
which is the desired formula.

 \subsection{ Proof of Proposition \ref{prop: BR-ADMM convergence} }\label{appen: proof of BR ADMM}

We provide the proof of the expected convergence of BR-ADMM here, as this proof is a slightly smaller subset of the proof of Theorem \ref{Thm 1}. We will just describe the necessary modifications. 

We only need to prove a similar version of Theorem \ref{Thm 2}, i.e., the spectral radius of the expected update matrix of BR-ADMM is less than 1.
Throughout the proof, we need to change the matrix $Q = \frac{1}{|\Gamma|} \sum_{\sigma \in \Gamma} Q_{\sigma} $ to another one  defined as
\begin{equation}\label{Q BR def}
Q^{\text{BR} }   \triangleq \frac{1}{| \Gamma^{\text{BR}}  |} \sum_{\sigma \in \Gamma^{\text BR}} Q_{\sigma} ,
\end{equation}
where $ \Gamma^{\text{BR}}$ denotes the set of all possible permutations according to the Bernoulli randomization rule.  It is easy to see that $|\Gamma^{\text{BR}}| = 2^n$.  Other matrices such as $M$ should be changed accordingly. 

The proof of Theorem \ref{Thm 2} mainly consists of Lemma \ref{lemma 1} and Lemma \ref{lemma 2}. Since Lemma \ref{lemma 1} has nothing to do with the specific expression of
$Q$, so we only need to prove Lemma \ref{lemma 2} for BR-ADMM, i.e., the matrix $ A Q^{\text{BR}}  A^T $ has all eigenvalues in the region $( 0, 4/3 )$. Following the proof of Lemma \ref{lemma 2}, we divide the proof into three steps.

\textbf{Step 0}: Expression of $ Z^{\text{BR}} \triangleq	I - A Q^{\text{BR}}  A^T $. 
In Claim \ref{lemma of projection expression}, we have prove the expression \eqref{AQA' and projection} that $	I - AQ_{\sigma} A^T = P_{\sigma_n} P_{\sigma_{n-1}} \dots  P_{\sigma_1} $  for any permutation $\sigma$, which implies  $$ Z^{\text{BR}} = 	I - A Q^{\text{BR}}  A^T 
\overset{\eqref{Q BR def}}{=} \frac{1}{2^n} \sum_{\sigma \in \Gamma^{\text BR}}  P_{\sigma_n} P_{\sigma_{n-1}} \dots  P_{\sigma_1} . $$

\textbf{Step 1}: Induction formula. 
Notice that a characteristic  of the Bernoulli randomization rule is: the first block is either updated first or updated last. For instance, when $n = 4$, $(1,3,4,2)$ is a feasible permutation in $ \Gamma^{\text{BR}} $  and $(3,4,2, 1)$ is also a feasible permutation, but $(3,1, 4,2)$ is not feasible. After removing the first block, the rest $n-1$ blocks
form a permutation in $  \hat{\Gamma}_{\text{BR} }, $ where $  \hat{\Gamma}_{\text{BR} } $ is the set of all permutation
of $2,3, \dots, n$ according to the Bernoulli randomization rule. 
In other words, we have $ \Gamma^{\text{BR}}  = \{  (1, \hat{\sigma}), (\hat{\sigma}, 1), \text{ where } \hat{\sigma} \in \hat{\Gamma}^{\text{BR}}   \} $. 
Thus we have an induction formula
\begin{equation}\label{BR induction}
Z^{\text{BR}} = \frac{1}{2^n} \sum_{\sigma = (\sigma_1 ,\dots, \sigma_{n-1}) \in  \hat{\Gamma}_{\text{BR} } }   (  P_{1} P_{\sigma_{n-1}} \dots  P_{\sigma_1} + 
P_{\sigma_{n-1}} \dots  P_{\sigma_1} P_{1} ) = \frac{1}{2} (  P_{1} \hat{Z}^{\text{BR}}  + \hat{Z}^{\text{BR}} P_1 ), 
\end{equation}
where $ \hat{Z}^{\text{BR}}  $ is the lower dimensional analog of $ Z^{\text{BR}} $ for the rest $n-1$ blocks (after removing the first block).

\textbf{Step 2}: Applying mathematical induction. This step is almost the same as Step 2 of the proof of Lemma \ref{lemma 2}. More specifically, combining the induction hypothesis that $ \text{eig} ( \hat{Z}^{\text{BR}} )   \in (-1/3, 1), $  Strang's result Lemma  \ref{Strang lemma} and \eqref{BR induction}, we obtain the desired result $
  \text{eig} ( Z^{\text{BR}}  ) \in (-1/3, 1). $ 
This finishes the proof.

\section{Proof of Technical Results for Expected Convergence Rates}

\subsection{Proof of Claim \ref{claim: RP-ADMM spectral radius}}\label{subsec: proof of Claim of spec radius}
Suppose all the distinct eigenvalues of $I - QA^TA$ are $ 0 <  \tau_{N'} <  \dots < \tau_1 < 4/3 $, where $1  \leq N' \leq N$.
Denote $\tau_{\min} = \tau_{N'}, \tau_{\max} = \tau_1. $
According to Lemma \ref{lemma 1}, the expected update matrix of RP-ADMM $M $ has $2N'$ distinct  eigenvalues $\lambda_{k,1}, \lambda_{k,2}$ given by 
$$
\lambda_{k,1} = 1 - \tau_k + \sqrt{\tau_k (\tau_k - 1) }, \quad \lambda_{k,2} = 1 - \tau_k - \sqrt{\tau_k (\tau_k - 1) }, \; k=1,\dots, N'.
$$

Suppose the integer $m \in [1, N'+1]$ satisfies $ \tau_{m} \leq 1 < \tau_{m-1} $. When $m = 1$, every $\tau_k \leq 1$; when $m = N'+1$, every $\tau_k > 1$. 

For $N' \geq k \geq m $, i.e., $\tau_k \leq 1$,  we have  $\tau_k (\tau_k - 1) \leq 0$,  thus the two corresponding eigenvalues of $M$ are
$$
\lambda_{k,1} = 1 - \tau \pm i \sqrt{\tau( 1 -\tau) }, \lambda_{k,2} = 1 - \tau \pm i \sqrt{\tau( 1 -\tau) },
$$
which implies $|\lambda_{k,1}| = |\lambda_{k,2}| = \sqrt{ (1 - \tau_k)^2 + \tau_k (1-\tau_k) } = \sqrt{1 - \tau_k} $.
Thus $\rho_1 = \max_{ N' \geq k \geq m } \{ |\lambda_{k,1}| , |\lambda_{k,2}| \}  = \sqrt{1 - \tau_{N'}} = \sqrt{1 - \tau_{\min}} $ if such $k $ exists;
when such $k$ does not exist, i.e., $ \tau_k > 1 \; \forall \; k$ we denote $\rho_1 = 0 $ which equals $\sqrt{ (1 - \tau_{\min})_+ } $. In summary, we have $\rho_1 = \sqrt{ (1 - \tau_{\min})_+ }$.

For $m - 1 \geq k \geq 1 $, i.e., $\tau_k > 1$, 
we have $\tau_k( \tau_k - 1) > 0$. It is easy to verify $\lambda_{k,1} > 0 > \lambda_{k,2}$ and
$$
|\lambda_{k,2}| =  \tau_k - 1  +  \sqrt{\tau_k( \tau_k - 1 ) } >  1 - \tau_k + \sqrt{\tau_k(\tau_k - 1) }  =  | \lambda_{k,1} |.
$$
Denote $\rho_2 =  \max_{ m-1 \geq k \geq 1 } \{ |\lambda_{k,1}| , |\lambda_{k,2}| \} $,
then $\rho_2 = \max_{ m-1 \geq k \geq 1 } \{ |\lambda_{k,2}| \} = \max_{ m-1 \geq k \geq 1 } \{ \tau_k - 1  +  \sqrt{\tau_k( \tau_k - 1 ) }  \} = \tau_{\max} - 1  +  \sqrt{\tau_{\max}( \tau_{\max} - 1 )} $ if such $k$ exists; when such $k$ does not exist, i.e., $\tau_k \leq 1 \; \forall \; k$, we denote $\rho_2 = 0$ which equals   $ ( \tau_{\max} - 1 )_+  +  
\sqrt{\tau_{\max}( (\tau_{\max} - 1)_+ )} $.

Combining the two scenarios, we have 
$ \rho(M) = \max_{ N' \geq k \geq 1 } \{|\lambda_{k,1}| , |\lambda_{k,2}|   \}
= \max\{\rho_1, \rho_2 \} = \max \{ \sqrt{ (1 - \tau_{\min})_+ },   \; ( \tau_{\max} - 1 )_+  +  
\sqrt{\tau_{\max}( (\tau_{\max} - 1)_+ )}  \} . $

Next, we prove
\begin{equation}\label{further bound rho M}
\begin{split}
( \tau_{\max} - 1 )_+ + \sqrt{ \tau_{\max} (\tau_{\max} -1 )_+ } 
& \leq \max \left\{ 1 - \frac{3}{4} (4 - 3 \tau_{\max}) , 0 \right\},    \\
\sqrt{ (1 - \tau_{\min})_+ }  &  \leq 1 - \frac{1}{2}\tau_{\min}. 
\end{split}
\end{equation}
In fact, when $4/3 \geq \tau \geq 1$, we have
$1 - ( \tau - 1  + \sqrt{ \tau (\tau -1 ) }) = 2 - \tau - \sqrt{ \tau (\tau -1 ) }
= \frac{ (2 - \tau )^2 -  \tau (\tau -1 ) }{ 2 - \tau + \sqrt{ \tau (\tau -1 ) }   } = \frac{ 3 - 4 \tau }{  2 - \tau + \sqrt{ \tau (\tau -1 ) }}	\geq \frac{3}{4} (3 - 4 \tau ) $, thus $ \tau - 1  + \sqrt{ \tau (\tau -1 ) } \leq 1 -  \frac{3}{4} (3 - 4 \tau )  .$
When $\tau < 1$, clearly $ \tau - 1  + \sqrt{ \tau (\tau -1 ) } = 0$. Thus
$ (\tau - 1 )_+  + \sqrt{ \tau (\tau -1 )_+ } \leq \max \{ 0,  1 - \frac{3}{4} (4 - 3 \tau)  \} .$
For the second relation, if $0 \leq \tau < 1$ then $ \sqrt{1 - \tau} = 1 - \frac{\tau}{1 + \sqrt{1-\tau}} \leq 1 - \frac{\tau }{2} $;
if $1 \leq  \tau \leq 4/3$ then $  \sqrt{(1 - \tau)_+}  = 0 < 1 - \frac{1}{2}\tau.  $	Thus $  \sqrt{(1 - \tau)_+}  \leq 1 - \frac{1}{2}\tau$ holds for any $\tau \in [0,4/3]$. 

Substituting \eqref{further bound rho M} into the expression of $\rho(M)$, we obtain the desired inequality
$$
\rho(M) \leq \max\left\{ 1 - \frac{3}{4} (4 - 3 \tau_{\max}) , \; \; 1 - \frac{1}{2}\tau_{\min} \right\}. 
$$ 

\subsection{Proof of Lemma \ref{lemma of stronger bound}}\label{subsec: lemma 5 proof}

This is one of the two  main lemmas of proving the expected convergence rate of RP-ADMM
(the other is the expected convergence rate of RP-CD).  

The proof outline of Lemma \ref{lemma of stronger bound} and the main techniques are described below.
The previous proof for the expected convergence of RP-ADMM in Section \ref{sec: proof of Lemma 2 for n-block} is not strong enough to prove a convergence rate. We have to obtain a more refined estimate of the spectral radius of $AQA^T$. 
To do so, we transform the induction formula in Proposition \ref{prop new: induction formula}
to a ``dual'' form: instead of $AQA^T$, we consider a similar matrix $QA^T A $. 
We then apply the two simple techniques used in the proof of Lemma \ref{lemma 1}: factorize and rearrange, and reduce the dimension by eliminating a variable from the eigenvalue equation.
We obtain a somewhat complicated inequality relating $  \lambda_{\max}(Q A^T A) $
and its lower-dimensional analog $ \lambda_{\max}( \hat{Q} \hat{A}^T \hat{A}) $.  Finally, we perform a detailed analysis of the inequality to prove the desired bound. 







\subsubsection{Step 1: Mathematical Induction and Induction Formula}

Define a sequence 
$\{ \alpha_k \} _{k=1}^{\infty} $ such that  
\begin{equation}\label{alpha recursion, 1st time}
\alpha_1 = 1/3, \quad \alpha_{ k + 1} =  h(\alpha_k) \triangleq  \frac{ \alpha_k}{ 8 } \frac{ 16 - 3 \alpha_k }{ 2 + 3 \alpha_k }.
\end{equation}
It is easy to verify that $0< \alpha_{k+1} <  \alpha_k \leq 1/3$ for all $k$. The following claim provides a bound of $\alpha_k$ (the proof will be given in Section  \ref{subsubsec: proof of Claim of alpha bound}). 
\begin{claim}\label{claim of distance sequence 4/3 - lambda} $\;$
	Suppose the sequence $\{ \alpha_k \} _{k=1}^{\infty} $ satisfies \eqref{alpha recursion, 1st time},
	then $\alpha_k \geq \frac{4}{9(k+1)}, \forall \; k \geq 1. $
\end{claim}

According to this claim, to prove the desired result $\lambda_{\max} (AQA^T) \leq \frac{4}{3} - \frac{4}{9(k+1) }  $,
we only need to prove the following result:
\begin{equation}\label{alpha lambda relation}
   \text{eig}(AQA^T) \subseteq (0, \frac{4}{3} -  \alpha_n ].  
\end{equation}

We prove this result by mathematical induction. 
When $ n =1 $,  
since $ A^TA= A_1^T A_1 = I $, we have $\lambda_{\min}( AQA^T  ) =  \lambda_{\max}( AQA^T  )  =  1 = \frac{4}{3} - \alpha_1 $.

Suppose the result holds for $n-1$, i.e., for a problem with $n-1$ blocks, the eigenvalues of the corresponding matrix $\hat{A} \hat{Q} \hat{A}^T$ 
lie in the region $( 0 ,  \frac{4}{3} - \alpha_{n-1} )$.

Next, we build the induction formula, which is the dual form of the one we derived before. 
According to  \eqref{Z induction}, we have 
$$
  2(	I - AQA^T )  = \frac{1}{n} \sum_{k=1}^n  ( P_k (  I -  \hat{A}_k  \hat{ Q}_k  \hat{A}_k ) + (  I -  \hat{A}_k  \hat{ Q}_k  \hat{A}_k ) P_k   ), 
$$
which can be rewritten as  
\begin{equation}\label{Z induction again}
AQ A^T  = \frac{1}{n} \sum_{k=1}^n  \left[ I -  \frac{1}{2}P_k (  I -  \hat{A}_k  \hat{ Q}_k  \hat{A}_k ) -  \frac{1}{2} (  I -  \hat{A}_k  \hat{ Q}_k  \hat{A}_k )P_k \right]
\end{equation}

Note that 
\begin{align*}
I - P_k (  I -  \hat{A}_k  \hat{ Q}_k  \hat{A}_k )
&  = I -  (I - A_k A_k^T ) (  I -  \hat{A}_k  \hat{ Q}_k  \hat{A}_k^T ) \\
& =  A_k A_k^T + \hat{A}_k  \hat{ Q}_k  \hat{A}_k^T 
- A_k A_k^T  \hat{A}_k  \hat{ Q}_k  \hat{A}_k^T    \\
& = [\hat{A}_k, A_k ] \begin{bmatrix}
\hat{ Q}_k   &  0  \\
-   A_k^T  \hat{A}_k  \hat{ Q}_k    &     I   \\
\end{bmatrix}  [\hat{A}_k, A_k ]^T  
\end{align*}

Thus the symmetrized version
\begin{align}
&  I -  \frac{1}{2}P_k (  I -  \hat{A}_k  \hat{ Q}_k  \hat{A}_k ) -  \frac{1}{2} (  I -  \hat{A}_k  \hat{ Q}_k  \hat{A}_k )P_k   \\
& = [\hat{A}_k, A_k ] \begin{bmatrix}
\hat{ Q}_k   &  -  \frac{1}{2}   \hat{ Q}_k^T \hat{A}_k^T  A_k   \\
-  \frac{1}{2}  A_k^T  \hat{A}_k  \hat{ Q}_k    &     I   \\
\end{bmatrix}  [\hat{A}_k, A_k ]^T    \\
&  =  \bar{A}_k Q_k  \bar{A}_k^T,   \label{temp each term AQA}
\end{align}
where  in the last step we use the definitions
\begin{equation}\label{Qk def}
\bar{A}_k \triangleq  [\hat{A}_k, A_k], \quad Q_k \triangleq \begin{bmatrix}
\hat{Q}_k &     -\frac{1}{2} \hat{Q}_k  W_k      \\
-\frac{1}{2}W_k^T \hat{Q}_k &    I_{d_k \times d_k}    \\
\end{bmatrix}
\end{equation}
Sum up \eqref{temp each term AQA} for $k = 1, \dots, n$ and applying \eqref{Z induction again}, we have
\begin{equation}\label{induction formula for genenral n, n-block}
AQ A^T  = \frac{1}{n} \sum_{k=1}^n   \bar{A}_k Q_k  \bar{A}_k^T.  
\end{equation}


Consequently,
\begin{equation}\label{decompose eig bound, block}
\frac{1}{n} \sum_{k=1}^n \lambda_{\min}(  \bar{A}_k Q_k  \bar{A}_k^T ) \leq  \lambda_{\min}( AQA^T ) \leq  \lambda_{\max}(  AQA^T ) \leq \frac{1}{n}\sum_{k=1}^n \lambda_{\rmax}( \bar{A}_k Q_k  \bar{A}_k^T).
\end{equation}

To prove $\text{eig}(AQA^T) \subseteq (0, \frac{4}{3} -  \alpha_n ] $, we only need to prove 
for any $k=1,\dots, n$,
\begin{equation}\label{bound eig for each k, n-block}
\text{eig}( \bar{A}_k Q_k  \bar{A}_k^T  )  \subseteq  (0,  \frac{4}{3} -  \alpha_n ) .
\end{equation}

Note that $\hat{Q}_k$ only depends on the entries of $\hat{A}_k^T \hat{A}_k \in \dR^{(N - d_k) \times (N - d_k)}$ which has $(n-1)\times (n-1)$ blocks,
thus by the induction hypothesis, we have
\begin{equation}\label{hypothesis of Qk, n-block}
\text{eig}( \hat{Q}_k \hat{A}_k^T \hat{A}_k ) \subseteq (0,  \frac{4}{3} - \alpha_{n-1} ] .
\end{equation}

\begin{prop}\label{prop 2: eigenvalue induction, n-block}  $\; $
	Suppose $A = [\hat{A}_n, A_{n}] \in \dR^{N \times N} $ is a non-singular matrix, where
	$\hat{A}_n \in \dR^{N \times (N-d_n)} $, and $A_n \in \dR^{ N \times d_n}$   
	satisfies $ A_n^T A_n = I_{d_n \times d_n}$.
	Suppose $ \hat{Q}_n \in \dR^{(N-d_n) \times (N-d_n)}  $ is symmetric, satisfying
	\begin{equation}\label{hypothesis, n-block}
\text{eig}(A  \hat{Q}_n A^T) \subseteq  (0,	\frac{4}{3} -  \alpha_{n-1} ],
	\end{equation} 
where $\{ \alpha_k \}$ is defined in \eqref{alpha recursion, 1st time}.	Define
	\begin{equation}\label{Qn def, in prop 2, n-block}
	W_n \triangleq   \hat{A}_n^T A_n \in \dR^{ (N-d_n) \times d_n }, \quad  Q_n \triangleq \begin{bmatrix}
	\hat{Q}_n &     -\frac{1}{2} \hat{Q}_n  W_n      \\
	-\frac{1}{2}W_n^T \hat{Q}_n &     I_{d_n \times d_n}    \\
	\end{bmatrix} .
	\end{equation}
	Then $ \text{eig}(A Q_n A^T) \subseteq ( 0,  \frac{4}{3} - \alpha_n ] $.
\end{prop}
The proof of Proposition \ref{prop 2: eigenvalue induction, n-block} will be divided into two parts, and
  given in Section \ref{sec: proof of prop 2, n-block} and Section \ref{sec: proof of prop 2, more precise bound}. 

We claim that \eqref{bound eig for each k, n-block}
follows from the induction hypothesis \eqref{hypothesis of Qk, n-block} and the expressions
of $\bar{A}_k $ and $Q_k$ in  \eqref{Qk def}. 
In fact, the above proposition directly proves \eqref{bound eig for each k, n-block} for $k = n$. If we replace $A, \hat{A}_n, A_n, \hat{Q}_n, Q_n $ by $\bar{A}_k, \hat{A_k},  A_k, \hat{Q}_k, Q_k$ respectively in the following proposition, we will obtain \eqref{bound eig for each k, n-block} for any $k$.
Finally, as mentioned earlier, the desired result $\text{eig}(AQA^T) \subseteq (0,  \frac{3}{4} - \alpha_n ] $ in Lemma \ref{lemma 2} follows immediately from \eqref{bound eig for each k, n-block} and \eqref{decompose eig bound, block}.


\subsubsection{Step 2:  Relation Between $\lambda_{\max}(  A_n Q A_n^T ) $ and its analog}\label{sec: proof of prop 2, n-block}

In this subsection, we provide a proof of a weaker result $ \text{eig}(A Q_n A^T) \subseteq ( 0,  \frac{4}{3}  ) $
under the conditions of Prop. \ref{prop 2: eigenvalue induction, n-block}; 
the proof of the desired result $ \text{eig}(A Q_n A^T) \subseteq ( 0,  \frac{4}{3} - \alpha_n ] $ will
be provided in the next subsection. 

For simplicity, throughout this proof, we denote
$$
W \triangleq W_n \in \dR^{ (N - d_n) \times d_n}, \ \hat{Q} \triangleq \hat{Q}_n \in \dR^{(N-d_n)\times (N-d_n)}, \ \hat{A} \triangleq \hat{A}_n \in \dR^{ N \times (N-d_n)}.
$$
According to the assumption of Prop. \ref{prop 2: eigenvalue induction, n-block}, we have
 \begin{equation}\label{lambda induction hypothesis}
    \hat{\lambda} \triangleq  \lambda_{\max} ( A  \hat{Q} A^T ) \in (0, \frac{4}{3} -  \alpha_{n-1}]. 
 \end{equation}

We first prove 
\begin{equation}\label{theta def, n-block}
0 \preceq \Theta  \triangleq W^T \hat{Q} W \prec \frac{4}{3}I .
\end{equation}
Since $ \text{eig}(\hat{Q} \hat{A}^T \hat{A} ) \subseteq (0, \infty) $  and $\hat{A}$ is non-singular,
thus $\hat{Q} \succ 0$. Then we have $\Theta =  W^T \hat{Q} W \succeq 0 $, which proves the first relation of \eqref{theta def, n-block}.
By the definition $W =  \hat{A}^T A_n $ we have
\begin{equation}\label{Theta spec bound}
\begin{split}
\rho( \Theta ) = \rho( A_n^T \hat{A} \hat{Q} \hat{A}^T A_n )
= \max_{v \in \dR^{ d_n \times 1 }, \| v\| = 1} v^T  A_n^T \hat{A} \hat{Q} \hat{A}^T A_n v  \\
\leq \rho( \hat{A} \hat{Q} \hat{A}^T  ) \max_{v \in \dR^{ d_n \times 1 }, \| v\| = 1} \| A_n v\|^2
= \rho( \hat{A} \hat{Q} \hat{A}^T  )  \| A_n\|^2
= \rho( \hat{A} \hat{Q} \hat{A}^T ) < \frac{4}{3} ,
\end{split}
\end{equation}
where the last equality is due to the assumption $ A_n^T A_n = I$, and the last inequality is due to the assumption \eqref{hypothesis, n-block}.
By \eqref{Theta spec bound} we have $\Theta \prec \frac{4}{3}I $, thus \eqref{theta def, n-block} is proved.

We apply a trick that we have previously used: factorize $Q_n$ and change the order of multiplication.
To be specific,  $Q_n$ defined in \eqref{Qn def, in prop 2, n-block} can be factorized as
\begin{equation}\label{Q factorize, n-block}
Q_n = \begin{bmatrix}
I &     0     \\
-\frac{1}{2} W^T &    I    \\
\end{bmatrix}
\begin{bmatrix}
\hat{Q} &     0      \\
0 &    I - \frac{1}{4}W^T \hat{Q} W   \\
\end{bmatrix}
\begin{bmatrix}
I &   -\frac{1}{2}  W     \\
0 &    I    \\
\end{bmatrix} = J  \begin{bmatrix}
\hat{Q} &     0      \\
0 &    C   \\
\end{bmatrix}
J^T ,
\end{equation}
where $J \triangleq \begin{bmatrix}
I &     0     \\
-\frac{1}{2} W^T &    I    \\
\end{bmatrix}$, $I$ in the upper left block denotes the $(N-d_n)$-dimensional identity matrix,
$I$ in the lower right block denotes the $d_n$-dim identity matrix,
and
\begin{equation}\label{def of c}
C \triangleq  I - \frac{1}{4} W^T \hat{Q} W \in \dR^{d_n \times d_n}.
\end{equation}

It is easy to prove
\begin{equation}\label{PD, n-block}
\text{eig}( A Q_n A^T ) \subseteq (0, \infty).
\end{equation}
In fact, we only need to prove $Q_n \succ 0$.
According to \eqref{Q factorize, n-block}, we only need to prove $  \begin{bmatrix}
\hat{Q} &     0      \\
0 &    C   \\
\end{bmatrix} \succ 0 .$
This follows from $\hat{Q} \succ 0 $ 
and the fact $ C =  I - \frac{1}{4}W^T \hat{Q} W \overset{\eqref{theta def, n-block}}{\succ} I - \frac{1}{3}I \succ 0. $
Thus \eqref{PD, n-block} is proved.

It remains to prove
\begin{equation}\label{4/3 side bound}
\rho(A Q_n A^T) < \frac{4}{3} .
\end{equation}
Denote $\hat{B} \triangleq \hat{A}^T \hat{A} \in \dR^{(N-d_n)\times (N-d_n)}$, then we can write $A^T A$ as
\begin{equation}
A^T A = \begin{bmatrix}
\hat{B} &  W      \\
W^T &    I    \\
\end{bmatrix}.
\end{equation}
We simplify the expression of $\rho(AQ_n A^T) $ as follows:  
\begin{equation}\label{transform simple, n-block}
\rho(AQ_n A^T ) = \rho \left(
A J
\begin{bmatrix}
\hat{Q} &     0      \\
0 &    C   \\
\end{bmatrix}
J^T A^T
\right)
=   \rho \left(
\begin{bmatrix}
\hat{Q} &     0      \\
0 &    C  \\
\end{bmatrix}
J^T A^T A J
\right) .
\end{equation}
By algebraic computation, we have
\begin{equation}
\begin{split}
J^T A^T A J   =
&
\begin{bmatrix}
I &   -\frac{1}{2}  W     \\
0 &    I    \\
\end{bmatrix}
\begin{bmatrix}
\hat{B} &  W      \\
W^T &    I    \\
\end{bmatrix}
\begin{bmatrix}
I &     0     \\
-\frac{1}{2} W^T &    I    \\
\end{bmatrix}
\\   
=
&
\begin{bmatrix}
I &   -\frac{1}{2}  W     \\
0 &    I       \\
\end{bmatrix}
\begin{bmatrix}
\hat{B} - \frac{1}{2} W W^T &     W     \\
\frac{1}{2} W^T &    I    \\
\end{bmatrix}
= \begin{bmatrix}
\hat{B} - \frac{3}{4} W W^T &  \frac{1}{2}  W     \\
\frac{1}{2} W^T &    I   \\
\end{bmatrix},
\end{split}
\end{equation}
thus
\begin{equation}\label{Z expression, n-block}
Y \triangleq \begin{bmatrix}
\hat{Q} &     0      \\
0 &    C  \\
\end{bmatrix}
J^T A^T A J
= \begin{bmatrix}
\hat{Q} &     0      \\
0 &    C  \\
\end{bmatrix}
\begin{bmatrix}
\hat{B} - \frac{3}{4} W W^T & \frac{1}{2}   W     \\
\frac{1}{2} W^T &    I   \\
\end{bmatrix}
=
\begin{bmatrix}
\hat{Q} \hat{B} - \frac{3}{4} \hat{Q} W W^T & \frac{1}{2}  \hat{Q} W     \\
\frac{1}{2}  C W^T &    C   \\
\end{bmatrix}.
\end{equation}

Suppose $\lambda > 0 $ is the maximal eigenvalue of $Y$. 
According to \eqref{transform simple, n-block} that $\rho(A Q_n A^T) = \rho(Y)$, we also have
$ \lambda = \lambda_{\max}(A Q_n A^T) $. 
To prove \eqref{4/3 side bound}, we only need to prove
\begin{equation}\label{lambda bound, n-block}
\lambda  < \frac{4}{3}.
\end{equation}

Suppose $v \in \dR^{N \times 1}\backslash\{0\}$ is the eigenvector corresponding to $\lambda$, i.e. $Zv = \lambda v$.
Partition $v$ into $v = \begin{bmatrix} v_1 \\ v_0 \end{bmatrix} $,
where $v_1 \in \dR^{N- d_n}, v_0 \in \dR^{d_n}$.
According to the expression of $Z$ in \eqref{Z expression, n-block}, $Z v = \lambda v$ implies
\begin{subequations}\label{Z eigen equations}
	\begin{align}
	(  \hat{Q} \hat{B} - \frac{3}{4} \hat{Q} W W^T ) v_1 +  \frac{1}{2}  \hat{Q} W  v_0 = \lambda v_1, \label{Z eigen equation a, n-block} \\
	\frac{1}{2}  C W^T v_1 + C v_0  = \lambda v_0.   \label{Z eigen equation b, n-block}
	\end{align}
\end{subequations}
If $\lambda I - C $ is singular, i.e. $\lambda$ is an eigenvalue of $C$, then by \eqref{theta def, n-block}
we have $\frac{2}{3} I \prec  C = 1 - \frac{1}{4}\Theta \preceq I $, which implies $ \lambda \leq 1$,
thus \eqref{lambda bound, n-block} holds.
In the following, we assume
\begin{equation}\label{lamI - C full rank}
\lambda I - C  \text{ is non-singular}.
\end{equation}
An immediate consequence is
$$
v_1 \neq 0,
$$
since otherwise  \eqref{Z eigen equation b, n-block} implies $ C v_0  = \lambda v_0 $, which combined with \eqref{lamI - C full rank} leads to $v_0 = 0$ and thus $v = 0$, a contradiction.

By \eqref{Z eigen equation b, n-block} we get
$$
v_0 = \frac{1}{2} (\lambda I - C)^{-1} C W^T v_1.
$$
Plugging into \eqref{Z eigen equation a, n-block}, we obtain
\begin{equation}\label{crucial eigen expression, n-block}
\lambda v_1 = (  \hat{Q} \hat{B} - \frac{3}{4} \hat{Q} W W^T ) v_1 + \frac{1}{2}  \hat{Q} W  \frac{1}{2} (\lambda I - C)^{-1} C W^T v_1
= ( \hat{Q} \hat{B} +  \hat{Q} W \Phi W^T   ) v_1,
\end{equation}
where
\begin{equation}\label{phi def, n-block}
\begin{split}
\Phi \triangleq - \frac{3}{4} I + \frac{1}{4} (\lambda I - C)^{-1} C = -I + \frac{1}{4} [ I + (\lambda I - C)^{-1} C ]   \\
= -I + \frac{\lambda}{4} (\lambda I - C)^{-1}
= -I + \lambda [ (4\lambda - 4)I + \Theta ]^{-1}.
\end{split}
\end{equation}
Here  we have used the definition $C = I - \frac{1}{4} W^T \hat{Q} W = I - \frac{1}{4} \Theta$. Since $\Theta$ is a symmetric matrix, $\Phi$ is also a symmetric matrix.

Define 
\begin{equation}\label{H def}
\tilde{H} \triangleq \hat{Q} W \Phi W^T   \in \dR^{(N-d_n) \times (N - d_n)}   ,
, H \triangleq  W^T  \hat{Q} W \Phi  = \Theta \Phi  \in \dR^{  d_n \times d_n  } .
\end{equation}
As a well-known linear algebra result, $\tilde{H}$ and $H$ have the same non-zero eigenvalues. Note
that $\lambda_{\max}(H)$ may not be equal to $\lambda_{\max}(\tilde{H})$ due to the possible zero eigenvalues.
Nevertheless,  we can define $ \lambda_{\max}^+(X) \triangleq \max \{ \lambda_{\max}(X ), 0  \} $, and then we have
$$  \lambda_{\max}^+(\tilde{H}) =  \lambda_{\max}^+(H ) . $$

According to \eqref{H def} and \eqref{phi def, n-block}, we know
\begin{align*}
 H   & = \Theta \Phi = \Theta ( -I + \lambda [ (4\lambda - 4)I + \Theta ]^{-1} ) \\
 & = -\Theta + \lambda \Theta [ (4\lambda - 4)I + \Theta ]^{-1}      \\
& = -\Theta + \lambda ( I - (4\lambda - 4)( [ (4\lambda - 4)I + \Theta ]^{-1}  ) )   \\
 & =   -\Theta + \lambda I - \lambda(4\lambda - 4)  [ (4\lambda - 4)I + \Theta ]^{-1} . 
\end{align*}
It is well-known that if $\alpha I+\Theta$ is invertible, then $\Theta$ has an eigenvalue $\theta$ iff   $ (\alpha I + \Theta )^{-1}$ has an eigevalue   $ (\alpha + \theta)^{-1} $, and the corresponding eigen-vectors are the same. Similarly, since we already assumed $ (4\lambda - 4)I + \Theta $ is invertible,    $\theta$ is an eigenvalue of $\Theta$ iff 
$ H = -\Theta + \lambda I - \lambda(4\lambda - 4)  [ (4\lambda - 4)I + \Theta ]^{-1}  $ has an eigenvalue
$  -\theta + \lambda - \lambda(4\lambda - 4)[ (4\lambda - 4) + \theta ]^{-1} $.
Recall that $ \Theta = W^T  \hat{Q} W  $ satisfies $ 0 \preceq \Theta \preceq \hat{\lambda} I $,
thus any eigenvalue $\theta$ satisfies $ 0 \leq \theta \leq \hat{\lambda} $. 
Therefore 
\begin{equation}\label{H eig bound}
\lambda_{\max}(H) \leq \max_{ \theta \in [0, \hat{\lambda} ] }  \{-\theta + \lambda -
\frac{ \lambda(4\lambda - 4)  }{ (4\lambda - 4) + \theta } \} \triangleq g(\theta). 
\end{equation}

Since $v_1 \neq 0$, without loss of generality, we can assume $\| v_1 \| = 1$. We have
\begin{equation}\label{lambda final bound, precise}
\begin{split}
\lambda = v_1 ^T \hat{Q} \hat{B}  v_1 + v_1^T \tilde{H} v_1 \leq \hat{\lambda} + v_1^T \tilde{H} v_1
\leq \hat{\lambda} + \lambda_{\max}^+(\tilde{H}) = \hat{\lambda} + \lambda_{\max}^+(H)  \\
\leq \hat{\lambda} + \max \{ 0,  \max_{ \theta \in [0, \hat{\lambda} ] }  \{-\theta + \lambda -
\frac{ \lambda(4\lambda - 4)  }{ (4\lambda - 4) + \theta } \}    \},
\end{split}
\end{equation}
where the first equality is due to \eqref{crucial eigen expression, n-block}, the first inequality is due to the induction hypothesis, the second inequality uses the obvious relation $ \lambda_{\max}( \tilde{H} ) \leq  \lambda_{\max}^+( \tilde{H} )   $, and the last inequality is due to \eqref{H eig bound}.

To prove \eqref{lambda bound, n-block}, we consider two cases.

Case 1: $ \max_{ \theta \in [0, \hat{\lambda}] } g(\theta) \leq 0  .$ In this case, $\lambda \leq \hat{\lambda} < 4/3$, 
where the first inequality is due to \eqref{lambda final bound, precise}, and the second inequality is due to the induction hypothesis. Thus in Case 1 \eqref{lambda bound, n-block} holds.

Case 2: $ \max_{ \theta \in [0, \hat{\lambda}] } g(\theta) > 0  .$ 
Then there exists some $\theta \geq 0$ such that $g(\theta) > 0$. Note that 
$g(\theta)  $ can also be expressed as $g(\theta) = \theta( - 1 + \frac{\lambda}{  (4\lambda - 4) + \theta } )$, 
thus 
\begin{equation}\label{-1 + ... >0}
- 1 + \frac{\lambda}{  (4\lambda - 4) + \theta } > 0.
\end{equation}
If $\lambda < 1$, then \eqref{lambda bound, n-block} already holds; so we can assume $\lambda > 1$.
Thus \eqref{-1 + ... >0} implies
$ 1 < \frac{\lambda}{  (4\lambda - 4) + \theta } \leq \frac{\lambda}{  4\lambda - 4 } $,
which leads to $\lambda < \frac{4}{3}$. Thus in Case 2 \eqref{lambda bound, n-block} also holds.
This finishes the proof of \eqref{lambda bound, n-block}.

Remark: The proof of this subsection can lead to an alternative proof of Lemma \ref{lemma 2}.
In particular, the induction step (Step 2) of Section \ref{sec: outline proof of Lemma 2} can be replaced 
by the proof here.
The proof presented here is more complicated and less intuitive than the one in Section 
\ref{sec: outline proof of Lemma 2} (which is just a straightforward application of Strang's result Lemma \ref{Strang lemma}, but the benefit is that it can help establish a stronger bound of $\lambda$, as done in the next subsection.

\subsubsection{Step 3:  More Precise Bound of $\lambda$}\label{sec: proof of prop 2, more precise bound}

We will continue the proof in Section \ref{sec: proof of prop 2, n-block}, to further prove 
\begin{equation}\label{lambda bound for n}
\lambda = \lambda_{\max}(A Q_n A^T) \leq  4/3 - \alpha_n. 
\end{equation}
We rewrite \eqref{lambda final bound, precise} as follows:
\begin{equation}\label{lambda final bound, precise rewrite}
  	\lambda  \leq \hat{\lambda} + \max \{ 0,  \max_{ \theta \in [0, \hat{\lambda} ] }  g(\theta)   \}, \text{ where } g( \theta ) = \lambda -	\frac{ \lambda(4\lambda - 4)  }{ 4\lambda - 4 + \theta } -\theta    .
\end{equation}
If $\lambda < 1$, then we are done since $1 \leq 4/3 - \alpha_n$.  Assume $1 \leq \lambda < 4/3 $ from now on. 

We first analyze the function $g(\theta)$. Taking the derivative of $g$, we get
$$
 g'(\theta) = \frac{ \lambda(4 \lambda -4) }{ (4\lambda - 4 + \theta )^2 } - 1
   = \frac{ ( \sqrt{\lambda(4 \lambda -4)} +  4\lambda - 4 + \theta )( \sqrt{\lambda(4 \lambda -4)} -  4\lambda + 4 - \theta  )  }{ (4\lambda - 4 + \theta )^2 } .
$$
Since $\lambda > 1$ and $\theta \geq 0$, the term in the first bracket in the numerator is positive. 
Define
$$
 \theta^* = \sqrt{ \lambda(4 \lambda -4) } -  4\lambda + 4  > 0,
$$
where the inequality holds due to $\lambda < 4/3.$
Then we have
$$
  g'(\theta) \begin{cases}
  \geq 0, &   \theta \leq \theta^*; \\
  \leq 0, &  \theta > \theta^*.  
  \end{cases}
$$
Therefore, $g(\theta)$ is increasing in $[0, \theta^*]$ and decreasing in $[\theta^*, \infty )$. This implies
\begin{equation}\label{g theta bound}
  g(\theta) \leq g(\theta^*), \; \forall \theta \; \geq 0.
\end{equation}
According to $0 < \lambda < 4/3$, we have $\lambda > \sqrt{ \lambda(4 \lambda -4) } = 4\lambda - 4 + \theta^* \Rightarrow 
-1 + \frac{ \lambda }{ 4\lambda - 4 + \theta^*} > 0 \Rightarrow
 g(\theta^*) > 0 .$ 
 Together with \eqref{g theta bound} we obtain $ \max \{ 0,  \max_{ \theta \in [0, \hat{\lambda} ] }  g(\theta) \} \leq g(\theta^*) $. Substituting into   \eqref{lambda final bound, precise rewrite}, we obtain
 $$
  \lambda \leq \hat{ \lambda} + g( \theta^* ).
 $$
We will derive an inequality on $\lambda$ and $\hat{\lambda}$ from the above relation as below. 
Substituting the expression of $g( \cdot )$ into the relation, we obtain
\begin{equation*}%
    \lambda \leq \hat{ \lambda} + \lambda -	\frac{ \lambda(4\lambda - 4)  }{ 4\lambda - 4 + \theta^* } -\theta^*
  \Longrightarrow  \hat{ \lambda}  \geq  \frac{ \lambda(4\lambda - 4)  }{ 4\lambda - 4 + \theta^* } + \theta^*
  = \sqrt{ \lambda(4\lambda - 4) } + \theta^* 
  =      2 \sqrt{ \lambda(4 \lambda -4) } -  4\lambda + 4.
\end{equation*}
This implies
\begin{subequations}
\begin{align}
& \hat{\lambda}^2 + (4\lambda -4)^2 + 2 \hat{ \lambda} (4\lambda -4 ) \geq 4\lambda(4\lambda -4)  \nonumber \\
\Longleftrightarrow \quad &    \hat{\lambda}^2 - \lambda^2 + 2 (\hat{ \lambda} - \lambda) (4\lambda -4 ) + (\lambda - (4\lambda -4))^2 \geq 0   \nonumber \\
\Longleftrightarrow \quad &    (\hat{\lambda} - \lambda)( \hat{ \lambda} + \lambda ) + 2 (\hat{ \lambda} - \lambda) (4\lambda -4 ) + (4 - 3 \lambda  )^2 \geq 0. \label{relation of lambda and lambda hat purely}
\end{align}
\end{subequations}

Define
\begin{equation}\label{delta def}
 \delta = 4/3 - \lambda \in (0, 1/3),  \quad \hat{ \delta } = 4/3 - \hat{ \lambda} \in (0, 4/3). 
\end{equation}
Substituting into \eqref{relation of lambda and lambda hat purely}, we obtain
\begin{align*}
 & ( \delta - \hat{ \delta} )(8/3 - \delta - \hat{ \delta} ) + (\delta - \hat{ \delta} )( 8/3 - 8 \delta ) +  9 \delta^2 \geq 0 \\
 \Longleftrightarrow \quad &  ( \delta - \hat{ \delta} )(16/3 - 9\delta - \hat{ \delta} )+  9 \delta^2 \geq 0 \\
 \Longleftrightarrow \quad &   \frac{16}{3} \delta - \frac{16}{3} \hat{\delta} + 8 \hat{ \delta} \delta +  \hat{ \delta}^2  \geq 0 \\
 \Longleftrightarrow \quad & \delta \geq \frac{ \hat{ \delta}(16 - 3 \hat{ \delta}) }{ 8 (2 + 3 \hat{ \delta}) } = h(\hat{\delta}) .
\end{align*} 

It is easy to verify that $h(t)$ is increasing in $t \in [0,4/3]$; in fact,
$h'(t) = \frac{36}{ (2+ 3t)^2 } - 1 = \frac{ (8 + 3t) (4-3t) }{ (2+3t)^2}  \geq 0 $ for $t \in [0, 4/3]$.  
According to \eqref{lambda induction hypothesis}, we have $\hat{ \delta } = 4/3 - \hat{ \lambda} 
\geq \alpha_{n-1}$. Applying the monotonicity of $h$, we have
$$
  \delta \geq h(\hat{ \delta }) \geq h( \alpha_{n-1} ) = \alpha_n,
$$
which combined with \eqref{delta def} leads to \eqref{lambda bound for n}. This finishes the proof of Proposition 
\ref{prop 2: eigenvalue induction, n-block}. 

\subsubsection{Proof of Claim \ref{claim of distance sequence 4/3 - lambda}}\label{subsubsec: proof of Claim of alpha bound}
Define another sequence as $ \omega_k = \frac{16}{ 3 \alpha_k }- 9 k $. Then $\alpha_k = \frac{16}{3} \frac{1}{ 9k + \omega_k }$ and $\omega_1 = 7$, $\omega_2 = 38/5.$
We then derive the recurrence equation of $\omega_k$. 
According to \eqref{alpha recursion, 1st time}, we have 
\begin{align*}
 &  \frac{16}{3} \frac{1}{ 9k + 9 + \omega_{k+1} } 
   = \frac{ 2}{ 3 } \frac{1}{ 9k + \omega_k}  \frac{ 16 - 16/( 9k + \omega_k) }{ 2 + 16/(9k + \omega_k) }
    = \frac{2}{3} \frac{1}{9k + \omega_k} \frac{16( 9k + \omega_k - 1 ) }{ 2( 9k + \omega_k + 8 ) }  \\
    \Longrightarrow   &     9k + 9 + \omega_{k+1}  =  \frac{(9k + \omega_k)( 9k + \omega_k + 8  )}{  9k + \omega_k - 1  } \\
      \Longrightarrow   &  \omega_{k+1}  = \omega_k + \frac{1}{  9k + \omega_k - 1  }[ (9k + \omega_k)( 9k + \omega_k + 8  ) - (9k + \omega_k - 1)(9k+9 + \omega_k)  ] \\
        \Longrightarrow   &  \omega_{k+1}  = \omega_k + \frac{ 9 }{  9k + \omega_k - 1  } .
\end{align*}

It is easy to see that $\omega_k > 0 \Rightarrow \omega_{k+1} > \omega_k > 0 $, thus $$
\omega_k > \omega_1 = 7,  \; \forall \; k.$$
Furthermore, 
 $  \omega_{k+1}  = \omega_k + \frac{ 9 }{  9k + \omega_k - 1  }  \leq  \omega_k  + \frac{ 1 }{ k } ,$ thus
$$
 \omega_k \leq \omega_1 + \sum_{j=1}^{k-1} \frac{1 }{ j } \leq  8  + \log ( k-1).
$$
The lower bound and upper bound on $\omega_k $ imply upper and lower bounds on $\alpha_k$:
\begin{equation}\label{alpha precise upp low bound}
    \frac{16 }{3} \frac{1}{9k + 7} \geq \alpha_k \geq    \frac{16 }{3} \frac{1}{9k + 8 + \log(k-1) }. 
\end{equation}
As a side comment, this implies that $\lim_{ k \rightarrow \infty} \alpha_k = \frac{16}{27 k} \approx \frac{ 0.59}{k} .$
For our purpose, we need a universal lower bound on $\alpha_k$. 
When $k \geq 3$, we have $ 3k \geq 8 + \log(k-1)$, thus $ 12k \geq 9k + 8 + \log(k-1) $, which further implies
$$
 \frac{16 }{3} \frac{1}{9k + 8 + \log(k-1) } \geq \frac{  4}{9k },  \; \forall k \geq 3.
$$
Combining with the bound \eqref{alpha precise upp low bound}, we obtain
$$
\alpha_k \geq \frac{  4}{9k } > \frac{4}{ 9 (k+ 1)}, \; \forall \; k \geq 3.
$$
Notice that $  \alpha_1 = \frac{1}{3} >  \frac{4}{9} \cdot \frac{1}{2}, $ and $\alpha_2 = \frac{5}{24} > \frac{4}{9} \cdot \frac{1}{3} $, we have
$ \alpha_k > \frac{4}{9(k+1)}$ for any $k \geq 1$. This finishes the proof of the claim.


\subsection{Proof of Lemma \ref{lemma of weak AMGM ineq} }\label{subsec: proof of AMGM weak}
We rewrite the lemma statement below. 
	Suppose $P_i \in \dR^{N \times N}, i=1,\dots, n $ are projection matrices, then the lemma claims that
	\begin{equation}\label{weak matrix AM-GM ineq}
	\frac{1}{n!}  \sum_{\sigma = (\sigma_1 ,\dots, \sigma_n) \in \Gamma}    P_{\sigma_n} P_{\sigma_{n-1}} \dots  P_{\sigma_1}  \preceq   \frac{1}{n} \sum_i P_i. 
	\end{equation}

We first prove the case $n=2$, $n=3$ and $n = 4$, then prove the general case
$n = 2k$ and $n = 2k+1$ separately.

When $n = 2$,  \eqref{weak matrix AM-GM ineq}
reduces to $  P_1 P_2 + P_2 P_1 \preceq P_1  + P_2 $.
Notice that $ P_i = P_i^2$ since $P_i$ is a projection matrix, we have
$  P_1  + P_2 - P_1 P_2 + P_2 P_1 = P_1^2  + P_2^2 - P_1 P_2 + P_2 P_1
= (P_1 - P_2)^2 = (P_1 - P_2)^T (P_1 - P_2) \succeq 0 $.

When $n = 3$, \eqref{weak matrix AM-GM ineq}
reduces to $ \frac{1}{6} \sum_{i, j,k \text{ are distinct}} P_i P_j P_k \preceq \frac{1}{3} (P_1  + P_2 + P_3) $.
Note that $ (P_i - P_k) P_j (P_i- P_k) \succeq 0 $, thus
$$
P_i P_j P_i + P_k P_j P_k \succeq  P_i P_j P_k + P_k P_j P_i.
$$
Summing up the above inequality for all possible triples $(i,j,k)$, we get
\begin{equation}\label{intermediate bound n=3}
\sum_{i \neq j} P_i P_j P_i  \succeq     \sum_{i, j,k \text{ are distinct}} P_i P_j P_k.
\end{equation}
We then need to bound the left-hand-side of the above inequality. Since 
$ I - P_j \succeq 0 $, we have $ P_i (I - P_j) P_i \succeq 0 $, which implies
$
P_i \succeq P_i P_j P_i.
$
Summing up this inequality for all pairs $ i\neq j $, we obtain
$  \frac{1}{6} \sum_{i \neq j} P_i P_j P_i   \preceq \frac{1}{3} (P_1  + P_2 + P_3) $.
Combining with \eqref{intermediate bound n=3}, we obtain the desired inequality
$ \frac{1}{6} \sum_{i, j,k \text{ are distinct}} P_i P_j P_k \preceq \frac{1}{3} (P_1  + P_2 + P_3)$. 

The proof for $n=4$ illustrates partially the gist of a general proof, so we present this proof.
When $n = 4$, \eqref{weak matrix AM-GM ineq}
reduces to $ \frac{1}{24} \sum_{i, j,k,l \text{ are distinct}} P_i P_j P_k P_l \preceq \frac{1}{4} (P_1  + P_2 + P_3 + P_4) $.
Similar to \eqref{intermediate bound n=3} in the $n=3$ case, we first prove
\begin{equation}\label{n=4 intermediate bound}
\frac{1}{24} \sum_{i, j,k,l \text{ are distinct}} P_i P_j P_k P_l 
\leq  \frac{1}{12} \sum_{i \neq j } P_i P_j P_i. 
\end{equation}
To prove this inequality, we need the following two basic inequalities:
\begin{align*}
(P_i - P_l) (P_j + P_k)^2 (P_i - P_l) \succeq 0,   \\
(P_i + P_l) (P_j - P_k)^2 (P_i + P_l) \succeq 0. 
\end{align*}
Summing up these two inequalities, we can eliminate terms like $P_i P_j P_k P_i$ (with three distinct subscripts) and keep the terms like $ P_i P_j P_i $ (with two distinct subscripts) and $ P_i P_j P_k P_l $ (with four distinct subscripts), to obtain
$$
P_i P_j P_i + P_i P_k P_i + P_l P_j P_l + P_l P_k P_l 
\succeq P_i P_j P_k P_l + P_i P_k P_j  P_l +  P_l P_j P_k P_i + P_l P_k P_j P_i . 
$$
Summing up this inequality for all possible  $(i,j,k,l)$ that are distinct, we obtain
\eqref{n=4 intermediate bound}.
Similar to the proof of  $n=3$ case, we have
$   \frac{1}{12} \sum_{i \neq j } P_i P_j P_i \leq \frac{1}{4} (P_1 + P_2 + P_3 + P_4) $,
thus combining with   \eqref{n=4 intermediate bound} we obtain the desired result.

We next prove the case $n = 2k $, where $k \geq 2 $ is a positive integer.  
We will prove that
\begin{equation}\label{key relation odd case}
E_{\sigma  \in \Gamma}   ( P_{\sigma_n} P_{\sigma_{n-1}} \dots  P_{\sigma_1} )  \preceq
E_{\pi  \in \Gamma_k }   ( P_{\pi_1 } \dots
P_{\pi_{k-1}}   P_{\pi_k} P_{\pi_{k-1}}  \dots P_{\pi_1} ),
\end{equation}
where $\Gamma_k$ is the set of $k$-permutations of $1,2,\dots, n$ (here, a $k$-permutation is
a permutation of $k$ distinct numbers chosen from $1,2,\dots, n$),
and $E_{\sigma  \in \Gamma}$ and $E_{\pi \in \Gamma_k}$ denote
the expectation over a uniform  distribution on $\Gamma$ and $\Gamma_k$ respectively.

To prove \eqref{key relation odd case},
we need the following fact:
for any $\epsilon = (\epsilon_1, \dots, \epsilon_{k}) \in \{1, -1 \}^{k}$, we have
\begin{equation}\label{key observation PSD ineq}
G_{ \sigma, \epsilon } \triangleq ( P_{\sigma_n} + \epsilon_1 P_{\sigma_1} )  ( P_{\sigma_{n-1}} +\epsilon_2 P_{\sigma_2} )
\dots  ( P_{\sigma_{k+1}} +  \epsilon_{k} P_{\sigma_{k}} ) 
( P_{\sigma_{k+1}} +  \epsilon_{k} P_{\sigma_{k}} )  \dots ( P_{\sigma_n} +  \epsilon_{1}P_{\sigma_1} ) \succeq 0 .
\end{equation}
This relation holds because for any positive-semidefinite matrix $X $ and 
any symmetric matrix $ Y $, we have $Y X Y = Y^T X Y \succeq 0 $. Applying this fact $k$ times leads to \eqref{key observation PSD ineq}.

The expression of $G_{\sigma, \epsilon}$ in \eqref{key observation PSD ineq} involves $2^{k}$ terms
in the form of $P_{i_1} P_{i_2} \dots P_{i_n}$.
To prove \eqref{key relation odd case}, only two terms are of interest to us.
The strategy is to pick $\epsilon_i$'s properly so that summing up a bunch
of relations of the form \eqref{key observation PSD ineq} will eliminate all but the two
desired terms. We elaborate this strategy below. 

Define 
\begin{align*}
\Lambda_k \triangleq \{ (\epsilon_1 ,\dots, \epsilon_{k} ) \in \{ 1, -1 \}^k \mid 
\text{ the number of } -1 \text{ in } \epsilon_1, \dots, \epsilon_k  \text{ is odd}    \} , \\
\Lambda_k^c = \{ (\epsilon_1 ,\dots, \epsilon_{k} ) \in \{ 1, -1 \}^k \mid 
\text{ the number of } -1 \text{ in } \epsilon_1, \dots, \epsilon_k  \text{ is even}    \} .
\end{align*}
For example, when $k = 3$, $\Lambda_3 = \{ (-1,1,1), (1,-1,1), (1,1,-1), (-1,-1,-1) \}$,
and the complement $\Lambda_3^c = \{ (1,1,1) , (-1, -1, 1 ), (-1,1,-1), (1,-1,-1) \} $.
As a well-known fact,
\begin{equation}\label{Lambda has same size as Lambda c}
|\Lambda_k^c| -  |\Lambda_k| = 
\sum_{i \text{ is even } , 0 \leq i \leq n,  } { { n }\choose{ i } }
-  \sum_{i \text{ is odd }, 0 \leq i \leq n } { { n }\choose{ i } }
= (1 - 1)^k = 0,
\end{equation}


This matrix $ G_{ \sigma, \epsilon } $ can be expressed as the sum
of $2^{k}$ terms, and each term is of the form $ \pm P_{\pi_1} \dots P_{\pi_n}$,
where $\pi_i \in \{ \sigma_i, \sigma_{n + 1 - i} \} $. 
For the fixed permutation $\sigma$, define a set
$$
\ \Omega(\sigma) = \{ (\pi_1, \dots, \pi_n) \mid \pi_i \in \{ \sigma_i, \sigma_{n + 1 - i} \} , \forall i   \}
$$
We partition the set into three subsets:
\begin{align*}
\Omega_0(\sigma) & =  \{ (\pi_1, \dots, \pi_n) \in \Omega \mid \pi_i = \pi_{n+1 -i} , \forall i   \},  \\
\Omega_1(\sigma) & = 
\{ (\pi_1, \dots, \pi_n) \in \Omega \mid \pi_i \neq \pi_{n+1 -i} , \forall i  \}, \\
\Omega_2(\sigma) & = \Omega \backslash (\Omega_0 \cap \Omega_1). 
\end{align*}
For most of the proof, we will use the abbreviation $\Omega_t = \Omega_t(\sigma), t=0,1,2.$
For any $\pi = (\pi_1, \dots, \pi_n) \in \Omega$, define an indicator vector of $\pi$ as
$
\delta(\pi) = (\delta_1 , \dots, \delta_k),
$ where each $\delta_i$ is determined by
\begin{equation} \label{def of indicator delta}
\delta_i = \mathbb{I}( \pi_i - \pi_{n+1-i} ) = \begin{cases}
0,  &  \pi_i = \pi_{n+1 - i}, \\
1,   &  \pi_i \neq \pi_{n+1 - i}, 
\end{cases}
\end{equation}
where $\mathbb{I}(z)$ equals $0$ if $z = 0$ and equals $1$ if $z \neq 0$,
For example, when $n= 6$ and $\pi = ( \sigma_1 , \sigma_2, \sigma_3, \sigma_4, \sigma_2, \sigma_6  )$, the corresponding indicator vector is $(0 ,0 ,1)  $; when $\pi = ( \sigma_1 , \sigma_2, \sigma_3, \sigma_3, \sigma_2, \sigma_1  )$, the indicator vector is $(0, 0, 0)$.
Clearly, we have 
\begin{equation}\label{delta Omega relation}
\delta (\pi) = (0,0,\dots, 0) , \; \forall \pi \in \Omega_0; \quad 
\delta(\pi) = (1,1,\dots ,1),  \; \forall \pi \in \Omega_1; \quad \delta(\pi) \notin \{ 0_k, 1_k \}, \; \forall \pi \in \Omega_2. 
\end{equation}

In the expression of $G_{\sigma, \epsilon}$, half of the terms have coefficient $1$
and the other half have coefficient $-1$.
To understand which terms have coefficient $1$ and which have coefficient $-1$,
consider a special $\epsilon = (-1, 1, \dots, 1)$, i.e., $\epsilon_1 = -1 $
and all other $\epsilon_i = 1$.
A term with coefficient $-1$
has the form $ P_{\sigma_1} P_{\pi_2} \dots P_{\pi_{n-1}} P_{\sigma_n}  $
or $ P_{\sigma_n} P_{\pi_2} \dots P_{\pi_{n-1}} P_{\sigma_1} $, i.e., with an indicator vector
whose first element $\delta_1 = 1 $,
and a term with coefficient $1$ has the form
$ P_{\sigma_1} P_{\pi_{n-1} } \dots P_{\pi_{2}} P_{\sigma_1}  $
or $ P_{\sigma_n} P_{\pi_{n-1} } \dots P_{\pi_{2}} P_{\sigma_n} $, i.e., with an indicator vector whose first element $\delta_1 = 0$.
We can see that the coefficient is in fact $ \epsilon_1^ {\delta_1} $. 
For general $\epsilon \in \Lambda$ and $\pi \in \Omega$, the coefficient of
$ P_{\pi_n} P_{\pi_{n-1} } \dots P_{\pi_{2}} P_{\pi_1 } $ in $G_{\sigma, \epsilon}$
is $ (\epsilon_1)^{ \delta_1 } \dots   (\epsilon_k)^{ \delta_k }$,
where $\delta = \delta(\pi)$ is defined as in \eqref{def of indicator delta}. We can then write the expression of $G_{\sigma, \epsilon}$ as
$$
G_{ \sigma, \epsilon} = \sum_{\pi \in \Omega} \epsilon_1^{\delta_1} \dots \epsilon_k^{\delta_k}  P_{\pi_n} P_{\pi_{n-1}} \dots  P_{\pi_1} . 
$$
Summing up this relation for all $\epsilon$ in $\Lambda_k$, we have
\begin{equation}\label{G delta sum}
\sum_{\epsilon \in \Lambda_k } G_{ \sigma, \epsilon }
= \sum_{ \epsilon \in \Lambda_k } \sum_{\pi \in \Omega} \epsilon_1^{\delta_1} \dots \epsilon_k^{\delta_k} P_{\pi_n} P_{\pi_{n-1}} \dots  P_{\pi_1}
= \sum_{\pi \in \Omega}P_{\pi_n} P_{\pi_{n-1}} \dots  P_{\pi_1} \left( \sum_{ \epsilon \in \Lambda_k } \epsilon_1^{\delta_1} \dots \epsilon_k^{\delta_k}  \right). 
\end{equation}
Note that in this expression, $\delta_1, \dots, \delta_k $ depend on $\pi$. 

Denote $0_k = (0,0, \dots, 0) \in \dR^k $, $1_k = (1,\dots, 1) \in \dR^k$.
Define 
\begin{align*}
g_k (\delta) \triangleq \sum_{ \epsilon \in \Lambda_k } \epsilon_1^{\delta_1} \dots \epsilon_k^{\delta_k} ,  \quad
h_k (\delta) \triangleq \sum_{ \epsilon \in \Lambda_k^c } \epsilon_1^{\delta_1} \dots \epsilon_k^{\delta_k}.
\end{align*}
For any $\delta \neq 0_k$, we have
$ g_k(\delta) + h_k(\delta) = \sum_{ \epsilon \in \{1,-1 \}^k } \epsilon_1^{\delta_1} \dots \epsilon_k^{\delta_k}  = ( 1^{\delta_1} + (-1)^{\delta_1} ) \dots ( 1^{\delta_k} + (-1)^{\delta_k} ) = 0  ,$
thus 
\begin{equation}\label{h is negative to g}
h_k(\delta) = -g_k(\delta), \quad \forall \delta \neq 0_k .
\end{equation}
It is easy to see that
\begin{equation}\label{g value for 0 and 1}
\frac{1}{|\Lambda_k|} g_k(\delta) =  \begin{cases}
1  , & \delta = (0,0,\dots, 0),  \\
-1,  &  \delta = (1,1, \dots, 1) . \\  
\end{cases}
\end{equation}

We will prove: for any $\delta \notin \{ 0_k, 1_k\}, $
\begin{equation}\label{g equal zero for any delta}
g_k (\delta) = \sum_{ \epsilon \in \Lambda_k } \epsilon_1^{\delta_1} \dots \epsilon_k^{\delta_k} = 0 , 
\end{equation}
We prove \eqref{g equal zero for any delta} by induction on $k$.
When $k = 2$, $\Lambda_2 = \{ (-1,1), (1,-1) \}$, we have:
\begin{align*}
\text{when } \delta = (0,1), & \quad  g_2(\delta) = (-1)^0 1^1 + 1^0 (-1)^1 = 1 -1 = 0,  \\
\text{when } \delta = (1,0), & \quad  g_2(\delta) = (-1)^1 1^0 + 1^1 (-1)^0 = -1 + 1 = 0. 
\end{align*}
Assume \eqref{g equal zero for any delta} holds for $k-1$, i.e., 
\begin{equation}\label{induction hypo g}
g_{k-1}( \hat{\delta} ) = 0, \;
\forall \; \hat{\delta}  \in \{0,1 \}^{k-1} \backslash \{ 0_{k-1}, 1_{k-1} \} .
\end{equation}
According to \eqref{h is negative to g}, we have
\begin{equation}\label{induction hypo h}
h_{k-1} ( \hat{\delta}  ) = 0, \forall \; \hat{\delta}  \in \{0,1 \}^{k-1} \backslash \{ 0_{k-1}, 1_{k-1} \} .
\end{equation}
Now consider $k$. Since $\delta \neq 0_k$, there must exist some $j$ such that $\delta_j =1$; without loss of generality, we assume 
\begin{equation}\label{delta k =1} 
\delta_k = 1. 
\end{equation}
Partition $\Gamma_k$ into two sets: 
\begin{equation}\label{epsilon k 1 or -1}
\Lambda_{k,1}  \triangleq \{ \epsilon \in \Lambda_k \mid \epsilon_k =1   \}, \quad 
\Lambda_{k,2}  \triangleq \{ \epsilon \in \Lambda_k \mid \epsilon_k =-1   \} .
\end{equation}
If $\epsilon$ contains an odd number of $-1$ and the last element $\epsilon_k = 1$ (or $\epsilon_k = -1$),
then the first $k-1$ elements contain an odd (or even) number of $-1$.
Thus
$$
\Lambda_{k,1} = \{ ( \hat{\epsilon}, 1) \mid \hat{\epsilon}\in \Lambda_{k-1}  \}, \quad
\Lambda_{k,2} = \{ (\hat{\epsilon}, -1) \mid \hat{\epsilon} \in \Lambda_{k-1}^c  \}.
$$
Split $g_k(\delta)$ into two parts 
$ g_k (\delta) = g_{k,1}(\delta) + g_{k,2}(\delta),
$
where
$$
g_{k,1}(\delta) = \sum_{ \epsilon \in \Lambda_{k,1} } \epsilon_1^{\delta_1} \dots \epsilon_k^{\delta_k},  \quad g_{k,2}(\delta) = \sum_{ \epsilon \in \Lambda_{k,2} } \epsilon_1^{\delta_1} \dots \epsilon_k^{\delta_k}.
$$
Denote $\hat{\delta} = (\delta_1, \dots, \delta_{k-1})$. 
We already assume $\delta \neq 1_k$ and $\delta_k = 1$, so we know 
\begin{equation}\label{induction not all 1}
\hat{\delta} \neq 1_{k-1}. 
\end{equation}
But it is possible that $\hat{\delta} = 0_{k-1}$. 
Consider two cases. 

\textbf{Case 1}: $\hat{\delta} = 0_{k-1}$, i.e.,  $\delta = (0_{k-1} , 1) $. 

In this case 
\begin{align*}
g_{k,1}(\delta) = \sum_{ \epsilon \in \Lambda_{k,1} } \epsilon_1^{\delta_1} \dots \epsilon_k^{\delta_k} = \sum_{ \epsilon \in \Lambda_{k,1} } \epsilon_1^{0} \dots
\epsilon_{k-1}^{0} \epsilon_k^{1} 
=   \sum_{ \epsilon \in \Lambda_{k,1} } \epsilon_k^{1}  
=   \sum_{ \epsilon \in \Lambda_{k,1} } 1^1
= |\Lambda_{k,1} | = |\Lambda_{k-1}|,   \\
g_{k,2}(\delta) = \sum_{ \epsilon \in \Lambda_{k,2} } \epsilon_1^{\delta_1} \dots \epsilon_k^{\delta_k} = \sum_{ \epsilon \in \Lambda_{k,2} } \epsilon_1^{0} \dots
\epsilon_{k-1}^{0} \epsilon_k^{1} 
=   \sum_{ \epsilon \in \Lambda_{k,2} } \epsilon_k^{1}  
=   \sum_{ \epsilon \in \Lambda_{k,2} } (-1)^1
= - |\Lambda_{k,2} | = - |\Lambda_{k-1}^c|, 
\end{align*}
Thus 
$$
g_k(\delta) =  g_{k,1}(\delta)  +  g_{k,2}(\delta) = |\Lambda_{k-1}| - |\Lambda_{k-1}^c| = 0,
$$
where the last step is due to \eqref{Lambda has same size as Lambda c}.

\textbf{Case 2}: $\hat{\delta} \neq 0_{k-1}$. Together with \eqref{induction not all 1}, we have
$$ 
\hat{\delta}  \notin \{0_{k-1}, 1_{k-1} \}.
$$   
which enables us to apply the induction hypothesis \eqref{induction hypo g} and its corollary \eqref{induction hypo h}. In fact,
\begin{align*}
g_{k,1}(\delta) = \sum_{ \epsilon \in \Lambda_{k,1} } \epsilon_1^{\delta_1} \dots \epsilon_k^{\delta_k} 
\overset{\eqref{delta k =1},\eqref{epsilon k 1 or -1}}{=}  \sum_{ \epsilon \in \Lambda_{k,1} } \epsilon_1^{\delta_1} \dots
\epsilon_{k-1}^{\delta_{k-1}} 1^1 
=  \sum_{ \hat{ \epsilon } \in \Lambda_{k-1} } 
\hat{ \epsilon }_1^{\delta_1} \dots
\hat{ \epsilon }_{k-1}^{\delta_{k-1}}
= g_{k-1}( \hat{\delta} )  \overset{\eqref{induction hypo g}}{=} 0,    \\
g_{k,2}(\delta) = \sum_{ \epsilon \in \Lambda_{k,2} } \epsilon_1^{\delta_1} \dots \epsilon_k^{\delta_k} 
\overset{\eqref{delta k =1},\eqref{epsilon k 1 or -1}}{=} \sum_{ \epsilon \in \Lambda_{k,2} } \epsilon_1^{\delta_1} \dots
\epsilon_{k-1}^{\delta_{k-1}} (-1)^{1} 
= - \sum_{ \hat{ \epsilon } \in \Lambda_{k-1}^c } 
\hat{ \epsilon }_1^{\delta_1} \dots
\hat{ \epsilon }_{k-1}^{\delta_{k-1}}
= h_{k-1}( \hat{\delta} )  \overset{\eqref{induction hypo h}}{=} 0.
\end{align*}
Thus $g_k(\delta) = g_{k,1}(\delta)  +  g_{k,2}(\delta) = 0 $.

In both cases, we have proved $g_k(\delta) = 0$, which finishes the induction step.
Therefore \eqref{g equal zero for any delta} holds for any $k$.

Next, we analyze the sum $\sum_{\epsilon \in \Lambda_k } G_{ \sigma, \epsilon }.$
According to \eqref{G delta sum}, we have
\begin{align*}
\sum_{\epsilon \in \Lambda_k } G_{ \sigma, \epsilon }
& = \sum_{\pi \in \Omega}P_{\pi_n} P_{\pi_{n-1}} \dots  P_{\pi_1} \left( \sum_{ \epsilon \in \Lambda_k } \epsilon_1^{\delta_1} \dots \epsilon_k^{\delta_k}  \right) \\
& = \sum_{\pi \in \Omega}P_{\pi_n} P_{\pi_{n-1}} \dots  P_{\pi_1} g_k(\delta(\pi) )  \\
& \overset{(i)}{=} \sum_{\pi \in \Omega_0 }P_{\pi_n} P_{\pi_{n-1}} \dots  P_{\pi_1}  \cdot g_k(0_k)
+ \sum_{\pi \in \Omega_1 }P_{\pi_n} P_{\pi_{n-1}} \dots  P_{\pi_1}  \cdot g_k(1_k)
+ \sum_{\pi \in \Omega_2 }P_{\pi_n} P_{\pi_{n-1}} \dots  P_{\pi_1}  \cdot g_k(\delta(\pi) ) \\
& \overset{ (ii) }{=} \sum_{\pi \in \Omega_0 }P_{\pi_n} P_{\pi_{n-1}} \dots  P_{\pi_1}  \cdot |\Gamma_k|
+ \sum_{\pi \in \Omega_1 }P_{\pi_n} P_{\pi_{n-1}} \dots  P_{\pi_1}  \cdot (-1)|\Gamma_k|
+ \sum_{\pi \in \Omega_2 }P_{\pi_n} P_{\pi_{n-1}} \dots  P_{\pi_1}  \cdot 0 \\
& = |\Gamma_k|  \left( \sum_{\pi \in \Omega_0 }P_{\pi_n} P_{\pi_{n-1}} \dots  P_{\pi_1}  
- \sum_{\pi \in \Omega_1 }P_{\pi_n} P_{\pi_{n-1}} \dots  P_{\pi_1} \right). 
\end{align*}
where (i) is due to \eqref{delta Omega relation} and (ii) is due to \eqref{g value for 0 and 1}, \eqref{g equal zero for any delta}.
According to \eqref{key observation PSD ineq}, any $G_{\sigma, \epsilon} \succeq 0$,
thus the above relation implies the following important relation
\begin{equation}\label{key inequality for weak AM GM}
\sum_{\pi \in \Omega_0 }P_{\pi_n} P_{\pi_{n-1}} \dots  P_{\pi_1}  
\succeq \sum_{\pi \in \Omega_1 }P_{\pi_n} P_{\pi_{n-1}} \dots  P_{\pi_1}
\end{equation}
Note that this relation holds for a fixed permutation $\sigma $ and the corresponding set
$\Omega_0 = \Omega(\sigma)$ and $\Omega_1(\sigma)$. 
Each $\pi \in \Omega_0 $ corresponds to
a $k$-permutation $\chi$ of $ (12\dots n) $ determined by  $\pi = (\chi_1 \dots \chi_{k-1} \chi_k \chi_k \chi_{k-1} \dots \chi_1)$
and each  $\pi \in \Omega_1$ corresponds to a permutation of $(12\dots n)$.
We rewrite \eqref{key inequality for weak AM GM} as
$$
\sum_{\pi \in \Omega_0(\sigma) }P_{\pi_n} P_{\pi_{n-1}} \dots  P_{\pi_1}  
\succeq \sum_{\pi \in \Omega_1(\sigma) }P_{\pi_n} P_{\pi_{n-1}} \dots  P_{\pi_1}
$$
and summing up this relation for all possible permutations $\sigma \in \Gamma$ leads to
$$
E_{\chi  \in \Gamma_k }   ( P_{\chi_1 } \dots
P_{\chi_{k-1}}   P_{\chi_{k}}  P_{\chi_{k}} P_{\chi_{k-1}}  \dots P_{\chi_1} ) \succeq
E_{\sigma  \in \Gamma}   ( P_{\sigma_n} P_{\sigma_{n-1}} \dots  P_{\sigma_1} ) ,
$$
which is exactly \eqref{key relation odd case}. 

It remains to prove
\begin{equation}\label{relax square to sum Pi}
E_{\chi  \in \Gamma_k }   ( P_{\chi_1 } \dots
P_{\chi_{k-1}}   P_{\chi_k} P_{\chi_{k-1}}  \dots P_{\pi_1} ) \preceq \frac{1}{n}\sum_i P_i.
\end{equation}
In fact, for any positive-semidefinite matrix $X $ and 
any symmetric matrix $ Y $, we have $Y X Y = Y^T X Y \succeq 0 $. Applying this fact $k-1$ times leads to \eqref{relax square to sum Pi}.

Combining \eqref{key relation odd case} and \eqref{relax square to sum Pi}, we immediately obtain the desired result \eqref{weak matrix AM-GM ineq} for the case $n = 2k$.

The case that $n = 2k-1$ is an odd number is almost the same, except that the key quantity $G_{ \sigma, \epsilon }  $ is now defined as
\begin{equation}\label{key observation PSD ineq, odd case}
G_{ \sigma, \epsilon } \triangleq ( P_{\sigma_n} + \epsilon_1 P_{\sigma_1} )  
\dots  ( P_{\sigma_{k+1}} +  \epsilon_{k-1} P_{\sigma_{k-1}} )
P_{\sigma_{k}} 
( P_{\sigma_{k+1}} +   \epsilon_{k-1} P_{\sigma_{k-1}} )  \dots ( P_{\sigma_n} +  \epsilon_{1}P_{\sigma_1} ) .
\end{equation}
In words, we pair $P_{\sigma_i}$ with $P_{\sigma_{n+1-i}}$ for $i=1,\dots, k-1$ and leave
$P_{\sigma_k}$ alone (following the same rule it would have been paired with itself).
The rest of the proof is almost the same as the even case, so we skip it. 
$\quad \quad $ \textbf{Q.E.D.}

{\black
\section{ Numerical Experiments }\label{sec: simulation}
In this section, we test the performance of cyclic ADMM and RP-ADMM for solving various kinds of linear systems.
As a benchmark, we also test the gradient descent method (GD) with a constant stepsize $\alpha = 1/\lambda_{\max}(A'A)$
for solving the least square problem $\min_{x\in \dR^N} \| Ax-b \|^2/2$.
Of course there are many other advanced algorithms for solving the least square problem such as the conjugate gradient method,
but we do not consider them since our focus is on testing the two ADMM algorithms.  
These two ADMM algorithms can be used to solve far more general problems than just linear systems,
and we believe that the performance comparison for solving linear systems can shed light on more general scenarios.

In the numerical experiments, we set $b = 0$, thus the unique optimal solution is $x^* = 0$.
 The coefficient matrix $A$ will be generated according to one of the random distributions below:
\begin{itemize}
  \item Gauss: independent Gaussian entries $A_{i,j} \sim \mathcal{N}(0,1)$.
  \item Log-normal: independent log-normal entries $A_{i,j} \sim \text{exp}( \mathcal{N} (0,1) )$. 
  \item Uniform: each entry is drawn independently from a uniform distribution on $[0,1]$.
  \item Circulant Hankel: circulant Hankel matrix with independent standard Gaussian entries. More specifically,
  generate $\delta_1, \delta_2, \dots, \delta_{N} \sim \mathcal{N}(0,1) $ and let $A_{i,j} = \delta_{i+j-1}$ (define $\delta_{k} = \delta_{k-N}$ if $k>N$).
  Note that the entries of the circulant Hankel matrix are not independent since one $\delta_i$ can appear in multiple positions.


\end{itemize}


For the two ADMM algorithms, we only consider the $n$-coordinate versions, i.e. each block consists of only one coordinate.
We let the three tested algorithms start from the same random initial point $y^0 = [x^0; \lambda^0] $ (GD will start from $x^0$).
To measure the performance, we define  the epoch complexity $k$ to be the minimum $k$ so that the relative error
$$ \| A x^k - b \| / \| A x^0 - b \| < \epsilon ,$$
where $\epsilon$ is a desired accuracy (we consider $10^{-2}$ and $10^{-3}$\footnote{For high accuracy such
	as $\epsilon = 10^{-6}$, it takes too many epochs for the algorithms to converge when $n=100$ as most matrices
	we generated are highly ill-conditioned, so we do not report the results. Based on the limited experiments for high accuracy, similar gaps between RP-ADMM and GD are observed. }). 
For the two ADMM algorithms, one epoch refers to one round of primal and dual steps; for GD, one epoch refers to one gradient step.
The total computation time should be proportional to the epoch complexity since  GD and the two ADMM variants have similar per-epoch cost\footnote{In matlab simulations each epoch of GD takes much less time than a round of ADMM because matlab implements matrix operations much faster than a ``for'' loop. For a more fair CPU time comparison, one should use other programming languages such as C. }: a gradient descent step $x^{k+1} = x^k - \alpha A^T(Ax-b)$ contains two matrix-vector multiplications and thus takes time $2N^2 + O(N)$, and an ADMM round also takes time $2N^2 + O(N)$ (the primal update step of ADMM takes time $2N^2 + O(N)$ and the dual update step of ADMM takes time $O(N)$).
We test 1000 random instances for $N \in \{3,10\}$ and $300$ random instances for $N=100$, and record the geometric mean of the number of epochs.
In the table, ``Diverg. Ratio'' represents the percentage of tested instances for which cyclic ADMM diverges and
``CycADMM'' represents ``cyclic ADMM'' (note that RP-ADMM converges in all instances we tested, so 
its divergence ratio is 0).
Note that for cyclic ADMM we only report the epoch complexity when it converges, while for RPADMM and GD we report the 
epoch complexity in all tested instances.
If restricting to the successful instances of cyclic ADMM, we find that the epoch complexity of RPADMM does not change too much, while the epoch complexity of GD will be reduced (significantly in some settings). 

The simulation results are summarized in Table \ref{table simulation}.
The main observations from the simulation are:
\begin{itemize}
\item For all random distributions of $A$ we tested, cyclic ADMM does not always converge even when $N$ is fixed to be $3$. For $N=100$ and many random distributions, cyclic ADMM diverges with probability $1$.
This means that the divergence of cyclic ADMM is not merely a ``worst-case'' phenomenon, but actually quite common.
When the dimension increases, the divergence ratio will increase.

\item For standard Gaussian entries, cyclic ADMM converges with high probability.
When cyclic ADMM converges, it converges faster than RP-ADMM and sometimes much faster.

\item RPADMM typically converges faster than the basic gradient descent method and sometimes more than $10$ times faster.
\end{itemize}




\begin{table}[!htbp]
\caption{ \emph{ Results of Solving Linear Systems by Cyclic ADMM, RP-ADMM and GD. For the two ADMM variants, one epoch refers to one round of primal and dual steps; for GD, one epoch refers to one gradient step.} }
\label{table simulation}
\centering
\begin{tabular}{|c|c|c|c|c|c|c|c|}
 \hline
 \multirow{2}{*}{ N } &  \multirow{2}{*}{ Diverg. Ratio  } &
 \multicolumn{3}{c|}{ Epochs for $\epsilon = 0.01 $  } &
 \multicolumn{3}{c|}{ Epochs for $\epsilon = 0.001 $  }    \\  
 \cline{3-8}  
   & & CycADMM\footnotemark & RPADMM & GD & CycADMM  & RPADMM & GD  \\ 
 \hline
 \hline
   \multicolumn{8}{|c|}{ Gaussian }  \\
    \hline
  3 & 0.7\%   & 1.4e01 & 3.4e01 & 5.0e01 & 3.2e01  & 8.8e01  & 1.4e02 \\
     \hline
  10 & 1.1\%   & 4.1e01 & 1.8e02  & 2.0e02 & 1.2e02  & 1.1e03  & 1.5e03 \\
       \hline
  100 & 3\%  & 1.7e02 & 4.3e02  & 3.6e02 & 1.0e03  & 7.4e03  & 6.5e03 \\
       \hline \hline
       \multicolumn{8}{|c|}{ Log-normal  }  \\
        \hline
 3 & 0.8\%   & 1.5e01 & 3.7e01 & 5.7e01 &  3.3e01  & 9.6e01  & 1.7e02 \\
     \hline
       10 & 39.2\%   & 1.2e02 & 3.4e02 & 6.4e02 &   3.2e02  & 2.4e03  & 6.3e03 \\
       \hline
  100 & 100\%  & N/A & 5.5e02 & 5.4e03  &  N/A  & 8.8e03  & 1.0e05 \\
       \hline         \hline
        \multicolumn{8}{|c|}{ Uniform }  \\
         \hline
 3 & 3.2\%   & 2.8e01 & 7.4e01 & 1.5e02 & 7.0e01  & 2.6e02  & 6.0e02 \\
     \hline
       10 & 83.0\%  & 2.1e02  & 4.1e02 & 1.2e03   & 5.2e02 & 3.0e03 & 9.1e03  \\
       \hline
  100 & 100\%  & N/A &  9.1e02 & 1.4e04  &   N/A   & 1.4e04  & 9.7e04 \\
       \hline   \hline
       \multicolumn{8}{|c|}{ Circulant Hankel  }  \\
        \hline
 3 & 5.6\%   & 1.2e01 & 1.7e01 & 1.5e01  & 1.7e01  & 2.8e01  & 2.6e01 \\
     \hline
       10 & 54.3\%   & 4.2e01 & 6.0e01  & 6.5e01 & 7.5e01  & 1.3e02 & 1.7e02 \\
       \hline
  100 & 100\%  & N/A & 1.3e02 & 1.7e02  & N/A  & 2.9e02 & 6.5e02 \\
       \hline
 \end{tabular}
 \end{table}
\footnotetext{For cyclic ADMM, only record the iteration complexity in convergent instances.}

We have also tested BR-ADMM for solving the same problems, though the simulation results are not listed in the above table. As expected, BR-ADMM also always converges for solving these linear systems. The convergence speed is usually slower than RP-ADMM. Nevertheless, BR-ADMM can save some sampling time compared to RP-ADMM, and may be more favorable if random permutation
is not available due to system architecture constraint. The detailed comparison of BR-ADMM and RP-ADMM, and the design of other randomized schemes or even deterministic schemes that
outperform RP schemes are left as future work.  

}

\section{ Concluding Remarks }
In this paper, we prove the expected convergence of randomly permuted ADMM (RP-ADMM) for solving a non-singular square system of equations (extension to non-square systems is straightforward). We also prove a bound on the expected convergence rate of RP-ADMM for solving linear systems and the expected convergence rate of RP-BCD for solving quadratic problems. 
The motivation is to resolve the divergence issue of cyclic multi-block ADMM.
 Our result shows that RP-ADMM may serve as a simple remedy, and we expect RP-ADMM to be one of the important solvers in large-scale optimization.
 One interesting finding along the path is that the update matrix of RP-BCD has spectrum lying in $(-1/3, 1)$ instead of the commonly seen $(-1,1)$.


 Randomly permutation is widely known to be empirically better than independently randomized versions, but little was known about its theoretical properties in general. 
Note that most existing analyses of BCD (e.g. \cite{tseng2001convergence,beck2013convergence,sun2015improved}) are applicable to both the cyclic update rule and the random permutation update rule. However, in light of a recent study which
established an up to $O(n^2)$ gap between cyclic CD and R-CD \cite{sun2016worst}, it is unlikely that RP-CD will have the same rate as cyclic CD. Our result in this paper established, for the first time, an $O(n)$ gap between RP-CD and cyclic-CD for general quadratic problems, making some progress towards the conjecture that RP-CD is faster than R-CD.


We emphasize that the convergence speed analysis of large-scale optimization has mostly been limited to independently randomized update order in the past decade. Going beyond independent randomized order is an important topic for enlarging the scope of large-scale optimization.  Not only the analysis of random permutation is quite challenging, even the analysis of the most classical cyclic order is highly nontrivial \cite{sun2016worst}. There are quite a few open questions regarding the convergence rate of non-independent-randomized order. 
Regarding the random permutation order, a very interesting open question is the worst-case convergence rate of RP-BCD for quadratic problems. Due to the close relation with matrix AM-GM inequality, this problem seems to be a quite fundamental problem. Moving to ADMM, the similar questions about the convergence rate of various variants of ADMM, including RP-ADMM and BR-ADMM, are also open.  

\section{ Acknowledgment}
We thank an anonymous reviewer for many helpful comments on the manuscript, which enabled us to improve the presentation of the paper.

\appendix

\bibliographystyle{unsrt} 

\bibliography{ref_test}

\end{document}